\documentclass[12pt]{article}             
\usepackage{graphicx} 
\usepackage{subcaption}
\usepackage{srcltx,epsfig}
\usepackage{amsmath, amssymb, amsthm}
\usepackage{color}
\usepackage{float}
\usepackage{xcolor}
\usepackage{multirow}
\usepackage[margin=1in]{geometry}
\usepackage{verbatim}
\usepackage{pgfplots}
\pgfplotsset{compat=1.18} 
\usepackage{mathrsfs}
\usepackage{array}
\usepackage{booktabs} 
\usepackage{placeins} 
\newcommand\pd[2]{\dfrac{\partial {#1}}{\partial {#2}}}

\newcommand{\mN}{\mathcal{N}}
\newcommand{\Kn}{\text{Kn}}
\newtheorem{remark}{Remark}[section]

\usepackage{bm}
\usepackage[colorlinks,linkcolor=blue,anchorcolor=blue,citecolor=red]{hyperref}
\usepackage[T1]{fontenc}
\usepackage{tikz}
\usetikzlibrary{shapes.geometric, arrows}
\usetikzlibrary{matrix}
\usepackage{algorithm}
\usepackage{algpseudocode}
\usepackage[numbers]{natbib}
\usepackage{makecell}
\usepackage{siunitx}
\usepackage{chngcntr}
\usepackage{diagbox}
\counterwithin{figure}{section}
\usepackage{indentfirst}
\numberwithin{equation}{section} 

\newcommand{\bmx}{\bm{x}}

\newcommand{\dd}{\mathrm{d}}
\newcommand{\mM}{\mathcal{M}}
\newcommand{\mE}{\mathcal{E}}

\newcommand{\bu}{\bm{u}}
\newcommand{\bv}{\bm{v}}
\newcommand{\bbR}{\mathbb{R}}

\graphicspath{{images/}}

\title{Learning a generalized Navier–Stokes model beyond the continuum regime}
\author{
   Sihong Shao\thanks{CAPT, LMAM, and School of Mathematical Sciences, Peking University, Beijing, 100871, China (sihong@math.pku.edu.cn).} \and 
    Yanli Wang\thanks{Beijing Computational Science Research Center, Beijing, 100193, China (ylwang@csrc.ac.cn).} \and
    Zhongwei Xu\thanks{School of Mathematical Sciences, Peking University, Beijing, 100871, China (xuzw@pku.edu.cn).}
    }

\begin{document}
\date{}
\maketitle

\begin{abstract}
Navier-Stokes (NS) equations lose accuracy in the transition regime when the Knudsen number is large, while the Boltzmann equation provides an adequate kinetic description with a substantially higher computational cost. In this work, a data-driven generalized Navier-Stokes (GNS) model is proposed to extend the applicability of the NS equations to the transition regime. A new relationship between the non-equilibrium variables, such as the stress tensor and the heat flux, and the equilibrium variables, including the density, macroscopic velocity, and the temperature, together with the Knudsen number, is learned from the data formed by the numerical solution to the Boltzmann equation. The standard end-to-end learning is utilized to construct the loss function, based on which, a long-term end-to-end learning method is proposed for the loss function to reduce accumulated error and improve long-term predictive accuracy. Several numerical examples, including one-dimensional wave and Riemann problems as well as two-dimensional isentropic vortex and Taylor-Green vortex problems, are studied to validate the efficiency of this new GNS model.

\end{abstract}

\section{Introduction}
Accurately simulating gas flows across different regimes remains a major computational challenge. Navier-Stokes equations are the most important governing equations in the continuum regime. However, they fail in the transition regime, not to mention the rarefied regime. 

These flow regimes are classified by the Knudsen number (Kn), which is defined as the ratio of the gas molecular mean free path to a characteristic length scale. In the continuum regime (Kn $\lesssim 0.01$), the macroscopic Navier-Stokes (NS) equations provide an accurate and computationally efficient model of the flow dynamics. As Kn increases, this classical macroscopic model becomes unreliable in the transition regime ($0.01 \lesssim$ Kn $\lesssim 10$), requiring a more fundamental kinetic description. The Boltzmann equation describes the evolution of the particle distribution function in phase space by accounting for molecular transport and collisions. However, the numerical simulation of the Boltzmann equation is generally considered computationally expensive due to its high dimensionality and the computational complexity of the quadratic collision term.

To reduce the computational cost associated with the high dimensionality of Boltzmann equations in the transition regime, one of the classical approaches is to deduce the macroscopic governing equations beyond the continuum regime. For example, in the moment equation, a series of macroscopic equations are derived by expanding the microscopic distribution function by a finite set of velocity moments, which can significantly reduce the computational cost. A pioneering example is Grad’s moment method \cite{grad1949kinetic}, in which the closure of the moment system is achieved by truncating the Hermite expansion. This simple closure leads to the loss of hyperbolicity and inaccurate simulations. To avoid this, several recent works have proposed closures that preserve the hyperbolicity of the moment system. In \cite{cai2014globally}, a globally hyperbolic regularization is constructed by modifying the governing equations of the highest-order moments with the characteristic speed correction. In \cite{koellermeier2014framework, koellermeier2014hyperbolic}, global hyperbolicity is guaranteed by employing quadrature-based projection methods. Besides, the most famous moment equations are the regularized 13 models, for which several different versions are derived \cite{ManuelR13, Cai2020} and have been validated in the high Mach flow problem \cite{Cai2020}. 

For the Navier-Stokes type models, extended macroscopic models are also derived by deducing new closure relationships for the non-equilibrium variables. For example, the coupled constitutive relations proposed in \cite{rana2018coupled} introduce coupling between the stress and heat flux to describe rarefaction effects not captured by the classical NS equations. More recently, nonlinear coupled constitutive relations have been applied to hypersonic thermochemical nonequilibrium flows, improving predictions of surface heat flux over the NS model \cite{zeng2023nonlinear}. Besides, the Chapman-Enskog method \cite{Chapman1970} is utilized to deduce Burnett and Super-Burnett models \cite{Bobylev1982}, which are not linearly stable, and limit their application in the simulation of fluid dynamics. Data from the direct simulation Monte Carlo (DSMC) method have also been used to modify the viscosity and thermal conductivity in the NS equations to account for local shear non-equilibrium, improving predictions of hypersonic flows in the near-continuum regime \cite{ou2019dsmc}. 

Recently, learning-based methods have been utilized to construct closed macroscopic systems from kinetic data. For example, a framework where an autoencoder constructs generalized moments is designed in \cite{han2019uniformly}, followed by neural networks that learn the closure terms from the kinetic equation to achieve uniform accuracy across a wide range of Knudsen numbers. In \cite{Li2023ipnc}, the neural network used to approximate the moment closure relation is specifically designed to preserve Galilean, reflection, and scaling invariance. Moreover, learned moment closures have been constrained to preserve hyperbolicity, improving the stability of the resulting macroscopic systems and supporting predictions beyond the training window \cite{christlieb2025hyperbolic, huang2026machine}. In addition, a differentiable programming framework combining automatic differentiation and adjoint methods is developed for the end-to-end optimization \cite{xiao2025differentiable}.

In this paper, a data-driven generalized Navier-Stokes (GNS) model is proposed for flows in the transition regime by learning a neural network closure that gives the relationship between the stress tensor, heat flux, and the macroscopic variables, namely density, macroscopic velocity, and temperature. Precisely, the spatial derivatives of the macroscopic velocity are utilized instead of the macroscopic velocity as the input variables to preserve Galilean invariance. The U-Net architecture is adopted to capture the multiscale spatial features, with techniques such as input normalization and instance normalization utilized to improve training stability and efficiency. For the loss function, the standard end-to-end learning method is utilized, along with a new sequential training strategy that progresses from the standard end-to-end learning method to long-term end-to-end learning proposed to improve long-term predictive accuracy. The efficiency of this new GNS model is validated through numerical experiments, including the 1D3V and 2D3V problems.

The structure of this paper is as follows. The macroscopic conservation laws derived from the Boltzmann equation and the associated closure problem are discussed in Sec. \ref{sec:Bol}, and the detailed methodology to obtain the GNS model is presented in Sec. \ref{sec:moment_closure}. In Sec. \ref{sec:nn-training}, the network architecture and the training methodology are proposed, and the study of the numerical examples is shown in Sec. \ref{sec:num}, with the conclusion in Sec. \ref{sec:con}.
\section{Macroscopic conservation laws and the closure problem}
\label{sec:Bol}
In this section, the macroscopic conservation laws and the related closure problem are discussed in detail. The unclosed form for the 3D macroscopic conservation laws has the following form 
\begin{equation}
\label{eq:NSF}
\left\{
\begin{aligned}
&\pd{\rho}{t} + \pd{\rho u_k}{x_k} = 0, \\
& \pd{\rho u_i}{t} + \pd{(\rho u_i u_k + p \delta_{ik}+ \sigma_{ik})}{x_k} = 0,\qquad i = 1, 2, 3.\\
& \pd{E}{t} + \pd{(E u_k + (p \delta_{ik} + \sigma_{ik})u_i +q_k)}{x_k} = 0, 
\end{aligned} \right.
\end{equation}
where $\rho$, $u_i$ and $E$ are the density, macroscopic velocity and total energy, respectively, which are also the equilibrium variables. $p$ is the pressure, while $\sigma_{ij}$ and $q_i$ are stress tensor and heat flux, which are also the non-equilibrium variables. These macroscopic variables have the relationship below 
\begin{equation}
\label{eq:sigmaq}
  E  = \frac{3}{2}\rho \theta + \frac{1}{2}\rho |\bu|^2, \qquad p_{ij}= p \delta_{ij} + \sigma_{ij}
\end{equation}
with $p_{ij}$ the pressure tensor. Thus, the ideal gas law gives the relationship between the pressure $p$ and the temperature $\theta$ as 
\begin{equation}
    \label{eq:p-T}
    p = \frac{1}{3}\sum_{i = 1}^{3}p_{ii} = \rho \theta.
\end{equation}
The conservation law \eqref{eq:NSF} is the basic system of the macroscopic system. However, these governing equations of the macroscopic variables $\rho, \bu$ and $\theta$ contain additional variables, such as the heat flux $q_i$ and the stress tensor $\sigma_{ij}$. Therefore, additional equations are required to obtain a closed system. The simplest method is to directly set $q_i$ and $\sigma_{ij}$ to zero, and the Euler equations are derived. If the Fourier law is utilized to obtain a closed system \cite{Struchtrup2005} as  
\begin{equation}
    \label{eq:Fourier-law}
    \sigma_{ij} = -\mu \left(\pd{u_i}{x_j} + \pd{u_j}{x_i}\right) + \frac{2}{3}\mu\delta_{ij}\pd{u_k}{x_k}, \qquad q_i = -\kappa \pd{\theta}{x_i}, \qquad  i,j\in\{1,2,3\},
\end{equation}
where $\mu$ is the viscosity and $\kappa$ is the thermal conductivity, the Navier-Stokes-Fourier (NSF) system is obtained, which is a first-order approximation of the initial Boltzmann equation related to the Knudsen number. However, the closure relationships \eqref{eq:Fourier-law}  lose accuracy as the Knudsen number increases, and the NSF system is limited only to the continuum regime. 

If a more generalized macroscopic system beyond the continuum regime is wanted, more accurate closure relationships for the stress tensor and heat flux are therefore needed in the transition regime. This closure relationship can also be obtained from the Boltzmann equation. To express this, we will begin by deducing the NS system from the Boltzmann equation. Let $\bm{x},\bm{v}\in \mathbb{R}^3$ denote the space variable and microscopic velocity. The Boltzmann equation takes the form 
\begin{equation}
\label{eq:bzm}
    \frac{\partial f}{\partial t}+\bm{v} \cdot \nabla_{\bm{x}} f = Q(f, f),
\end{equation}
where $t$ is time, and $f = f(\bm{x},\bm{v}, t)$ is the distribution function of particles in the phase space. Here, $Q(f,f)$ is the quadratic collision term
\begin{equation}
\label{eq:quad}
    Q(f,f) = \int_{\mathbb{R}^3} \int_{\mathbb{S}^2}\left(f_*^{\prime} f^{\prime}-f_* f\right)B(\cos{\chi},\left|\bm{v}-\bm{v}_*\right| )\dd \Omega \dd \bm{v}_*,
\end{equation}
where $f_*' = f(\bm{x},\bm{v}_*^{\prime}, t)$ with $f^{\prime}$ and $f_*$ defined in a similar way. Here, $(\bm{v},\bm{v}_*)$ and $(\bm{v}^{\prime}, \bm{v}_*^{\prime})$ are pre-collision and post-collision velocities, respectively, which have the following relationship
\begin{equation}
\bm{v}^{\prime}=\frac{\bm{v}+\bm{v}_*}{2}+\frac{\left|\bm{v}-\bm{v}_*\right|}{2} \Omega,\quad \bm{v}_*^{\prime}=\frac{\bm{v}+\bm{v}_*}{2}-\frac{\left|\bm{v}-\bm{v}_*\right|}{2} \Omega,
\end{equation}
where $\Omega \in \mathbb{S}^2$ is a unit vector. The deflection angle $\chi$ between the pre- and post-collision relative velocities satisfies 
\begin{equation}
    \label{eq:angle}
\cos{\chi} =\Omega\cdot\frac{\bm{v}-\bm{v}_*}{\left|\bm{v}-\bm{v}_*\right|}.    
\end{equation}
The variable hard sphere (VHS) model \cite{Bird1994} is utilized in this work, in which the collisional kernel is independent of the deflection angle $\chi$ and is given by
\begin{equation}
    B = C_\alpha \left|\bm{v}-\bm{v}_*\right|^\alpha,\qquad \alpha = \frac{\eta-5}{\eta-1}, \qquad \eta = \frac{2\omega+3}{2\omega-1}.
\end{equation}
Here, $\omega$ is the constant index related to the VHS model, and the constant $C_\alpha$ has the following form \cite{wu2013deterministic}
\begin{equation}
\label{eq:coeC}
    C_\alpha = \frac{5}{64\sqrt{\pi}\sqrt{2}^{\alpha}\Gamma(\frac{\alpha+3}{2})\mu_{\mathrm{ref}}}, \qquad \mu_{\mathrm{ref}} = \frac{15\sqrt{\pi}}{2(5-2\omega)(7-2\omega)}\mathrm{Kn}, 
\end{equation}
with the reference viscosity coefficient $\mu_{\mathrm{ref}}$ . Moreover, the collision term $Q(f,f)$ conserves total mass, momentum, and energy as 
\begin{equation}
\label{eq:con}
\int_{\mathbb{R}^3}\psi(\bv) Q(f,f) \dd\bm{v} = \bm{0}, \qquad \psi(\bv) = \left(1, \bm{v}, \frac12 |\bv|^2\right)^T.
\end{equation}

Multiplying \eqref{eq:bzm} by $\psi(\bv)$ and integrating over $\bm{v}\in\mathbb{R}^3$, the macroscopic system \eqref{eq:NSF} is derived. Moreover, the relationships between the macroscopic variables and the distribution function are listed below
{\small
\begin{equation}
    \label{eq:rel}
    \begin{gathered}
    \rho = \int_{\bbR^3} f \dd \bv, \qquad \rho \bu = \int_{\bbR^3} \bv f \dd \bv, \qquad  
    E = \frac{1}{2} \int_{\bbR^3}\bv^2 f \dd \bv \\
        p_{ij}=   \int_{\bbR^3}(v_i - u_i)(v_j - u_j) f \dd \bv, \qquad 
        q_i = \frac{1}{2}\int_{\mathbb{R}^3}|\bm{v}-\bm{u}|^2(v_i-u_i)f\dd\bm{v},\quad i,j\in\{1,2,3\}.
\end{gathered}
\end{equation}
}
Besides, the steady-state solution to the Boltzmann equation is called Maxwellian,  and has the following form 
\begin{equation}
    \label{eq:Maxwellian}
    \mM(\bv) = \frac{\rho}{(2 \pi \theta)^{3/2}} \exp\left(-\frac{|\bv - \bu|^2}{2 \theta} \right).
\end{equation}

Generally speaking, different collision models bring different closure relationships even with the same closure method. The deduction of the closure relationships is quite complicated for complex collision models \eqref{eq:quad}, and the Fourier law \eqref{eq:Fourier-law} can be deduced from the BGK collision model \cite{Struchtrup2005}.

In this work, a data-driven machine learning method is proposed to obtain the closed Navier-Stokes equations from the data formed by the numerical solution to the Boltzmann equation, in which a neural network is utilized to learn new closure relations for the stress tensor and heat flux. This method is expected to have a large computational regime and can describe the behavior of the gas in the transition regime.

\section{Data-driven moment closure}
\label{sec:moment_closure}

In this section, a generalized NS (GNS) model is built by learning the stress tensor and the heat flux from the Boltzmann equation. The main innovation point is in designing the loss function, where the neural network is trained by minimizing the difference between the numerical solution of the GNS and that of the NS equations derived from the Boltzmann equation. We will begin from the 1D3V Boltzmann equation in Sec. \ref{sec:1D-dim}, and the extension to the multi-dimensional problem is presented in Sec. \ref{sec:multi-dim}.

\subsection{1D3V Boltzmann equation moment closure}
\label{sec:1D-dim}
For the 1D3V problem, the conservation law \eqref{eq:NSF} is reduced into 
\begin{equation}
\label{eq:1D3V}
    \left\{
    \begin{aligned}
        & \frac{\partial\rho}{\partial t}+\frac{\partial(\rho u_1)}{\partial x}  = 0, \\
        & \frac{\partial\rho u_1}{\partial t}+\frac{\partial(\rho u_1^2+p+\sigma_{11})}{\partial x} = 0, \\
        &\frac{\partial E}{\partial t}+\frac{\partial(E u_1+(p+ \sigma_{11})u_1+q_1)}{\partial x}  = 0,
    \end{aligned}
\right.    
\end{equation}
where the related macroscopic variables \eqref{eq:sigmaq} are reduced into
\begin{equation}\label{eq:moment_sq}
    p+\sigma_{11} = \int_{\mathbb{R}^3}(v_1-u_1)^2 f\dd\bm{v},\quad p = \rho \theta, \quad q_1 = \frac{1}{2}\int_{\mathbb{R}^3}((v_1-u_1)^2+v_2^2+v_3^2)(v_1-u_1)f\dd\bm{v}.
\end{equation}
In general, to obtain a closed system of \eqref{eq:1D3V} with a machine learning method, a neural network is designed for $\sigma_{11}$ and $q_1$ as  
\begin{equation}
    \label{eq:neural_newtork}
    (\bar{\sigma}_{11}, \bar{q}_1) := \mathcal{N}_{\xi}(\bm{U}, \Kn),
\end{equation}
where $\xi$ represents the parameters in the neural network. $\bar{\sigma}_{11}$ and $\bar{q}_1$ are neural network approximations for $\sigma_{11}$ and $q_1$ in \eqref{eq:1D3V}, and the input of this neural network is the macroscopic variables $\bm{U} = (\rho,  u_1, \theta)$ along with the Knudsen number $\Kn$, which is added to indicate the fluid regime. With conservative variables defined as $\bm{W}=(\rho,\rho u_1,E)^{\top}$, the system \eqref{eq:1D3V} is rewritten as
\begin{equation}
\label{eq:1D3V_NN}
    \frac{\partial \bm{W}}{\partial t}+\frac{\partial\bm{F}^{\text{Euler}}\left(\bm{W}\right)}{\partial x}+\frac{\partial\bm{F}^{\text{NN}}\left(\bm{W},\text{Kn}\right)}{\partial x} = 0,
\end{equation}
with 
\begin{equation}
     \left\{\begin{array}{l}
    \bm{F}^{\text{Euler}}\left(\bm{W}\right)=\left(\rho u_1,\rho u_1^2+p,\left(E+p\right)u_1\right)^{\top}, \\        \bm{F}^{\text{NN}}\left(\bm{W}, \text{Kn}\right)=\left(0,\bar{\sigma}_{11},\bar{\sigma}_{11}u_1+\bar{q}_1\right)^{\top}.
\end{array}\right.
\end{equation}
The finite volume method is utilized to numerically solve \eqref{eq:1D3V_NN}, where the Kurganov and Tadmor's (KT) central scheme \cite{kt2000central} is adopted for the numerical flux. Precisely, let
\begin{equation}
    \tilde{\bm{F}} :=\bm{F}^{\text{Euler}}+\bm{F}^{\text{NN}},
\end{equation}
and the  numerical flux of \eqref{eq:1D3V_NN} at $x_{i+1/2}$ is given as
\begin{equation}
\label{eq:num_flux}
\begin{split}
&    \tilde{\bm{F}}_{i+1/2} := \tilde{\bm{F}}\left(\bm{W}_{i+1/2}^{+}, \bm{W}_{i+1/2}^{-},\text{Kn}\right) \\
    &\quad = \frac{1}{2}\left(\bm{F}^{\text{Euler},+}_{i+1/2} + \bm{F}^{\text{Euler},-}_{i+1/2} + \bm{F}^{\text{NN},+}_{i+1/2} + \bm{F}^{\text{NN},-}_{i+1/2}- a_{i+1/2}(\bm{W}_{i+1/2}^{+} - \bm{W}_{i+1/2}^{-})\right),
\end{split}
\end{equation}
with 
\begin{equation}
    \label{eq:flux}
   \bm{F}^{\text{Euler},\pm}_{i+1/2} =  \bm{F}^{\text{Euler}}(\bm{W}_{i+1/2}^{\pm}), \qquad\bm{F}^{\text{NN}, \pm}_{i+1/2} =  \bm{F}^{\text{NN}}(\bm{W}_{i+1/2}^{\pm},\text{Kn})
\end{equation}
where $\bm{W}_{i+1/2}^{\pm}$ denote the reconstructed numerical solution at the cell interface $x_{i+1/2}$ from the left and right sides. Here, the piecewise linear reconstruction with the van Leer limiter is utilized. $a_{i+1/2}$ is the maximum local wave speed in the original Kurganov and Tadmor's (KT) central scheme, and a constant large enough is utilized here for easy computation, which only gives a more dissipative numerical flux. 
With this numerical flux, \eqref{eq:1D3V_NN} is updated as 
\begin{equation}\label{eq:macroeq}
    \bm{W}_{{\rm GNS}, i}^{n+1}=\bm{W}_{{\rm GNS}, i}^n+\frac{\Delta t}{\Delta x}\left(\tilde{\bm{F}}_{i-1 / 2}^{n}-\tilde{\bm{F}}_{i+1 / 2}^{n}\right),
\end{equation}
with $\bm{W}_{{\rm GNS}, i}^{n}$ the numerical solution of GNS \eqref{eq:1D3V_NN} at $t^n$ on the $i$-th spatial mesh.

To obtain the target variables $\bm{W}$ utilized in the loss function, instead of directly obtaining them from the distribution function by \eqref{eq:rel} at each time step, they are obtained by numerically solving the NS system \eqref{eq:1D3V} with the information from the Boltzmann equation. Precisely, the macroscopic variables at time $t=t^n$ and the $i$-th spatial mesh are obtained by the finite volume method as 
\begin{equation}
\label{eq:microeq}
    \bm{W}_{{\rm NS}, i}^{n+1}=\bm{W}_{{\rm NS}, i}^n+\frac{\Delta t}{\Delta x}\left(\bm{F}_{i-1 / 2}^n-\bm{F}_{i+1 / 2}^n\right). 
\end{equation}
Here, $\bm{W}_{{\rm NS}, i}^{n}$ is the numerical solution of NS \eqref{eq:1D3V} at $t^n$ on the $i$-th spatial mesh. $\bm{F}_{i-1 / 2}^n,\bm{F}_{i+1 / 2}^n$ are the numerical fluxes across the cell interfaces obtained as
\begin{equation}
    \label{eq:flux_NS}
    \bm{F}_{i\pm 1 / 2}^n=\int_{\mathbb{R}^3}(1,v_1,\frac{1}{2}|\bm{v}|^2)^{\top}v_1f_{i\pm 1 / 2}^n\dd\bm{v},
\end{equation}
where $f_{i\pm 1 / 2}^n$ are the reconstructed numerical solutions of the Boltzmann equation at the cell interfaces $x_{i\pm 1/2}$. Here, the piecewise linear reconstruction with the van Leer limiter is utilized. 

\begin{remark}
$\bm{W}_{{\rm GNS},i}^n$ and $\bm{W}_{{\rm NS},i}^n$ are utilized in the loss function, and the detailed construction of the loss function is proposed in Sec. \ref{sec:nn-training}. For the input data in the neural network \eqref{eq:neural_newtork}, the macroscopic variables $\bm{U}$ are obtained by $\bm{W}_{{\rm GNS},i}^n$ in GNS. 
\end{remark}





\subsection{Extension to spatial multi-dimension problem}
\label{sec:multi-dim}
The extension of GNS to the spatial multi-dimension problem is natural. Defining $\bm{W}=(\rho,\rho u_1,\rho u_2,E)^{\top}$, the GNS model \eqref{eq:1D3V_NN} is changed into \eqref{eq:2D3V_NN} for the 2D3V Boltzmann equation as below
\begin{equation}\label{eq:2D3V_NN}
    \frac{\partial \bm{W}}{\partial t}+\frac{\partial\bm{F}^{\text{Euler}}\left(\bm{W}\right)}{\partial x}
    +\frac{\partial\bm{G}^{\text{Euler}}\left(\bm{W}\right)}{\partial y}+
    \frac{\partial\bm{F}^{\text{NN}}\left(\bm{W},\text{Kn}\right)}{\partial x} + 
    \frac{\partial\bm{G}^{\text{NN}}\left(\bm{W},\text{Kn}\right)}{\partial y}= 0,
\end{equation}
with 
\begin{equation}
     \left\{\begin{array}{l}
         \bm{F}^{\text{Euler}}\left(\bm{W}\right):=\left(\rho u_1,\rho u_1^2+p,\rho u_1u_2,\left(E+p\right)u_1\right)^{\top}, \\        
         \bm{G}^{\text{Euler}}\left(\bm{W}\right):=\left(\rho u_2,\rho u_1u_2,\rho u_2^2+p,\left(E+p\right)u_2\right)^{\top},\\
        \bm{F}^{\text{NN}}\left(\bm{W},\text{Kn}\right):= \left(0,\bar{\sigma}_{11},\bar{\sigma}_{12},\bar{\sigma}_{11} u_1+\bar{\sigma}_{12} u_2+\bar{q}_{1}\right)^{\top},\\
         \bm{G}^{\text{NN}}\left(\bm{W},\text{Kn}\right):= \left(0,\bar{\sigma}_{21},\bar{\sigma}_{22},\bar{\sigma}_{21} u_1+\bar{\sigma}_{22} u_2+\bar{q}_{2}\right)^{\top}.
\end{array}\right.
\end{equation}
For the two-dimensional case, the neural network \eqref{eq:neural_newtork} is changed into 
\begin{equation}
    \label{eq:net_2D}
    (\bar{\sigma}_{ij}, \bar{q}_i) := \mathcal{N}_{\xi}(\bm{U}, \Kn), \qquad  i,j = 1, 2,
\end{equation}
where $\bar{\sigma}_{ij}$ and $\bar{q_i}$ denote the neural network approximations of the stress tensor $\sigma_{ij}$ and heat flux $q_i$ respectively, and $\bm{U} = (\rho, u_1, u_2, \theta)$ are the input variables. Since the physical stress tensor satisfies $\sigma_{12}=\sigma_{21}$, only three independent stress quantities need to be predicted, and $\bar{\sigma}_{21}$ is accordingly taken as $\bar{\sigma}_{12}$.  With the application of the KT central scheme to each spatial dimension, the numerical solution to \eqref{eq:2D3V_NN} becomes 
\begin{equation}
 \label{eq:macroeq2D}
    \bm{W}_{{\rm GNS}, i,j}^{n+1}=\bm{W}_{{\rm GNS},i,j}^n+\frac{\Delta t}{\Delta x}\left(\tilde{\bm{F}}_{i-1 / 2,j}^{n}-\tilde{\bm{F}}_{i+1 / 2,j}^{n}\right)+\frac{\Delta t}{\Delta y}\left(\tilde{\bm{G}}_{i,j-1 / 2}^{n}-\tilde{\bm{G}}_{i,j+1 / 2}^{n}\right),
\end{equation}
where $\tilde{\bm{F}}_{i\pm1 / 2,j}$ and $\tilde{\bm{G}}_{i,j\pm1 / 2}$ denote the corresponding numerical fluxes on the cell interfaces at $(i\pm 1 / 2,j)$ and $(i,j\pm 1 / 2)$ respectively. The same local Lax-Friedrichs numerical flux and piecewise linear reconstruction as in \eqref{eq:num_flux} are utilized here to obtain $\tilde{\bm{F}}_{i\pm1 / 2,j}$ and $\tilde{\bm{G}}_{i,j\pm1 / 2}$.

To obtain the target data $\bm{W}_{{\rm NS},i,j}^n$, a similar method as in \eqref{eq:microeq} is utilized as 
\begin{equation}
\label{eq:microeq2D}
\bm{W}_{{\rm NS},i,j}^{n+1}=\bm{W}_{{\rm NS}, i,j}^n+\frac{\Delta t}{\Delta x}\left(\bm{F}_{i-1/2,j}^n-\bm{F}_{i+1/2,j}^n\right)+\frac{\Delta t}{\Delta y}\left(\bm{G}_{i,j-1/2}^n-\bm{G}_{i,j+1/2}^n\right).
\end{equation}
Here, $\bm{F}_{i\pm 1/2,j}^n$ and $\bm{G}_{i,j\pm 1/2}^n$ are the numerical fluxes obtained as 
\begin{equation*}
\label{eq:num_flux_2D}
\bm{F}_{i\pm 1/2,j}^n=\int_{\mathbb{R}^3}\phi v_1f_{i\pm 1/2,j}^n\,\dd\bm v,\qquad 
\bm{G}_{i,j\pm 1/2}^n=\int_{\mathbb{R}^3}\phi v_2f_{i,j\pm 1/2}^n\,\dd\bm v, 
\end{equation*}
where $\phi = (1, v_1, v_2, \frac{1}{2}|\bm v|^2)^{\top}$. $f_{i\pm 1 / 2,j}^n$ and $f_{i,j\pm 1 / 2}^n$ are the reconstructed distribution functions on the cell interfaces at $(i\pm 1 / 2,j)$ and $(i,j\pm 1 / 2)$ respectively. 

\begin{remark}
    For the input data in the neural network \eqref{eq:net_2D}, it is also obtained by $\bm{W}_{{\rm GNS}, i,j}^n$. 
\end{remark}

For now, the method to obtain the GNS model by the data-driven method is proposed, and the only remaining point is how to design the loss function, which will be introduced in detail in Sec. \ref{sec:nn-training}.
\section{Neural network architecture and training methodology}
\label{sec:nn-training} 
In this section, the neural network architecture and training methodology are detailed for constructing the generalized Navier-Stokes model (GNS). First, the network architecture, including the physical constraints that the neural network closure model should satisfy to preserve Galilean invariance, the U-Net architecture utilized, and the techniques to stabilize and accelerate the learning process, are presented in Sec. \ref{sec:network}. The detailed strategy to design the loss function is then introduced in Sec. \ref{sec:training}, including the loss function obtained by standard end-to-end learning in Sec. \ref{sec:short-term} and a sequential long-term learning strategy proposed in Sec. \ref{sec:long-term} to reduce the accumulated errors and improve long-term predictive accuracy.

\subsection{Neural network architecture}
\label{sec:network}
In this section, the neural network architecture to obtain GNS is proposed. To obtain a closed system, 
$\bar{\sigma}_{ij}$ and $\bar{q}_{i}$ are taken as the output features in the neural network, while the macroscopic variables $\rho,\bu,\theta$ and the Knudsen number are taken as the input features, which have the following form 
\begin{equation}\label{eq:closure_form}
    (\bar{\sigma}_{ij}, \bar{q}_{i}) = \mathcal{N}_{\xi}(\rho,\bu,\theta,\Kn).
\end{equation}
Before introducing the detailed architecture of $\mN_{\xi}(\cdot)$, the constraints that the neural network closure should satisfy to maintain Galilean invariance are discussed. Galilean invariance is a fundamental physical property of the Boltzmann equation, and preserving this invariance has been proved important when designing the closed moment system with machine learning methods \cite{han2019uniformly, Li2023ipnc}. 

\paragraph{Preserving Galilean invariance}
Galilean invariance requires that the equation remains unchanged under the Galilean group of transformations, which includes translation, rotation, and the Galilean boost. In the present GNS model, translational invariance and Galilean boost invariance are considered. The former can be naturally maintained by a fully convolutional neural network, in which the same convolutional kernels are shared across different spatial locations.

For Galilean boost invariance, if $f(\bmx,\bv,t)$ is a solution of the Boltzmann equation, then for any constant vector $\bu_0 \in \bbR^3$,
\begin{equation}
    \label{eq:tilde_f}
    \tilde{f}(\bmx, \bv, t) = f(\bmx-t\bu_0,\bv-\bu_0,t)
\end{equation}
is also a solution to the Boltzmann equation. This property imposes constraints on the macroscopic variables of $\tilde{f}(\bmx, \bv, t)$ and $f(\bmx, \bv, t)$ as 
\begin{gather}
\label{eq:tilde_macro}
\tilde{\rho}(\bmx,t) = \rho(\bmx-t\bu_0, t), \qquad  \tilde{\bu}(\bmx,t) = \bu(\bmx-t\bu_0, t) + \bu_0, \\
 \tilde{\theta}(\bmx,t)   = \theta(\bmx-t\bu_0, t) , \qquad 
 \tilde{\sigma}_{ij}(\bmx,t)  = \sigma_{ij}(\bmx-t\bu_0,t), \qquad 
 \tilde{q}_{i}(\bmx,t) = q_{i}(\bmx-t\bu_0,t).
\end{gather}
Therefore, to preserve invariance under Galilean boost, the network $\mN_{\xi}(\cdot)$ should satisfy the following constraints:
$\forall \bu_0 \in \bbR^3,$
\begin{equation}
\label{eq:constraint}
    \mN_{\xi}(\rho,\bu,\theta, \Kn) = \mN_{\xi}(\rho,\bu-\bu_0,\theta,\text{Kn}).
\end{equation}
To achieve this, the input feature $\bu$ is replaced by $\bu_x$ and $\bu_y$ instead as in \cite{Li2023ipnc, han2019uniformly}.

Moreover, $\ln{\Kn}$ instead of $\Kn$  is utilized in the input feature because the Knudsen number spans several orders of magnitude, and the logarithmic form has been empirically demonstrated to improve the prediction accuracy \cite{han2019uniformly, Li2023ipnc}. Specifically, take the 2D3V problem as an example. For the neural network $\mN_{\xi}(\cdot)$, the input $\mathbf{X}$ has the form as 
\begin{equation}
\label{eq:input2D}
    \mathbf{X}_{k,l} = [\rho,\bu_x,\bu_y,\theta,\ln{\text{Kn}}]_{k,l},
    \qquad k = 1,\cdots,N_x,\quad l = 1,\cdots,N_y,
\end{equation}
which is a tensor with size $N_{\rm in} \times N_x\times N_y$, where $N_x$ and $N_y$ are the spatial mesh sizes in the \(x\)- and \(y\)-directions. $N_{\rm in}$ is the number of the input features. Here, $N_{\rm in} = 7$ for the 2D3V problem, and $N_{\rm in} = 4$ for the 1D3V problems. Here, $\bu_x$ and $\bu_y$ are  approximated with the discretized differences in the numerical implementation.

\paragraph{Input normalization}
Input normalization is essential due to the large magnitude differences among the input features. Each input feature is normalized with its mean and variance over the entire training dataset. For the training dataset, it is derived from the numerical solution of the Boltzmann equation with various initial conditions. Here,  this numerical solution is obtained by the fast Fourier spectral method \cite{jingweihu2017}. Let the density $\rho$ be an example. Assume $N_t$ is the total time step, and $N_{\mathrm{ic}}$ is the total number of initial condition samples. Then the total training data set for the density is as follows
 \begin{equation}
     \label{eq:train_rho}
     \rho_{i,k,l}^j, \qquad i = 1, \cdots, N_{\rm ic}, ~ j = 1, \cdots, N_t, ~ k = 1, \cdots, N_x, ~ l = 1, \cdots, N_y
 \end{equation}
Thus, the input normalization $\mathcal{P}(\cdot)$ at spatial mesh $(k,l)$, time step $j$ with the $i$-th initial condition is calculated as 
\begin{equation}
\label{eq:normalized_rho}
\mathcal{P}(\mathbf \rho):= (\rho_{i,k,l}^j)_{\rm norm} = \frac{\rho_{i,k,l}^j-\Bar{\rho}}{\sigma_{\rho}}, 
\end{equation}
with
\begin{equation}
\bar{\rho} = \Pi \rho_{i,k,l}^j,\quad  \sigma^2_{\rho} = \Pi
    (\rho_{i,k,l}^j-\bar{\rho})^2+\epsilon,\qquad \Pi :=\frac{1}{N_{\mathrm{ic}}N_tN_xN_y}\sum\limits_{i=1}^{N_{\mathrm{ic}}}\sum\limits_{j=1}^{N_t}\sum\limits_{k=1}^{N_x}\sum\limits_{l=1}^{N_y}. 
\end{equation}
The parameter $\epsilon$ is a small constant that prevents division by zero. 

\paragraph{Instance normalization}
Instance normalization  (IN) was introduced by \cite{ulyanov2016instance} and has been shown to stabilize and accelerate the training of deep neural networks \cite{huang2023normalization}. Unlike input normalization, which is applied only once in the total learning process, IN is performed after each convolution except for the first one. The detailed learning process will be introduced in the following paragraphs.  

For each implementation, IN is performed independently for each intermediate feature channel over the spatial mesh. Take the spatially 2D case as an example. Assume $\bm{z}=(z_{k,l})_{1\leqslant k\leqslant N_x,\,1\leqslant l\leqslant N_y}$ is an intermediate output over the spatial space for one output feature. Instance normalization $\mathcal{I}(\cdot)$ is defined as
\begin{equation}
    \mathcal{I}(\bm{z})=\frac{\bm{z}-\mu(\bm{z})}{s(\bm{z})}, 
    \label{eq:IN}
\end{equation}
where
\begin{equation}
 \mu(\bm{z})=\frac{1}{N_xN_y}\sum_{k=1}^{N_x}\sum_{l=1}^{N_y} z_{k,l},
    \qquad
    s^2(\bm{z})=\frac{1}{N_xN_y}\sum_{k=1}^{N_x}\sum_{l=1}^{N_y}\left(z_{k,l}-\mu(\bm{z})\right)^2+\epsilon.
\end{equation}
$\epsilon$ is a small number added to prevent division by zero. Here, IN is applied on each output feature individually.

\begin{figure}[htbp]
    \centering
    \includegraphics[width=0.5\textwidth]{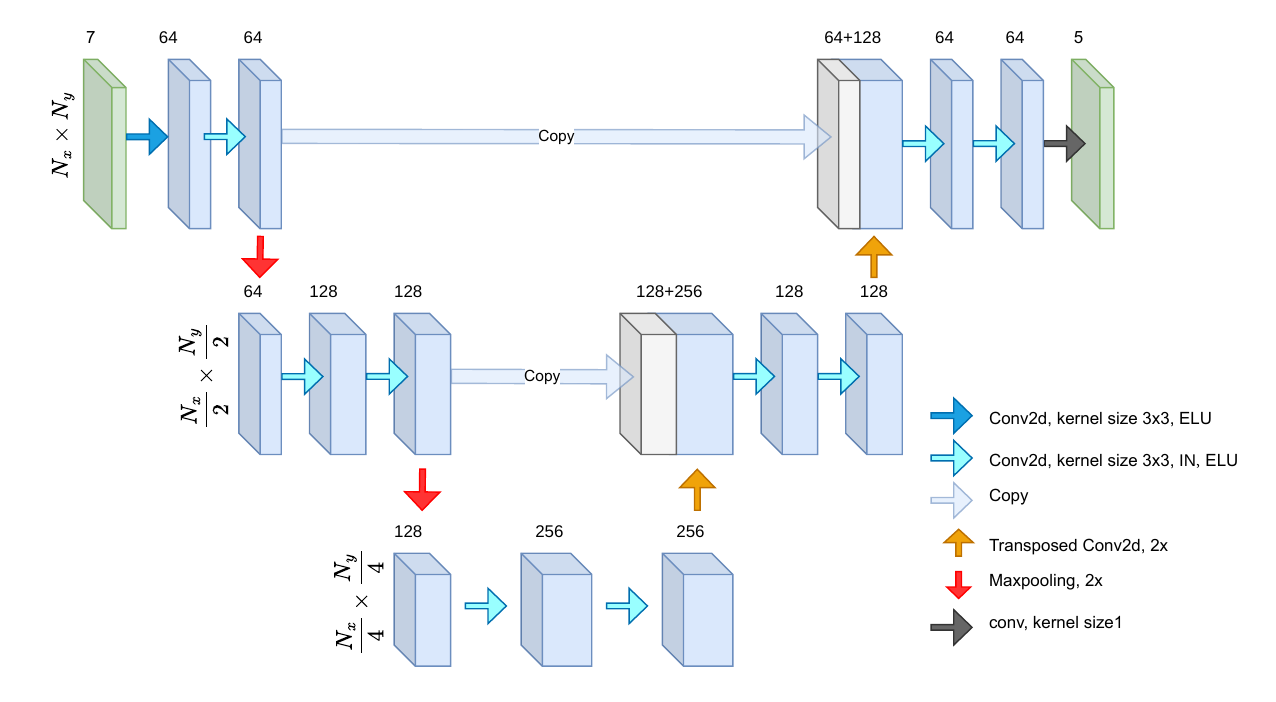}
    \caption{The U-Net architecture. The blocks denote the feature tensors, while the numbers above the blocks indicate the corresponding numbers of channels, and the spatial resolutions at different levels are shown on the left.}
    \label{fig:network_arch}
\end{figure}

\paragraph{Network architecture} 
For the network architecture, the U-Net \cite{ronneberger2015u} is utilized to map the normalized macroscopic input tensor to the output stress and heat flux. U-Net is a fully convolutional neural network architecture that contains three parts, namely the encoder, the decoder, and the skip connection. The U-Net architecture employed in the numerical experiments for the 2D3V problem is shown in Fig. \ref{fig:network_arch} as a representative example.

The basic operation in U-Net is convolution. Here, the convolution with the $3\times 3$ kernel is utilized. Precisely, for an input tensor ${\mathbf Z} \in \bbR^{d_{\rm in} \times N_x \times N_y}$ and an output tensor ${\mathbf Y} \in \bbR^{d_{\rm out} \times N_x \times N_y}$, let $\widetilde{\mathbf Z}$ denote the feature tensor $\mathbf Z$ after padding operation \cite{dumoulin2016guide}. The padding strategy is chosen according to the boundary condition, for example, circular padding is used for periodic boundary conditions. Thus, this 2D convolution has the form as 
\begin{equation}
   \label{eq:conv2D}
    Y_{n,k,l}=b_n+\sum_{d=1}^{d_{\rm in}}\sum_{r=-1}^{1}\sum_{s=-1}^{1}A_{n,d,r,s} \widetilde{Z}_{d,k+r,l+s},
    \end{equation}
where $1\leqslant n\leqslant d_{\rm out}, 1\leqslant k \leqslant N_x, 1 \leqslant l \leqslant N_y$. $A_{n,d,r,s}$ and $b_{n}$ are the trainable convolution kernel and bias. $d_{\rm in}$ and $d_{\rm out}$ are the numbers of input channels and output channels in one block. Here, the convolution is labeled as 
\begin{equation}
    \label{eq:oper_con}
    {\mathbf Y} = \mathcal{L}{\mathbf Z}. 
\end{equation}
Generally speaking, IN is applied after the convolution \eqref{eq:oper_con}. Since $\ln\Kn$ in \eqref{eq:input2D} is spatially constant, the application of IN after the first convolution will remove its contribution to the neural network. Therefore, the application of IN will be omitted after the first convolution. Consequently, the application of the U-Net neural network is summarized as below
\begin{enumerate}
    \item  first step: 
    \begin{equation}
        \label{eq:first_step}
      {\mathbf Z} = \phi(\mathcal{L} \mathcal{P} ({\mathbf X})),
    \end{equation}
    where  ${\mathbf X}$ is the initial training data, and  $\phi$ is the activation function. Here, ELU \cite{clevert2016fast} is utilized as 
    \begin{equation}
        \label{eq:ELU}
        \operatorname{ELU}(x,\alpha) = \begin{cases}
            x,  &{\rm if } ~x > 0, \\
            \alpha(e^{x}  -1 ), & {\rm if }~ x \leqslant 0.
        \end{cases}
    \end{equation}
    \item intermediate step: 
    \begin{equation}
        \label{eq:sec-step}
    {\mathbf Z} = \phi(\mathcal{I} (\mathcal{L} \mathbf Z)),         
    \end{equation}
where the operator $\mathcal L$ first pads the feature tensor $\mathbf Z$ to obtain $\widetilde{\mathbf Z}$ and then applies the convolution \eqref{eq:conv2D}. 
\item final step: in the final step, a kernel-one convolution is utilized to map the feature tensor to five output features, and we will omit its operator form here. 

Thus, the total process of applying the U-Net neural network  is summarized as 
\begin{equation}
    \left[\bar{\sigma}_{11},\bar{\sigma}_{12},\bar{\sigma}_{22},\bar q_1,\bar q_2\right] = \mathcal{N}_{\xi}(\mathbf X) \in\mathbb R^{5\times N_x\times N_y}.
    \label{eq:output2D}
\end{equation}
\end{enumerate}
As shown in Fig.~\ref{fig:network_arch}, the encoder progressively reduces the spatial resolution of the input through a series of convolutional and pooling layers, while simultaneously increasing the number of feature channels to capture more complex features. Symmetrically, the decoder progressively increases the spatial resolution through a series of convolutional and transposed convolutional layers to reconstruct an output with the same resolution as the input. 

The skip connections \cite{taghanaki2019} link the output of each encoder level directly to the corresponding decoder level, allowing the decoder to utilize the high-resolution spatial information and alleviating the vanishing and exploding gradient problems.



\subsection{Designing the loss function}\label{sec:training}
In this subsection, the detailed strategy to design the loss function for GNS will be introduced. Minimizing the distance over the entire time interval is an intuitive choice for the loss function, but it creates a very deep computational graph. Backpropagation through such a deep computational graph is computationally expensive and can also suffer from vanishing or exploding gradients \cite{gomez2017reversible, pascanu2013difficulty}. Therefore, a methodology called short-term end-to-end learning is first introduced in Sec.~\ref{sec:short-term}, followed by long-term end-to-end learning in Sec.~\ref{sec:long-term} to further control the accumulated prediction error over the long-time evolution.

Before introducing the loss function, we first define the evolution operator $\mathcal{T}_\xi$ in one time step as 
\begin{equation}
    \bm{U}^{j+1}=\mathcal{T}_\xi(\bm{U}^{j}, \Kn, \Delta x, \Delta y, \Delta t)\triangleq \mathcal{T}_\xi(\bm{U}^{j}).
    \label{eq:evolution_operator}
\end{equation}
The computational procedure of $\mathcal{T}_\xi$ is summarized in Alg. \ref{alg:evolution_algorithm}, together with the corresponding computational workflow. This procedure includes the transform between the primitive variables $\bm{U}^{j}$ and the conservative variables $\bm{W}^{j}$, the update of $\bm{W}^{j+1}$ by \eqref{eq:macroeq2D}, and the transform between $\bm{W}^{j+1}$ and $\bm{U}^{j+1}$. The corresponding multi-step evolution operator is defined recursively as
\begin{equation}
\label{eq:multi-oper}
    \mathcal{T}_\xi^1:=\mathcal{T}_\xi,\qquad \mathcal{T}_\xi^l:=\mathcal{T}_\xi\circ\mathcal{T}_\xi^{l-1},\quad l\geqslant 2.
\end{equation}

 \begin{algorithm}[H]
            \footnotesize
            \caption{One-step GNS evolution operator $\mathcal{T}_\xi$}
            \label{alg:evolution_algorithm}
            \begin{algorithmic}[1]
                \Statex \textbf{Input:} 
                $\bm{U}^{j}$, $\Kn$, $\Delta x$, $\Delta y$, $\Delta t$
                \Statex \textbf{Output:} 
                $\bm{U}^{j+1}$
                \State Obtain the network input $\mathbf{X}^{j}$ from $\bm{U}^{j}$ by \eqref{eq:input2D}.
                \State Apply input normalization on $\mathbf {X}$ \eqref{eq:normalized_rho} as $\mathcal{P}(\mathbf X)$. 
                \State Obtain  $(\bar{\bm{\sigma}}^{j},\bar{\bm{q}}^{j})$ by the U-Net neural network as \eqref{eq:output2D}.
                \State Incorporate $(\bar{\bm{\sigma}}^{j},\bar{\bm{q}}^{j})$ into the GNS fluxes \eqref{eq:2D3V_NN} and compute the numerical fluxes by \eqref{eq:num_flux}.
                \State Convert the macroscopic variables $\bm{U}^{j}$ to the conservative variables $\bm{W}^{j}$ and update them to $\bm{W}^{j+1}$ by \eqref{eq:macroeq2D}.
                \State Convert $\bm{W}^{j+1}$ back to $\bm{U}^{j+1}$.
            \end{algorithmic}
        \end{algorithm}   

\subsubsection{Short-term end-to-end learning and GNSS model}
\label{sec:short-term}
End-to-end learning is also known as "trajectory fitting" \cite{melchers2023comparison}, "solver-in-the-loop" \cite{um2020solver} or "embedded training" \cite{macart2021embedded}. In this work, this method is called short-term end-to-end learning because the distance between the predicted trajectories and the reference trajectories is computed over many short time fragments instead of the entire time step. Here, different trajectories mean different initial conditions, and the macroscopic model obtained by this short-term end-to-end learning is called GNSS for short. 
Let 
\begin{equation}
    \label{eq:index-set}
    \mathcal{I}=\{1,2,\ldots,N_{\mathrm{ic}}\}, \qquad \mathcal{J}=\{0,1,\ldots,N_t-1\}
\end{equation}
denote the index set of the different initial conditions and the time index set, respectively. The total training data set is combined as 
\begin{equation}
\label{eq:train_set}
\mathcal{D} = \left\{ \bm{U}_i^j: i\in\mathcal{I},\ j\in\mathcal{J} \right\},
\end{equation}
where $\bm{U}_i^j$ denotes the macroscopic variables of the
reference solution at time step $j$ for the $i$-th initial condition. For a trajectory length $K_s$, the admissible starting indices are
\begin{equation}
\mathcal{A}_{\mathrm{short}} = \left\{ (i,j)\in\mathcal{I}\times\mathcal{J}: j+K_s\in\mathcal{J} \right\}.
\end{equation}
In each optimization step, a mini-batch 
\begin{equation}
    \label{eq:mini_bash}
\mathcal{B}=\{(i_m,j_m)\}_{m=1}^{M}
\end{equation}
is randomly sampled from the admissible starting set $\mathcal{A}_{\mathrm{short}}$, following standard stochastic mini-batch optimization \cite{bottou2018optimization}. The short-term loss function is defined as
\begin{equation}\label{eq:L_short}
L_{\mathrm{short}} (\xi;\mathcal{B}) = \frac{1}{M} \sum_{m=1}^{M} \sum_{k=1}^{K_s} \left\| \mathcal{T}_\xi^k \left(\bm{U}_{i_m}^{j_m}\right) - \bm{U}_{i_m}^{j_m+k} \right\|_{{x}}^2.
\end{equation}
Here, the discrete spatial norm $||\bm{U}||_x$ is defined as 
\begin{equation}
\label{eq:discrete_spatial_norm}
    \|\bm{U}\|_{{x}}^2=\frac{1}{N}\sum_{l=1}^{N}\|\bm{U}_l\|_2^2,
\end{equation}
where $N$ is the total mesh size in the spatial space, and $\bm{U}_l$ is the numerical solution of $\bm{U}$ at the $l$-th mesh. 
This loss function \eqref{eq:L_short} computes the error over the trajectory segments of length $K_s$. 
\begin{remark}
    The choice of $K_s$ must balance computational efficiency and predictive accuracy. On one hand, $K_s$ needs to be relatively large to capture the intrinsic evolutionary behavior and predict future dynamics. On the other hand, an excessively large $K_s$ incurs high computational costs and gradient instabilities associated with a deep computational graph. 
\end{remark}

\subsubsection{Long-term end-to-end learning and the GNSL model}
\label{sec:long-term}
For short-term end-to-end learning, it is always expensive to capture long-term behavior, and the generalization ability is limited. Thus, long-term end-to-end learning is proposed in this section, which aims to minimize the error over a long, or even the entire time step, which is also expected to capture the long time evolutionary behavior. Accordingly, the obtained macroscopic model by the long-term end-to-end learning is called GNSL for short. 

To construct the loss function, let $\hat{\bm{U}}_i^{j,l}$ denote the predicted states generated by the evolution operator 
\begin{equation}
\label{eq:pre_state}
    \hat{\bm{U}}_i^{j,l} \triangleq \mathcal{T}_\xi^l(\bm{U}_i^{j}),\quad i\in \mathcal{I},\ j\in\mathcal{J},\ l\in \mathbb{N}_+,\ j+l\in\mathcal{J},
\end{equation}
where $\bm{U}_i^{j}$ represents the macroscopic variables at the $j$-th time step with the $i$-th initial condition. For each sampled pair $(i_m,j_m) \in \mathcal{B}$, let $\mathcal{L}_{i_m,j_m}\subseteq\mathbb{N}_+$ be a set of admissible evolution lengths satisfying $j_m+l+K_l\in\mathcal{J}$, where $K_l$ denotes the trajectory length used in $L_{\text{long}}$. The loss function is then defined as 
\begin{equation}
\label{eq:L_long}
    \begin{aligned}
    L_{\text{long}}\left(\xi;\mathcal{B},\{\mathcal{L}_{i_m,j_m}\}_{m=1}^{M}\right)
    ={}& \frac{1}{\sum\limits_{m=1}^{M}|\mathcal{L}_{i_m,j_m}|}
    \sum\limits_{m=1}^{M}\sum\limits_{l\in\mathcal{L}_{i_m,j_m}}\sum\limits_{k=1}^{K_l}
    \left\|\mathcal{T}_\xi^k\left(\operatorname{sg}\left(\hat{\bm{U}}_{i_m}^{j_m,l}\right)\right)-\bm{U}_{i_m}^{j_m+k+l}\right\|_x^2.
    \end{aligned}
\end{equation}
Here, $\operatorname{sg}(\cdot)$ denotes the stop-gradient operator \cite{vandenoord2017neural} satisfying 
\begin{equation}
    \label{eq:oper_sg}
    \operatorname{sg}(\bm{z})=\bm{z}, \qquad \partial\operatorname{sg}(\bm{z})/\partial\bm{z}=\bm{0}.
\end{equation}
It should be emphasized that the definition of the loss function does not change the forward prediction or the value of the trajectory error. In fact,
\begin{equation}
    \mathcal{T}_\xi^k\left(\operatorname{sg}\left(\hat{\bm{U}}_i^{j,l}\right)\right)
    = \mathcal{T}_\xi^k\left(\hat{\bm{U}}_i^{j,l}\right)
    = \hat{\bm{U}}_i^{j,l+k}
    = \mathcal{T}_\xi^{l+k}\left(\bm{U}_i^j\right).
\end{equation}
The key distinction in the long-term end-to-end learning is the usage of these intermediate predicted states $\hat{\bm{U}}_i^{j,l}$ as starting points when calculating trajectory error. Inspired by truncated backpropagation through time (TBPTT) \cite{williams1990efficient}, the gradient backpropagation is restricted to a finite number of evolution steps. Consequently, although $\hat{\bm{U}}_i^{j,l}$ is computationally derived from $\bm{U}_i^{j}$ through $l$ applications of $\mathcal{T}_\xi$,  the computational graph is detached at $\hat{\bm{U}}_i^{j,l}$ and the gradients are backpropagated only through the subsequent $K_l$ steps. To make this explicit, let 
\begin{equation}
    \label{eq:bmz-bmy}
\bm{z}(\xi)=\hat{\bm{U}}_i^{j,l}, \qquad \bm{y}(\xi,\bm{z})=\mathcal{T}_\xi^k\left(\operatorname{sg}(\bm{z})\right). 
\end{equation}
For the trajectory error 
\begin{equation}
    \label{eq:traj-error}
\ell=\|\bm{y}(\xi, \bm{z}) - \bm{U}_i^{j+l+k}\|_x^2,
\end{equation}
the gradient $\mathrm{d}\ell/\mathrm{d}\xi$ with respect to $\xi$ satisfies
\begin{equation}
    \frac{\mathrm{d}\ell}{\mathrm{d}\xi}= \frac{\partial \ell}{\partial \bm{y}}
    \left(\frac{\partial \bm{y}}{\partial \xi}+\frac{\partial \bm{y}}{\partial \operatorname{sg}(\bm{z})}\frac{\partial \operatorname{sg}(\bm{z})}{\partial \bm{z}}\frac{\mathrm{d}\bm{z}}{\mathrm{d}\xi}\right).
\end{equation}
Since $\partial\operatorname{sg}(\bm{z})/\partial\bm{z}=\bm{0}$, the gradient $\mathrm{d}\ell/\mathrm{d}\xi$ utilized in \eqref{eq:L_long} is reduced into
\begin{equation}
    \frac{\mathrm{d}\ell}{\mathrm{d}\xi}
    = \frac{\partial \ell}{\partial \bm{y}}
      \frac{\partial \bm{y}}{\partial \xi},
\end{equation}
Hence, the truncated gradient omits the contribution transmitted through the first $l$ evolution steps and provides an approximation to the full history gradient. Therefore, the depth of the computational graph involved in backpropagation is bounded by $K_l$, rather than $l+K_l$, which reduces the computational burden and alleviates the gradient instability mentioned in the previous section.

\begin{remark}
The motivation for this design is to control the accumulated prediction error over long-time evolution. Unlike short-term learning, $L_{\text{long}}$ starts from the intermediate predicted states $\hat{\bm{U}}_i^{j,l}$, which contain errors accumulated during the preceding $l$ evolution steps. If these states remain close to the corresponding reference states $\bm{U}_i^{j+l}$, they can be regarded as perturbations of the reference trajectory. Minimizing \eqref{eq:L_long} therefore trains the model on states containing accumulated prediction errors and improves the accuracy of the subsequent evolution from such perturbed states, thereby enhancing long-term predictive accuracy.

This interpretation requires a well-pretrained model. Consequently, a sequential training strategy is adopted in the following numerical experiments. Short-term end-to-end learning is first utilized to obtain a pretrained model, followed by long-term end-to-end learning to further reduce accumulated errors.
\end{remark}

 With a well-pretrained model, the following loss function is recommended to control the trajectory error over the entire time interval in the numerical experiments. At each optimization step, a mini-batch $\mathcal{B}_0=\{i_m\}_{m=1}^{M}$ is randomly sampled from $\mathcal{I}$ defined in \eqref{eq:index-set}. For positive integers $K_l$ and $S$ satisfying $K_l<S\leqslant N_t-1$, define the set of starting indices of the sliding windows as
 \begin{equation}
     \label{eq:slid-win}
     \mathcal{L}_0 =\left\{S \cdot n - K_l \mid n=1,\cdots,\lfloor\frac{N_t-1}{S}\rfloor \right\}.
 \end{equation}
 For each sampled trajectory, all evolution lengths in $\mathcal{L}_0$ are included in the loss, which is defined as
\begin{equation}\label{eq:L_long_experiment}
    L_{\text{long}}(\xi;\mathcal{B}_0,\mathcal{L}_0)
    =\frac{1}{M|\mathcal{L}_{0}|}\sum\limits_{m=1}^{M}\sum\limits_{l\in\mathcal{L}_0}\sum\limits_{k=1}^{K_l}\left\|\mathcal{T}_\xi^k\left(\operatorname{sg}\left(\hat{\bm{U}}_{i_m}^{0,l}\right)\right)-\bm{U}_{i_m}^{k+l}\right\|_x^2,
\end{equation}
Thus, the mini-batch sampling is performed only over the reference trajectories, while all time indices in $\mathcal{L}_0$ are utilized for each sampled trajectory in the loss. This setting ensures that supervision is applied uniformly throughout the time interval while maintaining computational efficiency. 
\section{Numerical experiments} 
\label{sec:num}
In this section, several numerical tests are studied to validate this new GNS model, including the 1D3V wave and Riemann problems, and 2D3V isentropic vortex and Taylor-Green vortex problems. 

We follow the training procedure described in Sec.~\ref{sec:long-term}. Initially, the GNS model is trained with the short-term loss function $L_{\text{short}}$, as defined in \eqref{eq:L_short}, until convergence, and this model is called GNSS for short. Then, the training continues with the long-term loss function $L_{\text{long}}$ defined in \eqref{eq:L_long_experiment}, and this model is called GNSL. Unless otherwise specified, $K_s$ is set as $K_s = 4$ in \eqref{eq:L_short} and $(K_l,S)=(2,10)$ in \eqref{eq:L_long_experiment} throughout the numerical experiments. The AdamW optimizer \cite{loshchilov2019decoupled} is utilized for both loss functions. The GNS model is implemented in PyTorch \cite{paszke2019pytorch} and trained on a single NVIDIA A100 GPU with 80~GB of memory.

The U-Net architecture with ELU activation introduced in Sec.~\ref{sec:network} is utilized for both the 1D3V and 2D3V problems. For 1D3V problems, 4 input channels and 2 output channels, with 128, 256, and 512 feature channels at the three resolution levels, are adopted. For 2D3V problems,  7 input channels and 5 output channels, with 64, 128, and 256 feature channels, are adopted at the three resolution levels. 

Moreover, to evaluate the performance of the different models,  the relative error $\mathcal{E}_{\mathrm{s}}^t$ for the density $\rho$, macroscopic velocity $\bu=(u_1,\ldots,u_D)$, and temperature $\theta$ is defined as
{\small 
\begin{equation}\label{eq:relative_error}
    \mathcal{E}_{\mathrm{s}}^t
    =
    \sqrt{
    \frac{1}{2+D}
    \left(
    \frac{\|\rho_{\mathrm{s}}^t-\rho^t\|_{{x}}^2}
         {\|\rho^t\|_{{x}}^2}
    +
    \sum_{i=1}^{D}
    \frac{\|u_{\mathrm{s},i}^t-u_i^t\|_{{x}}^2}
         {\|u_i^t\|_{{x}}^2+c}
    +
    \frac{\|\theta_{\mathrm{s}}^t-\theta^t\|_{{x}}^2}
         {\|\theta^t\|_{{x}}^2}
    \right)
    },
    \quad
    {\rm s} = {\rm GNSL, GNSS, NSF}.
\end{equation}
}
Here, ${(\cdot)}^t_{\rm s}$ is the numerical solution obtained by the model s at time $t$, and ${(\cdot)}^t$ is the reference solution usually obtained by directly solving the Boltzmann equation. $D$ denotes the spatial dimension, and the discrete spatial norm ${\|\cdot\|}_{{x}}$ for a scalar mesh function $\psi$ is defined in \eqref{eq:discrete_spatial_norm}. $c$ is a problem-dependent small constant added to the denominator when the reference discrete norm may approach zero. For the 1D3V problems and the 2D3V Taylor-Green vortex problem, the norm of the velocities may approach zero, and we take $c = 1$, while $c = 0$ for the 2D3V isentropic vortex problem.

\subsection{1D3V wave problem}\label{sec:wave}
In this section, the proposed GNS model is tested on the problem with smooth initial conditions governed by the 1D3V Boltzmann equation, a setup similar to that in \cite{han2019uniformly}. Then, the generalization ability of GNS is validated on the problem with discontinuous initial conditions in Sec.~\ref{sec:discontinuous_test}. 

\subsubsection{Training GNSL}
\label{sec:wave_training}
For this smooth initial condition problem, the computational domain in the spatial space is set to $[-0.5,0.5]$ with a periodic boundary condition. The initial condition is a linear combination of two Maxwellian functions:
\begin{equation}
\label{eq:ex1_ini}
    f(x,\bm{v},0)=\frac{\alpha_1 \mathcal{M}[\bm{U}_1](x,\bm{v})+\alpha_2 \mathcal{M}[\bm{U}_2](x,\bm{v})}{\alpha_1+\alpha_2+10^{-6}},
\end{equation}
where $\alpha_1$ and $\alpha_2$ are uniformly sampled from $[0,1].$ The Maxwellian function $\mathcal{M}[\bm{U}_1]$ for the macroscopic variables $\bm{U}_1=(\rho,u,\theta)$ is given by
\begin{equation}
    \mathcal{M}[\bm{U}_1](x,\bm{v})=\frac{\rho(x)}{\sqrt{2 \pi \theta(x)}^3} \exp \left(-\frac{|v_1-u(x)|^2+|v_2|^2+|v_3|^2}{2 \theta(x)}\right),
\end{equation}
with   
\begin{equation}\label{eq:wave_init}
 \rho(x)=a_\rho \sin \left(2 k_\rho \pi x / \ell+\phi_\rho\right)+b_\rho, \qquad
 u(x)=0, \qquad
 \theta(x)=a_\theta \sin \left(2 k_\theta \pi x / \ell+\phi_\theta\right)+b_\theta,
\end{equation}
where the parameters $(a_\rho,\phi_\rho,b_\rho)$ and $(a_\theta,\phi_\theta,b_\theta)$ are uniformly and independently sampled from $[0.2,0.3]\times[0,2\pi]\times[0.5,0.7],$ while $k_\rho $ and $k_\theta$ are uniformly and independently sampled from the set $\{1,2,3,4\}.$ The parameter $\ell$ represents the length of the spatial domain. The Maxwellian function $\mathcal{M}[\bm{U}_2]$ is randomly generated in the same way.

For the training data, the fast Fourier spectral method \cite{wu2013deterministic} is utilized to solve the full Boltzmann equation \eqref{eq:bzm}. The mesh size in the spatial space is $N_x = 100$, and the microscopic velocity space is truncated to $[-10, 10]^3$ with mesh size $120 \times 40 \times 40$. The Knudsen number $\Kn$ is sampled from a log-uniform distribution over $[0.001, 10]$. Altogether, $500$ samples based on the initial condition \eqref{eq:wave_init} and the sampled $\Kn$ are generated and evolved to $t = 0.1$ with the time step length $\Delta t = 0.0005$ to obtain the training data. The macroscopic variables $\rho, u$, and $\theta$ at each time step for all the $500$ samples constitute the training data set. 

\begin{table}[htbp]
    \centering
    \def\arraystrech{1.5}
   {\footnotesize
    \begin{tabular}{lrrrr}
        \toprule
        \diagbox[width=6em]{model}{time}  & $t=0.05$ & $t=0.10$ & $t=0.15$ & $t=0.20$ \\
        \midrule
        NSF & $23.94$& $22.81$ & $20.55$ & $18.29$ \\
        \midrule
       GNSS & $1.76$ & $2.69$ & $4.20$ & $6.29$ \\
        \midrule
       GNSL & $\bm{1.42}$ & $\bm{1.93}$ & $\bm{2.48}$ & $\bm{3.62}$ \\
        \bottomrule
    \end{tabular}}
     \caption{(1D3V wave problem in Sec.~\ref{sec:wave}) The average relative error $\mE_{\rm s}^{t}$ (in percentages) of the testing data for different models at $t = 0.05, 0.1, 0.15$ and $0.2$.}
    \label{tab:wave_error}
\end{table}

\begin{figure}[htbp]
    \centering
    \begin{subfigure}[b]{0.24\textwidth}
        \centering
        \includegraphics[width=\linewidth]{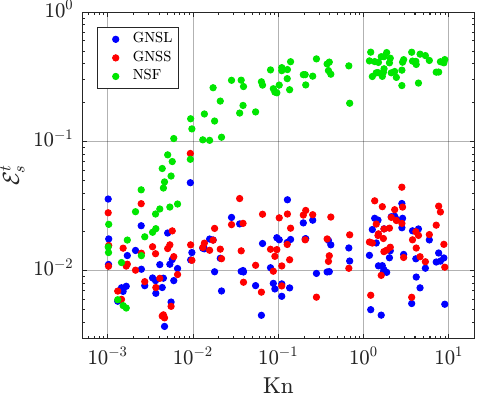}
        \caption{$t = 0.05$}
    \end{subfigure}
    \begin{subfigure}[b]{0.24\textwidth}
        \centering
        \includegraphics[width=\linewidth]{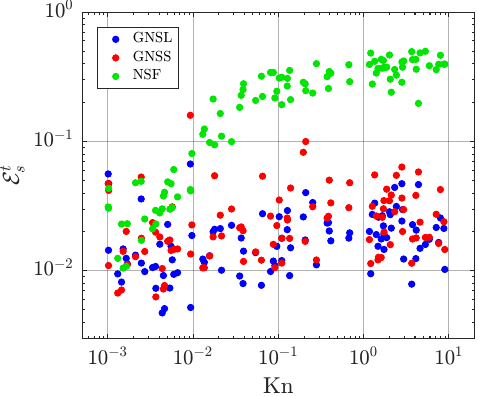}
        \caption{$t = 0.10$}
    \end{subfigure} 
    \begin{subfigure}[b]{0.24\textwidth}
        \centering
        \includegraphics[width=\linewidth]{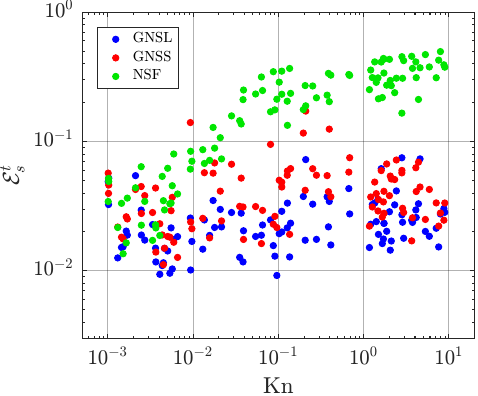}
        \caption{ $t = 0.15$}
    \end{subfigure}
    \begin{subfigure}[b]{0.24\textwidth}
        \centering
        \includegraphics[width=\linewidth]{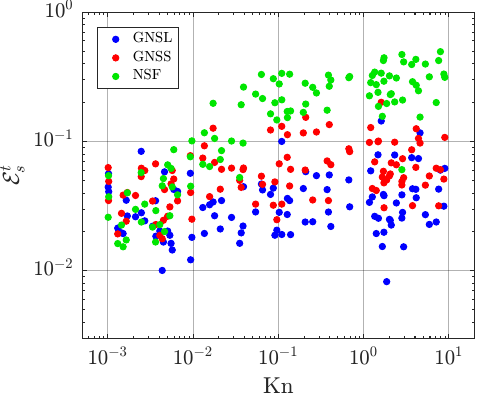}
        \caption{$t = 0.20$}
    \end{subfigure}    
    \caption{(1D3V wave problem in Sec.~\ref{sec:wave}) The distribution of the relative error $\mathcal{E}_{\rm s}^t$ (log-log form) with respect to the Knudsen number for NSF, GNSS, and GNSL at $t=0.05,0.1,0.15$ and $0.2$.}
    \label{fig:1D3Vwave}
\end{figure}

For the testing data, $100$ samples are generated with random initial conditions and Knudsen numbers following the same distribution as the training data. The final time is extended to $t = 0.2$ to test the generalization ability of the new GNS. The distribution of the corresponding relative errors \eqref{eq:relative_error} with respect to the Knudsen numbers for NSF, GNSS, and GNSL at $t = 0.05, 0.1, 0.15$ and $0.2$ is presented in Fig. \ref{fig:1D3Vwave}. The results indicate that there is only a small gap between the numerical results of the NSF equations and both GNS models for small Knudsen numbers. However, for large Knudsen numbers, GNSS and GNSL significantly outperform the NSF equations. Moreover, GNSL, which is obtained with subsequent training with $L_{\text{long}}$, effectively reduces the relative error, especially for $t = 0.15$ and $0.2$, where the generalizability of the GNS model is more important. The average relative errors for the NSF equations, GNSS, and GNSL are presented in Tab.~\ref{tab:wave_error}, where the decrease of the relative error is shown quantitatively.


\begin{figure}[htbp]
    \centering
    \begin{subfigure}[b]{0.3\textwidth}
        \centering
        \includegraphics[width=\textwidth]{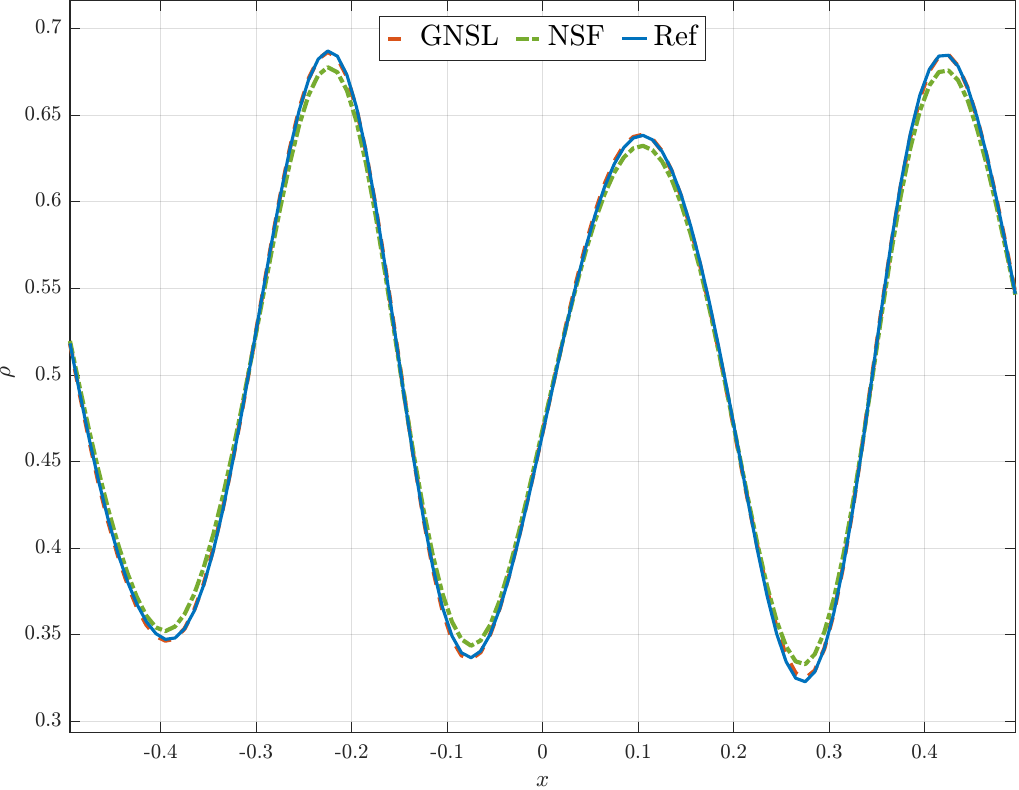}
        \caption{$\rho (\Kn = 0.0017, t = 0.05)$}
    \end{subfigure}
    \hfill
    \begin{subfigure}[b]{0.3\textwidth}
        \centering
        \includegraphics[width=\textwidth]{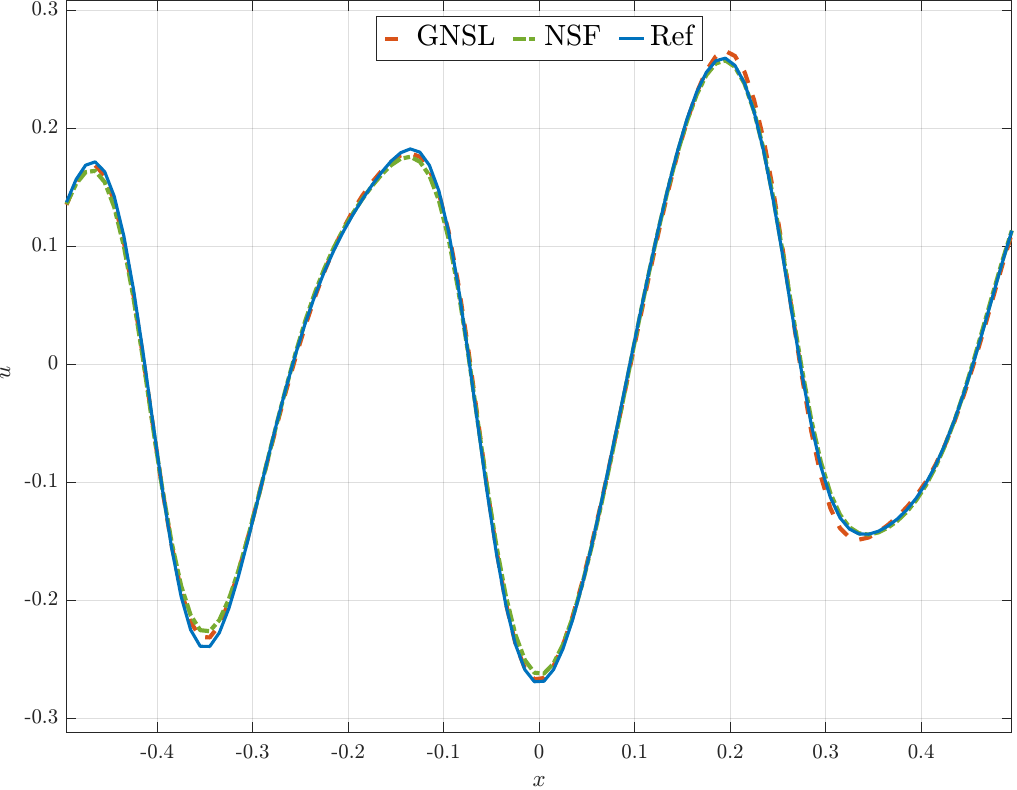}
        \caption{$u (\Kn = 0.0017, t = 0.05)$}
    \end{subfigure}
    \hfill
    \begin{subfigure}[b]{0.3\textwidth}
        \centering
        \includegraphics[width=\textwidth]{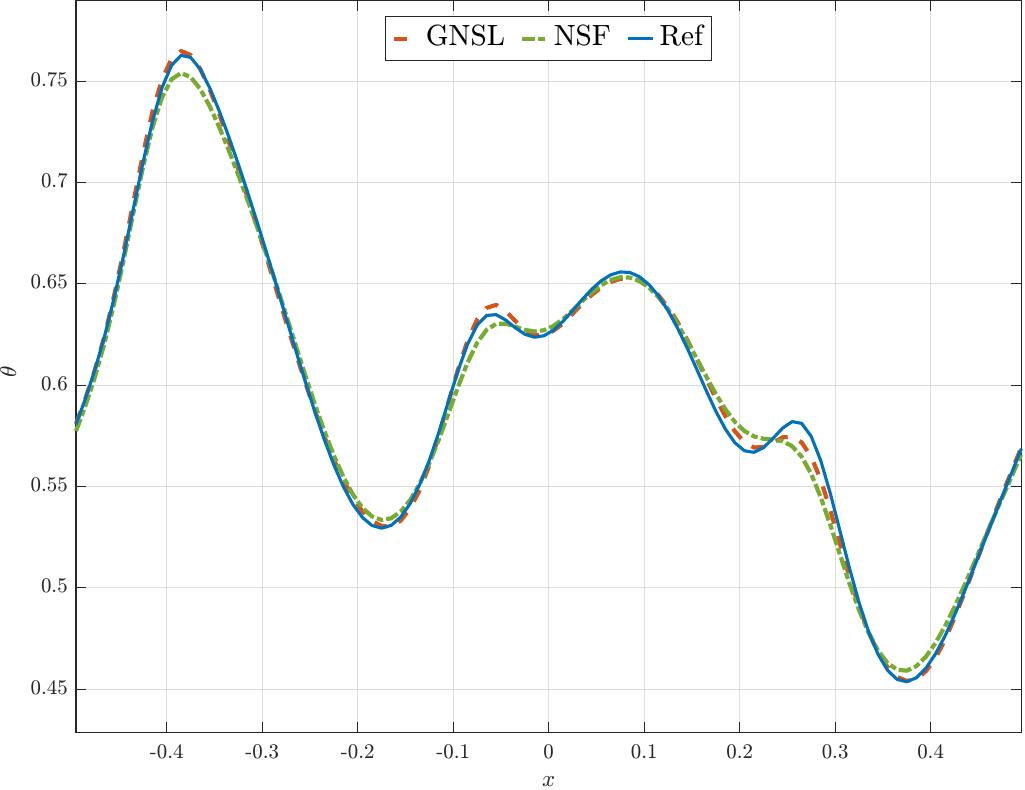}
        \caption{$\theta(\Kn = 0.0017, t = 0.05)$}
    \end{subfigure} \\
    \begin{subfigure}[b]{0.3\textwidth}
        \centering
        \includegraphics[width=\textwidth]{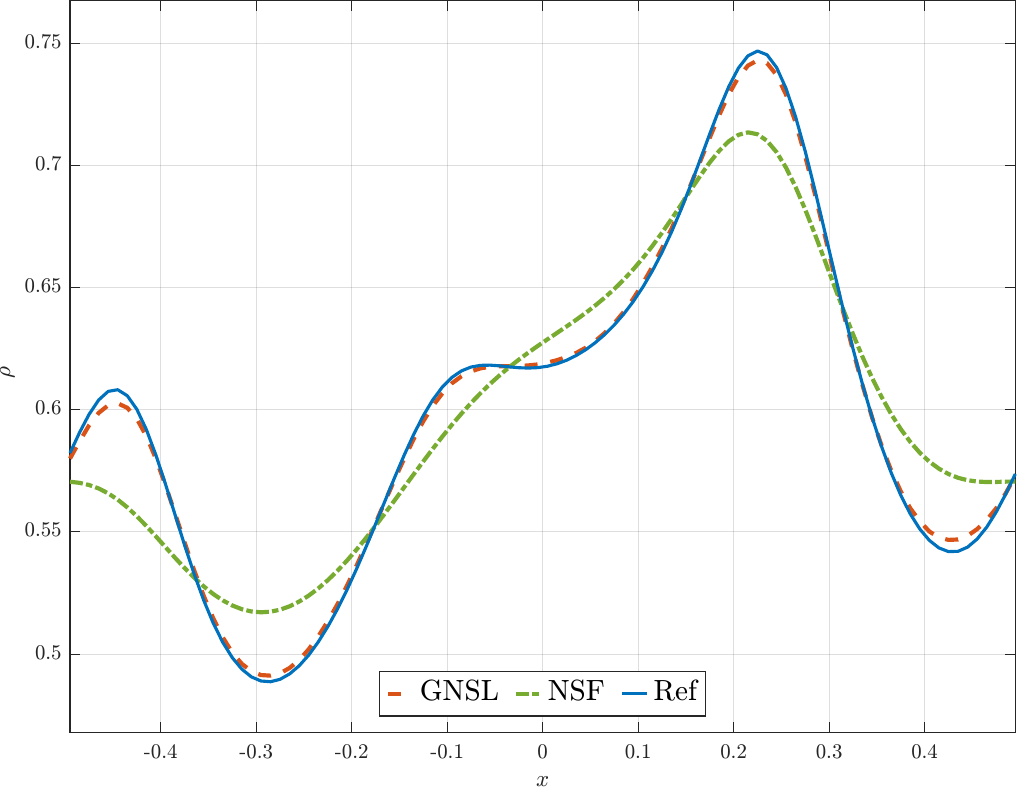}
        \caption{$\rho(\Kn = 0.0179, t = 0.1)$}
    \end{subfigure}
    \hfill
    \begin{subfigure}[b]{0.3\textwidth}
        \centering
        \includegraphics[width=\textwidth]{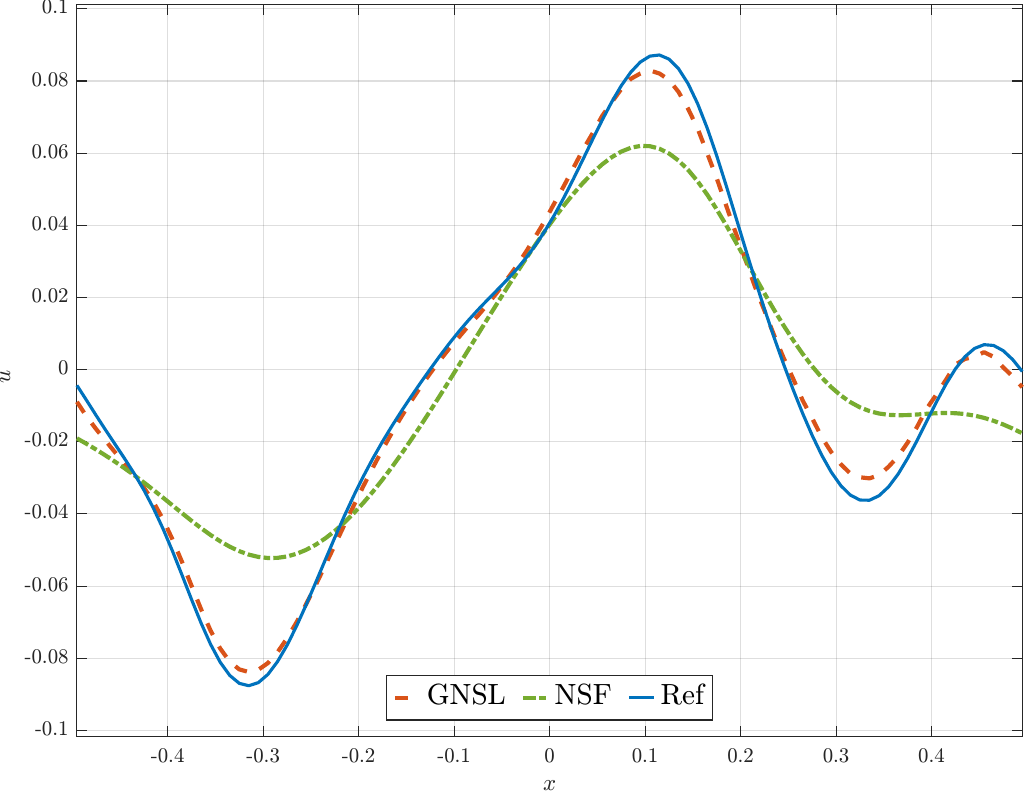}
        \caption{$u (\Kn = 0.0179, t = 0.1)$}
    \end{subfigure}
    \hfill
    \begin{subfigure}[b]{0.3\textwidth}
        \centering
        \includegraphics[width=\textwidth]{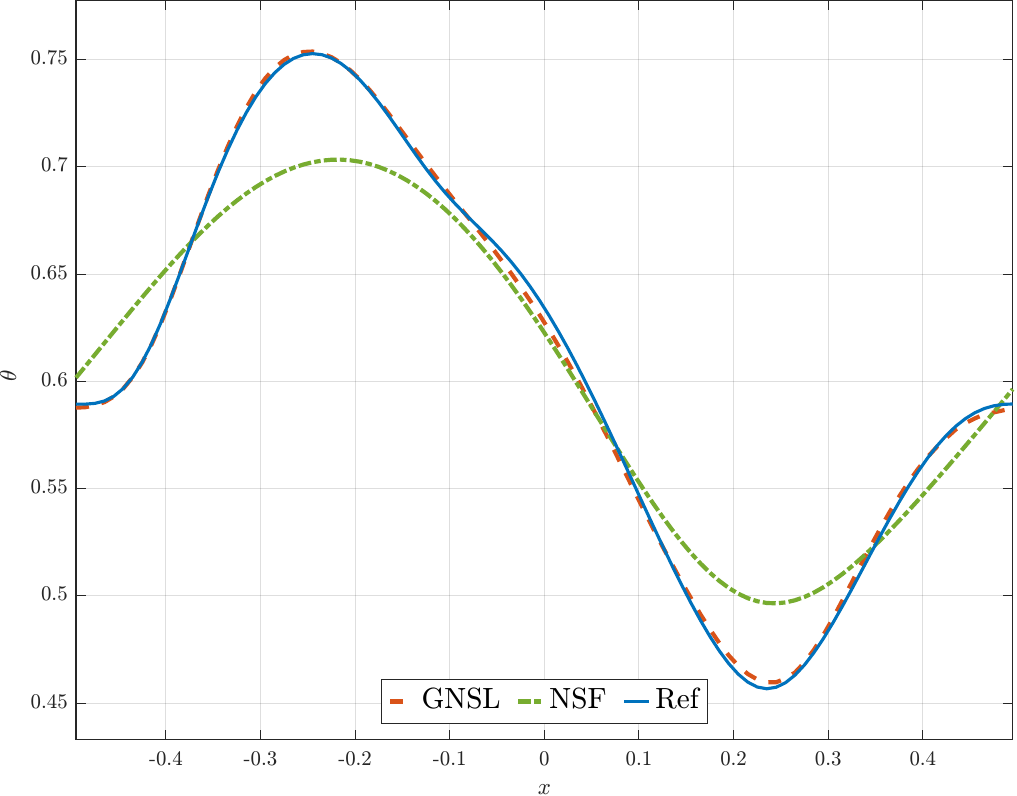}
        \caption{$\theta(\Kn = 0.0179, t = 0.1)$}
    \end{subfigure}
    \caption{(1D3V wave problem in Sec.~\ref{sec:wave}) The numerical solutions of the density $\rho$, macroscopic velocity $u$, and temperature $\theta$ obtained by GNSL, NSF, and the Boltzmann equation for $\Kn = 0.0017$ at $t = 0.05$, and $\Kn =  0.0179$ at $t = 0.1$.}
    \label{fig:wave_example_1}
\end{figure}

\begin{figure}[!htbp]
    \centering
    \begin{subfigure}[b]{0.3\textwidth}
        \centering
        \includegraphics[width=\textwidth]{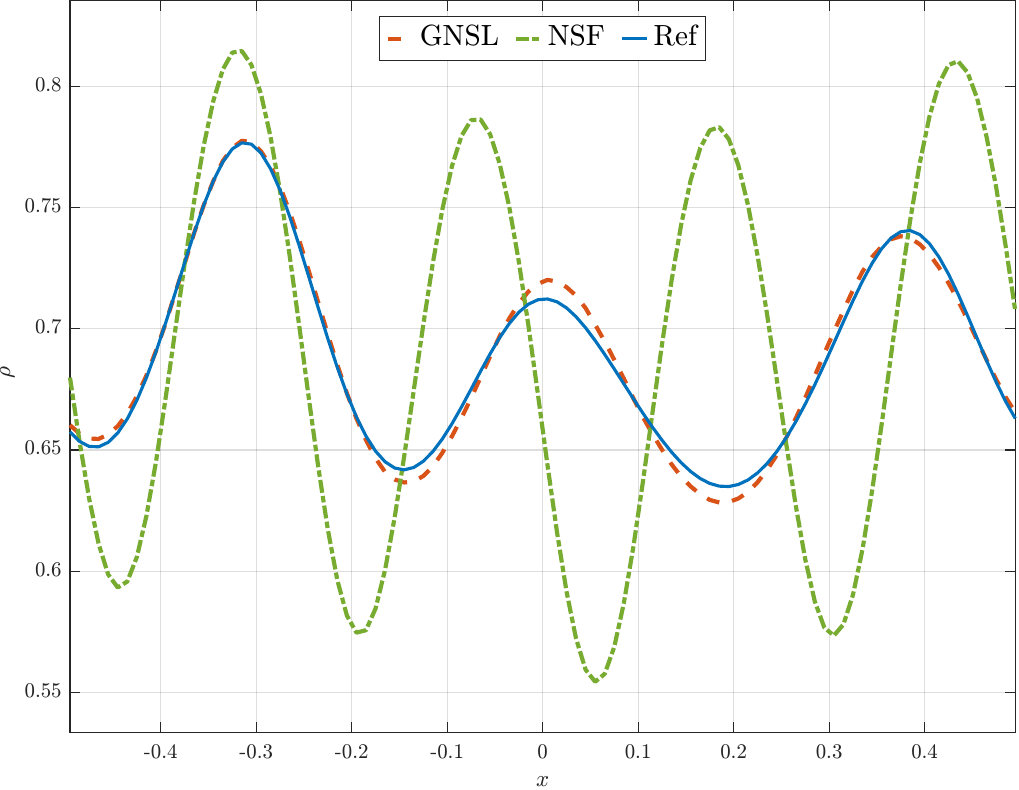}
        \caption{$\rho(\Kn = 0.1391, t = 0.15)$}
    \end{subfigure}
    \hfill
    \begin{subfigure}[b]{0.3\textwidth}
        \centering
        \includegraphics[width=\textwidth]{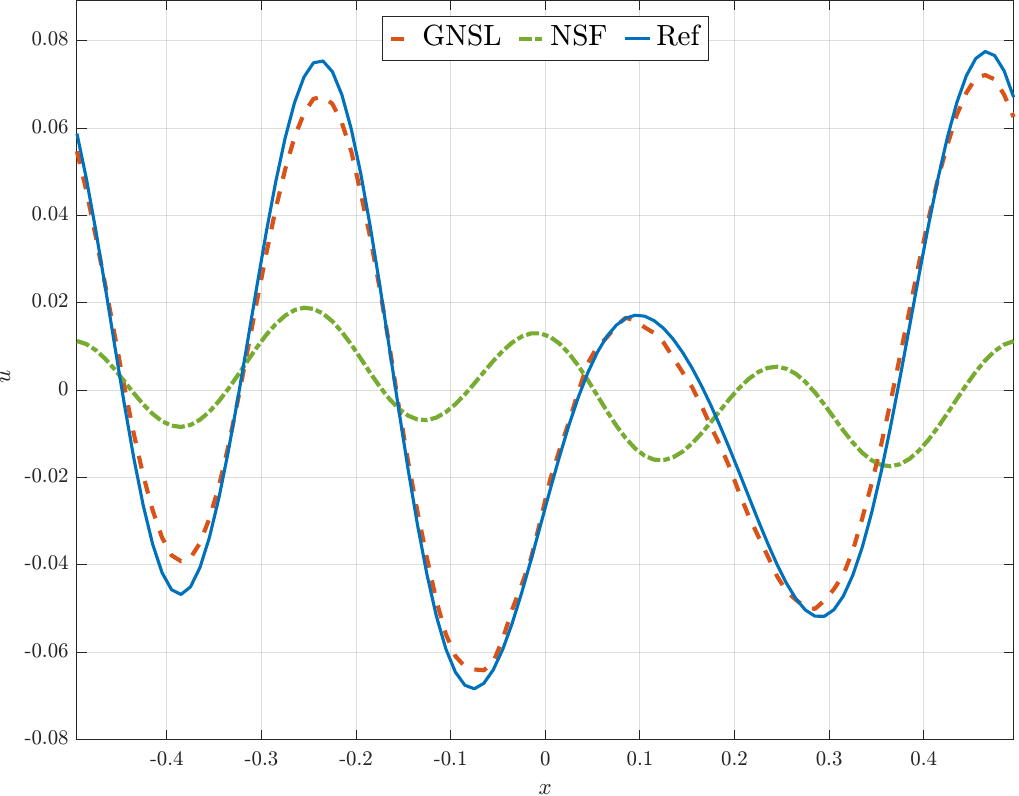}
        \caption{$u(\Kn = 0.1391, t = 0.15)$}
    \end{subfigure}
    \hfill
    \begin{subfigure}[b]{0.3\textwidth}
        \centering
        \includegraphics[width=\textwidth]{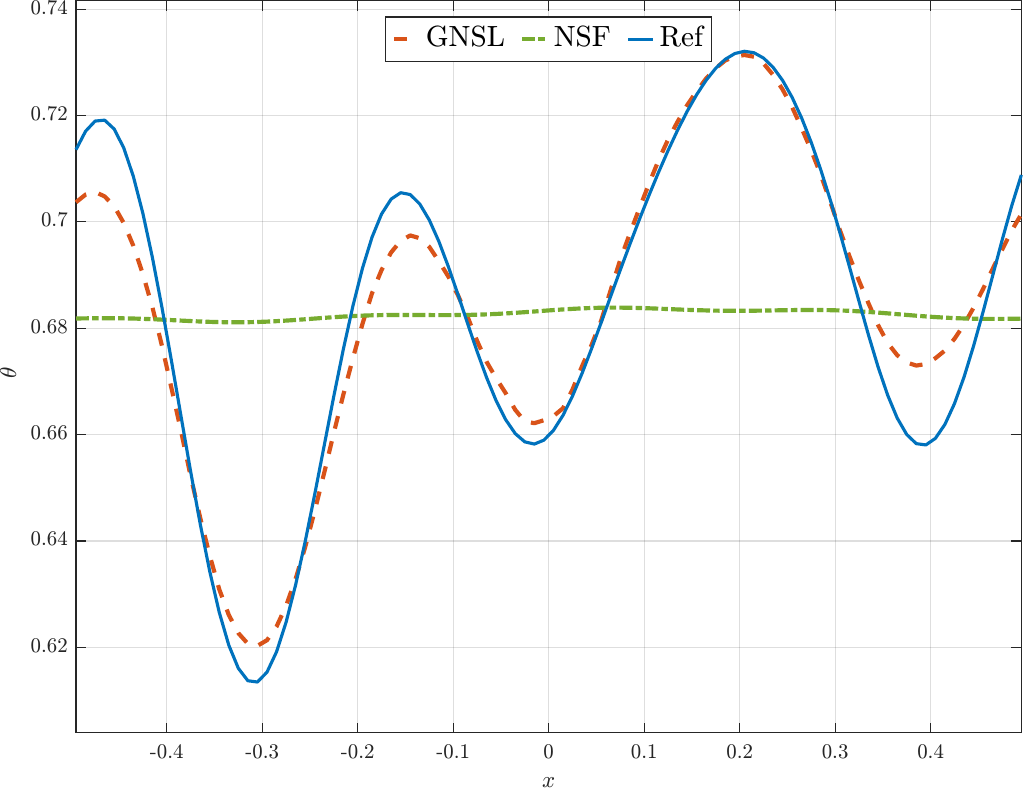}
        \caption{$\theta(\Kn = 0.1391, t = 0.15)$}
    \end{subfigure} \\
    \begin{subfigure}[b]{0.3\textwidth}
        \centering
        \includegraphics[width=\textwidth]{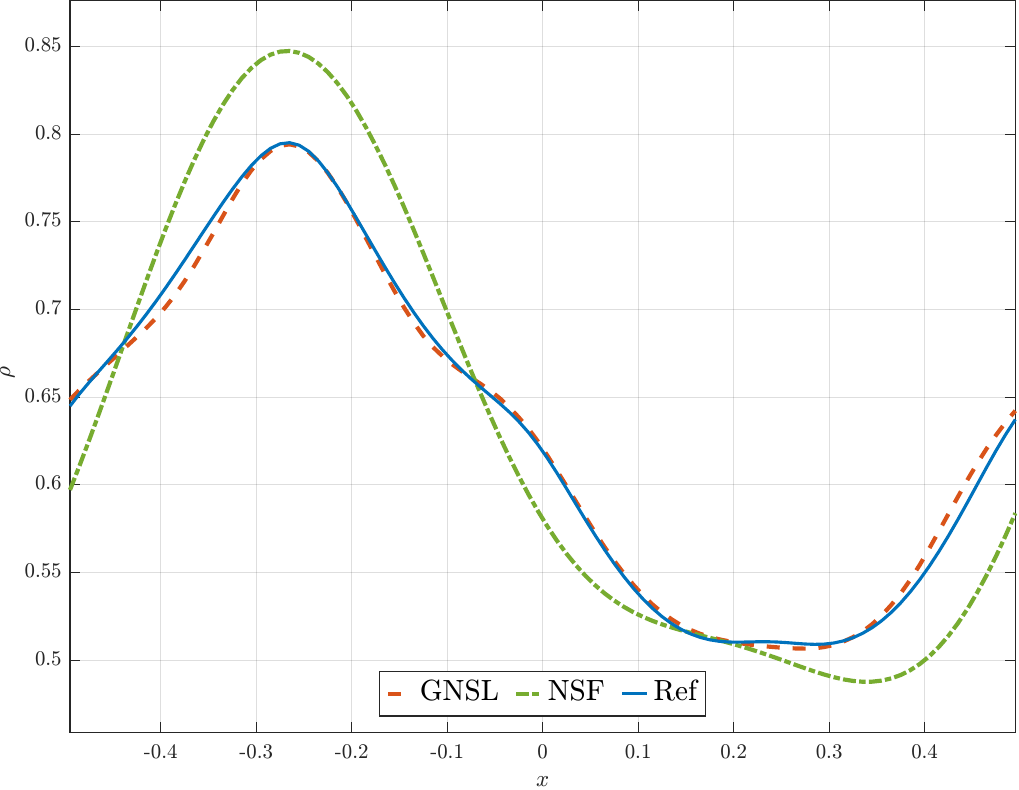}
        \caption{$\rho(\Kn = 1.2641, t = 0.2)$}
    \end{subfigure}
    \hfill
    \begin{subfigure}[b]{0.3\textwidth}
        \centering
        \includegraphics[width=\textwidth]{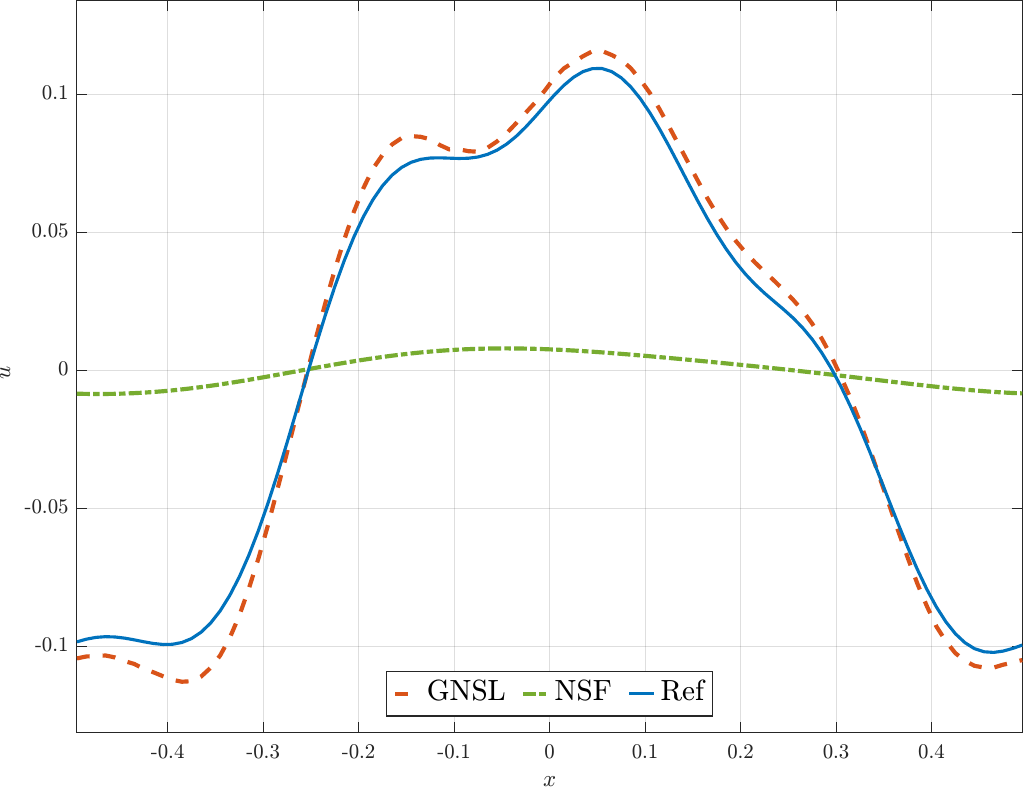}
        \caption{$u(\Kn = 1.2641, t = 0.2)$}
    \end{subfigure}
    \hfill
    \begin{subfigure}[b]{0.3\textwidth}
        \centering
        \includegraphics[width=\textwidth]{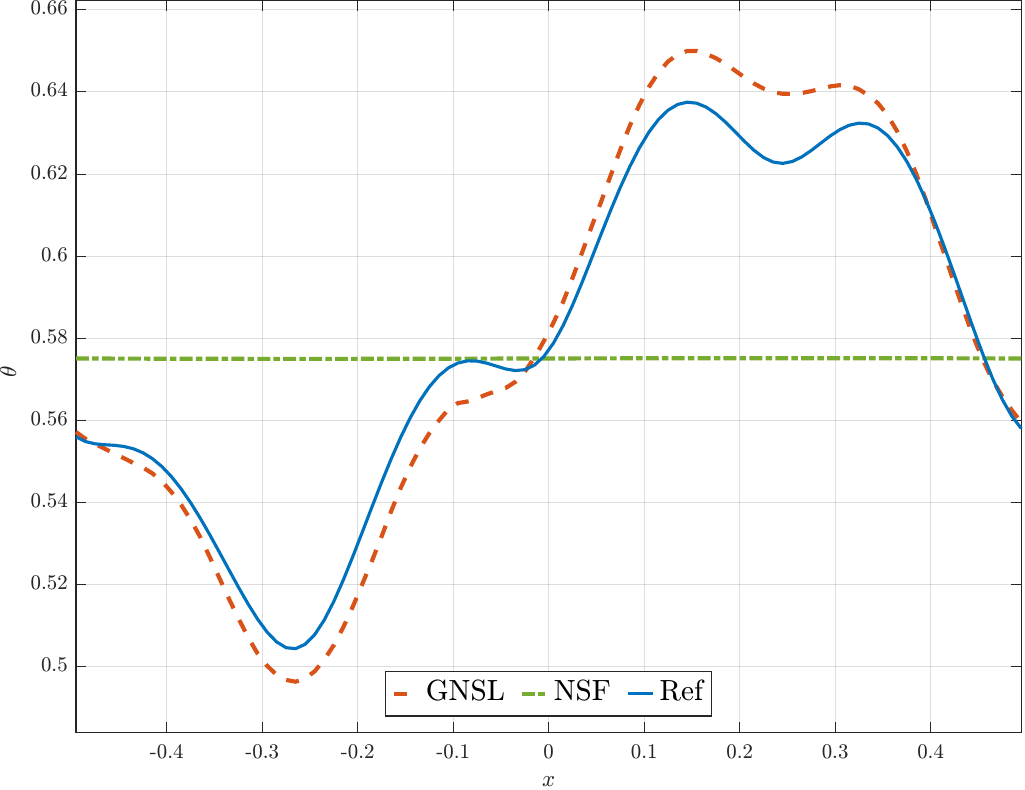}
        \caption{$\theta(\Kn = 1.2641, t = 0.2)$}
    \end{subfigure}
    \caption{(1D3V wave problem in Sec.~\ref{sec:wave}) The numerical solutions of the density $\rho$, macroscopic velocity $u$, and temperature $\theta$ obtained by GNSL, NSF, and the Boltzmann equation for $\Kn = 0.1391$ at $t = 0.15$, and $\Kn = 1.2641$ at $t = 0.2$.}    
    \label{fig:wave_example_2}
\end{figure}

To show the performance of GNSL more clearly, four specific samples with different Knudsen numbers at $t = 0.05, 0.1, 0.15$, and $0.2$ are selected, and the corresponding numerical solutions are plotted in Fig.~\ref{fig:wave_example_1} and \ref{fig:wave_example_2}. For these four samples, their relative errors $\mE_{\rm GNSL}^t$, which are listed in Tab.~\ref{tab:wave_error_1}, are all close to the average error levels, as shown in Tab.~\ref{tab:wave_error}. Fig.~\ref{fig:wave_example_1} shows that when $\Kn$ is small, the numerical solutions of GNSL and NSF all fit well with those of the Boltzmann equation. However, when $\Kn$ is increased to $\Kn = 0.0179$, there exist clear differences between NSF and Boltzmann, but those of GNSL and Boltzmann are still on top of each other. The generalization ability of GNSL is illustrated in Fig.~\ref{fig:wave_example_2}, where the numerical solutions at $t = 0.15$ and $t = 0.2$ are plotted. It shows that NSF can not capture the evolution of the macroscopic variables at such large Knudsen numbers, but the error between the numerical solutions of GNSL and the Boltzmann equation is still quite small, indicating the powerful generalization ability of GNSL at large Knudsen numbers.

\begin{table}[htbp]
    \centering
    \def\arraystrech{1.5}
   {\footnotesize
    \begin{tabular}{crrrr}
        \toprule
         & $t=0.05$ & $t=0.10$ & $t=0.15$ & $t=0.20$ \\
        \midrule
        $\Kn$ & $0.0017$& $0.0179$ & $0.1391$ & $1.2641$ \\
        \midrule
       $\mE_{\rm GNSL}^t(\%)$ & $1.31$ & $1.85$ & $2.38$ & $3.39$ \\
        \bottomrule
    \end{tabular}}
     \caption{(1D3V wave problem in Sec. \ref{sec:wave}) The detailed information of the four selected samples from the testing data.}
    \label{tab:wave_error_1}
\end{table}

\subsubsection{Validating GNSL with the discontinuous initial conditions}
\label{sec:discontinuous_test}
In this subsection, the generalization ability of the GNS model trained in Sec. \ref{sec:wave} on continuous initial conditions is tested by applying it to a problem of discontinuous initial conditions without any retraining. The spatial space is the same as in Sec. \ref{sec:wave}, and the discontinuous initial condition is defined as
\begin{equation}
\label{eq:ex2_ini}
    f(x,\bm{v},0)=\mathcal{M}[\bm{U}_{\text{shock}}](x,\bm{v}),
\end{equation}
where the macroscopic variables $\bm{U}_{\text{shock}}$ satisfy 
\begin{equation}
    \bm{U}_{\text{shock}} = \begin{cases}
        &(\rho_l,0,\theta_l),\quad x\in [-0.5,x_1) \text{ or } (x_2,0.5], \\
        &(\rho_r,0,\theta_r),\quad x\in[x_1,x_2].
    \end{cases}
\end{equation}
The parameter pairs $(\rho_l,\rho_r)$ and $(\theta_l,\theta_r)$ are sampled independently and uniformly from $[0.7,0.8]\times[0.4,0.5]$. The locations $x_1$ and $x_2$ are sampled uniformly from $[-0.3,-0.1]$ and $[0.1,0.3]$, respectively.

To validate the performance of the new GNS model across different flow regimes, five samples with Knudsen numbers $\Kn \in \{0.001, 0.01, 0.1, 1.0, 10\}$ are tested, where the initial conditions are randomly chosen from \eqref{eq:ex2_ini}. The reference solution is obtained from the Boltzmann equation with the mesh size $N_x = 100$ for the spatial space, and $120\times 40 \times 40$ for the truncated microscopic velocity space $[-10, 10]^3$.

\begin{figure}[htbp]
    \centering
    \begin{subfigure}[b]{0.3\textwidth}
        \centering
        \includegraphics[width=\textwidth]{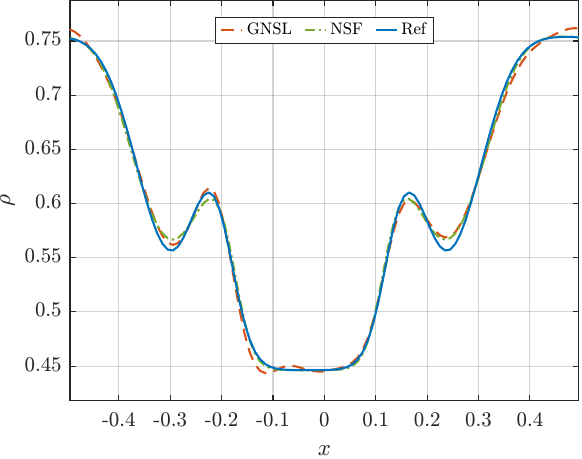}
        \caption{$\rho (\Kn = 0.001)$}
    \end{subfigure}
    \hfill
    \begin{subfigure}[b]{0.3\textwidth}
        \centering
        \includegraphics[width=\textwidth]{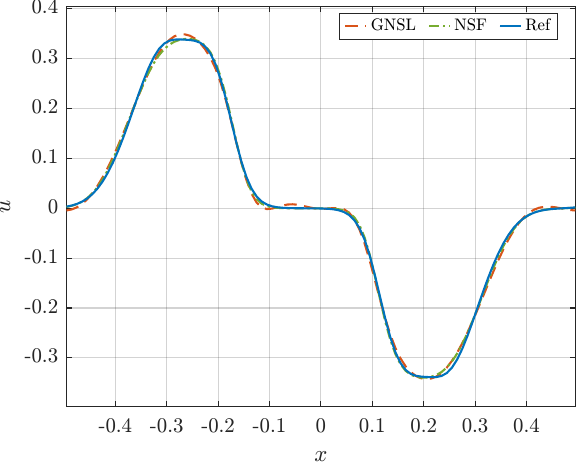}
        \caption{$u (\Kn = 0.001)$}
    \end{subfigure}
    \hfill
    \begin{subfigure}[b]{0.3\textwidth}
        \centering
        \includegraphics[width=\textwidth]{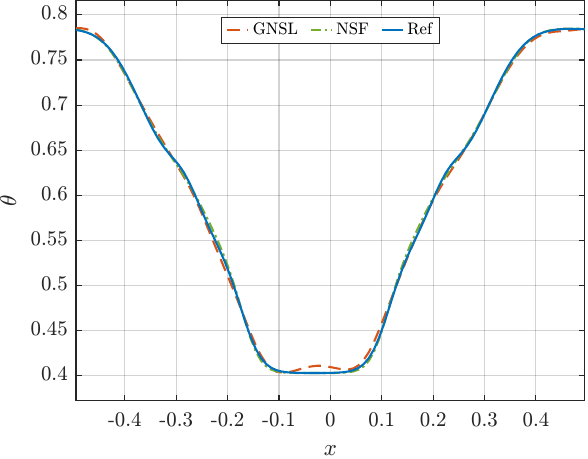}
        \caption{$\theta (\Kn = 0.001)$}
    \end{subfigure} \\
    \begin{subfigure}[b]{0.3\textwidth}
        \centering
        \includegraphics[width=\textwidth]{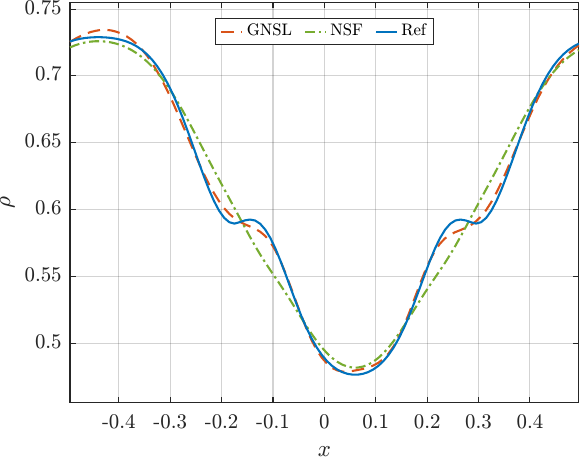}
        \caption{$\rho (\Kn = 0.01)$}
    \end{subfigure}
    \hfill
    \begin{subfigure}[b]{0.3\textwidth}
        \centering
        \includegraphics[width=\textwidth]{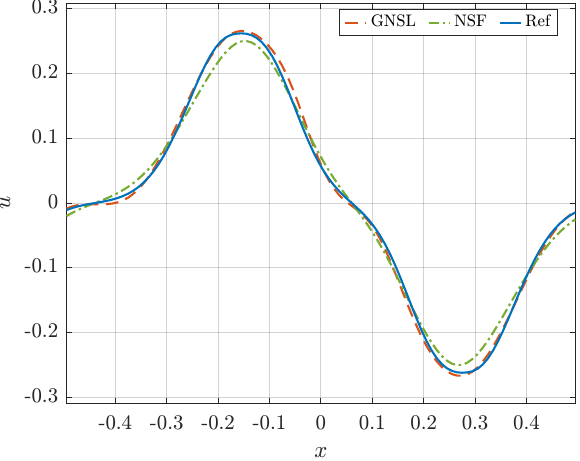}
        \caption{$u (\Kn = 0.01)$}
    \end{subfigure}
    \hfill
    \begin{subfigure}[b]{0.3\textwidth}
        \centering
        \includegraphics[width=\textwidth]{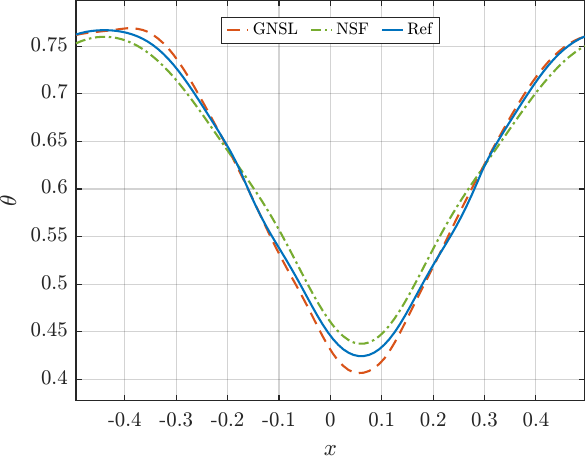}
        \caption{$\theta (\Kn = 0.01)$}
    \end{subfigure}
    \caption{(1D3V problem with discontinuous initial conditions in Sec.~\ref{sec:discontinuous_test}) The numerical solutions of the density $\rho$, macroscopic velocity $u$, and temperature $\theta$ obtained by GNSL, NSF, and the Boltzmann equation for $\Kn = 0.001$ and $0.01$ at $t = 0.1$.}
    \label{fig:discont_example_1}
\end{figure}

\begin{figure}[htbp]
    \centering
    \begin{subfigure}[b]{0.3\textwidth}
        \centering
        \includegraphics[width=\textwidth]{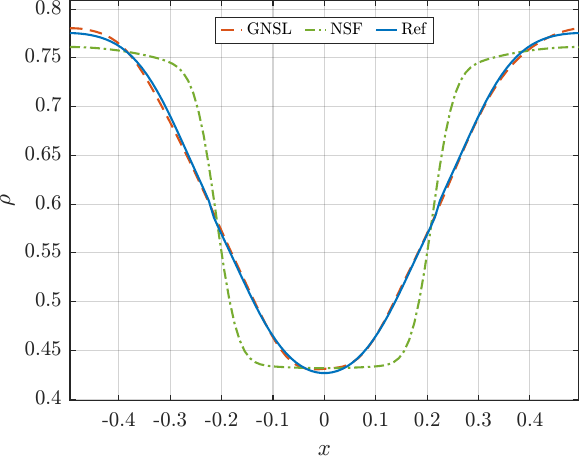}
        \caption{$\rho (\Kn = 0.1)$}
    \end{subfigure}
    \hfill
    \begin{subfigure}[b]{0.3\textwidth}
        \centering
        \includegraphics[width=\textwidth]{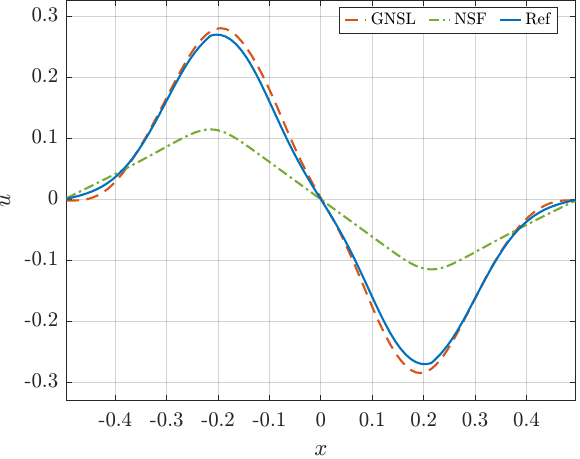}
        \caption{$u (\Kn = 0.1)$}
    \end{subfigure}
    \hfill
    \begin{subfigure}[b]{0.3\textwidth}
        \centering
        \includegraphics[width=\textwidth]{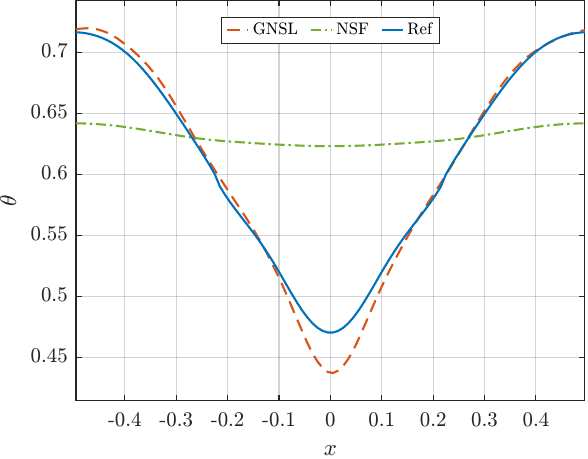}
        \caption{$\theta(\Kn = 0.1)$}
    \end{subfigure} \\
    \begin{subfigure}[b]{0.3\textwidth}
        \centering
        \includegraphics[width=\textwidth]{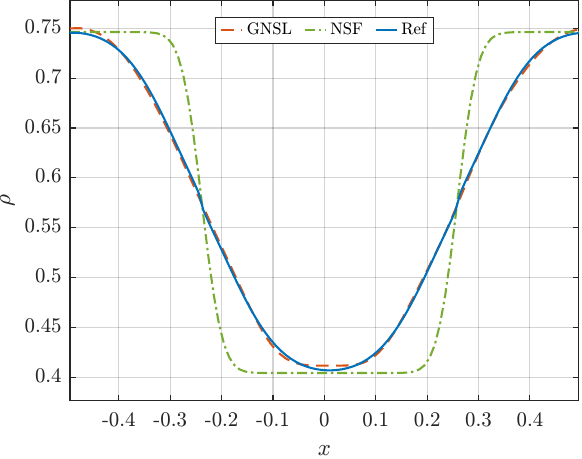}
        \caption{$\rho(\Kn = 1.0)$}
    \end{subfigure}
    \hfill
    \begin{subfigure}[b]{0.3\textwidth}
        \centering
        \includegraphics[width=\textwidth]{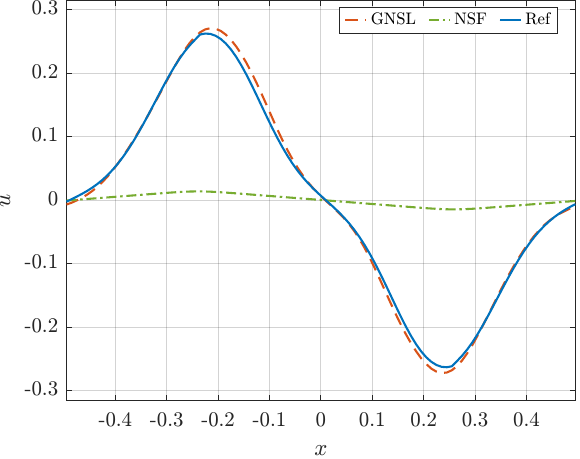}
        \caption{$u (\Kn = 1.0)$}
    \end{subfigure}
    \hfill
    \begin{subfigure}[b]{0.3\textwidth}
        \centering
        \includegraphics[width=\textwidth]{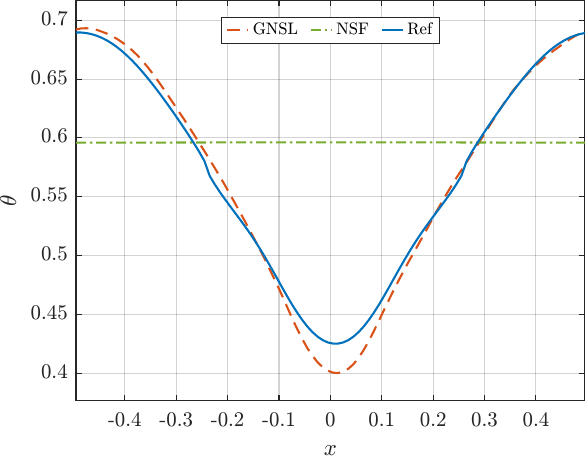}
        \caption{$\theta(\Kn = 1.0)$}
    \end{subfigure}\\
    \begin{subfigure}[b]{0.3\textwidth}
        \centering
        \includegraphics[width=\textwidth]{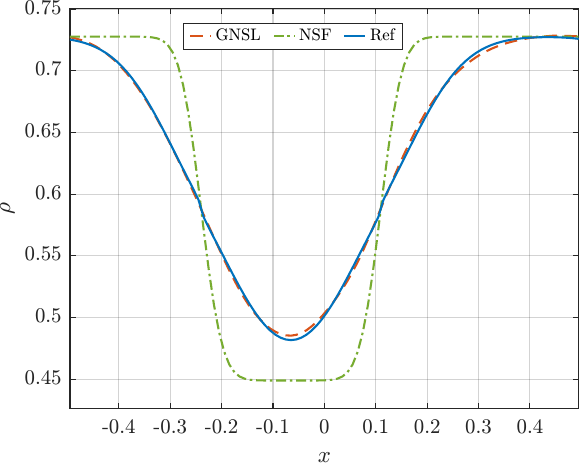}
        \caption{$\rho (\Kn = 10.0)$}
    \end{subfigure}
    \hfill
    \begin{subfigure}[b]{0.3\textwidth}
        \centering
        \includegraphics[width=\textwidth]{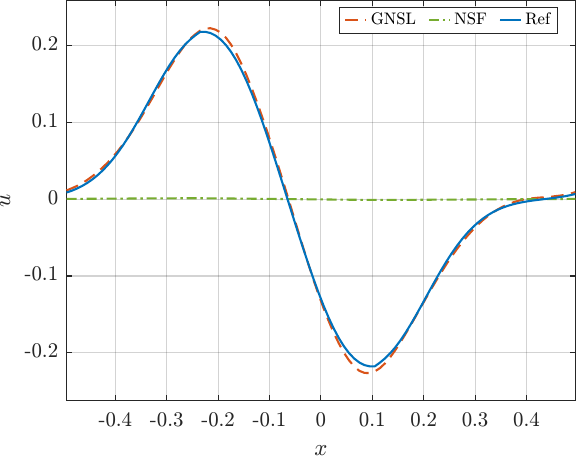}
        \caption{$u(\Kn = 10.0)$}
    \end{subfigure}
    \hfill
    \begin{subfigure}[b]{0.3\textwidth}
        \centering
        \includegraphics[width=\textwidth]{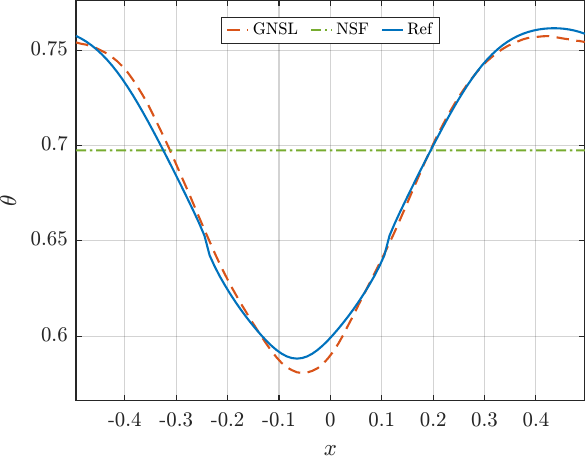}
        \caption{$\theta(\Kn = 10.0)$}
    \end{subfigure}
    \caption{(1D3V problem with discontinuous initial conditions in Sec.~\ref{sec:discontinuous_test}) The numerical solutions of the density $\rho$, macroscopic velocity $u$, and temperature $\theta$ obtained by GNSL, NSF, and the Boltzmann equation for $\Kn = 0.1, 1.0$ and $10.0$ at $t = 0.1$.}
    \label{fig:discont_example_2}
\end{figure}

The numerical solutions of density $\rho$, macroscopic velocity $u$, and temperature $\theta$ obtained by GNSL, NSF, and the Boltzmann equation at $t = 0.1$ are plotted in Fig.~\ref{fig:discont_example_1} and \ref{fig:discont_example_2}. It shows that in the near continuum regime $\Kn = 0.001$ and $0.01$, the numerical solutions for both GNSL and NSF are close to the reference solution as in Fig. \ref{fig:discont_example_1}. However, in the more rarefied regimes where $\Kn = 0.1, 1$ and $10$, the numerical solution of NSF deviates clearly, while that of GNSL still preserves reasonable accuracy. To show the different behavior of GNSL and NSF quantitatively, the relative errors $\mE_{\rm s}^t$ for GNSL and NSF are listed in Tab.~\ref{tab:wave_applyto_discont_error}. It clearly illustrates that for NSF, $\mE_{\rm NSF}^{t}$ is increasing quickly with increasing $\Kn$, while that of GNSL remains small for all $\Kn$, indicating that the GNS model trained on smooth initial conditions can effectively generalize to discontinuous initial conditions.


\begin{table}[htbp]
    \centering
    \def\arraystrech{1.5}
    {\footnotesize
    \begin{tabular}{lrrrrr}
        \toprule
        \diagbox[width=6em]{model}{Kn}
        & $0.001$ & $0.01$ & $0.1$ & $1.0$ & $10.0$ \\
        \midrule
        NSF
        & ${1.32}$ & $4.42$ & $28.83$ & $46.99$ & $45.76$ \\
        \midrule
        GNSL
        & $1.94$ & ${2.35}$ & ${3.30}$ & ${2.72}$ & ${1.70}$ \\
        \bottomrule
    \end{tabular}}
    \caption{(1D3V problem with discontinuous initial conditions in Sec.~\ref{sec:discontinuous_test}) The relative errors $\mE_{\rm s}^{t}$ (in percentages) of the problem with sampled discontinuous initial conditions for $\Kn=0.001$, $0.01$, $0.1$, $1.0$, and $10$ at $t = 0.1$.}
    \label{tab:wave_applyto_discont_error}
\end{table}

\subsection{1D3V Riemann problem}
\label{sec:Riemann}
In this section, GNS is validated on the Riemann problem, which is defined on the physical domain $[-1,1]$ with a zero gradient boundary condition. The initial condition is a Maxwellian $\mathcal{M}[\bm{U}]$, where the discontinuous macroscopic variables $\bm{U}=(\rho,u,\theta)$ are given by
\begin{equation}\label{eq:Riemann_init1}
    (\rho,u,\theta) = \begin{cases}
        &(\rho_L,u_L,\theta_L),\quad x\leqslant 0, \\
        &(\rho_R,u_R,\theta_R),\quad x>0,
    \end{cases}
\end{equation}
with
\begin{equation}\label{eq:Riemann_init2}
\begin{gathered}
 \rho_L(x)=1+a_L, \quad u_L(x)=0, \quad \theta_L(x)=1+b_L, \\
 \rho_R(x)=0.2+a_R,\quad u_R(x)=0,\quad \theta_R(x)=0.4+b_R,
\end{gathered}
\end{equation}
where the parameters $a_L,a_R$, and $b_L, b_R$ are uniformly and independently sampled from $[0,0.5]$.

A setting similar to that in Sec.~\ref{sec:wave_training} is utilized to generate the training data with the spatial mesh size $N_x = 400$, and the microscopic velocity space is truncated to $[-10,10]^3$ with a mesh size of $160\times 40\times 40$. The Knudsen number is sampled from a log-uniform distribution over $[0.001, 10]$. $500$ samples based on the initial condition \eqref{eq:Riemann_init1} and the sampled Knudsen numbers are generated, and they are evolved to $t = 0.1$ with the time step length $\Delta t = 0.00025$. The training data set consists of all macroscopic variables, including density $\rho$, macroscopic velocity $u$, and temperature $\theta$, at each time step.

\begin{figure}[htbp]
    \centering
    \begin{subfigure}[b]{0.24\textwidth}
        \centering
        \includegraphics[width=\linewidth]{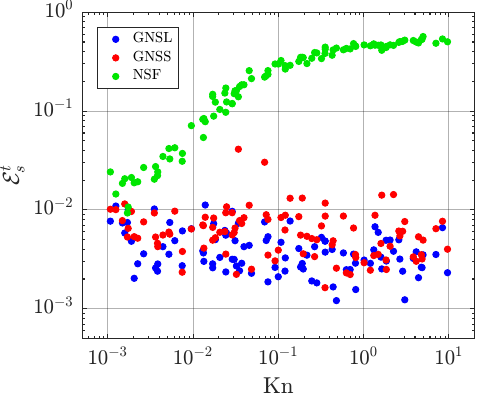}
        \caption{$t = 0.05$}
    \end{subfigure}
    \begin{subfigure}[b]{0.24\textwidth}
        \centering
        \includegraphics[width=\linewidth]{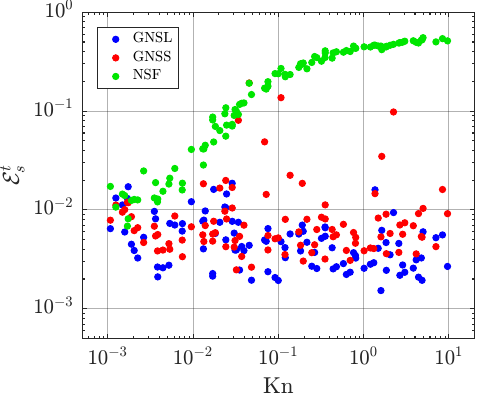}
        \caption{$t = 0.10$}
    \end{subfigure} 
    \begin{subfigure}[b]{0.24\textwidth}
        \centering
        \includegraphics[width=\linewidth]{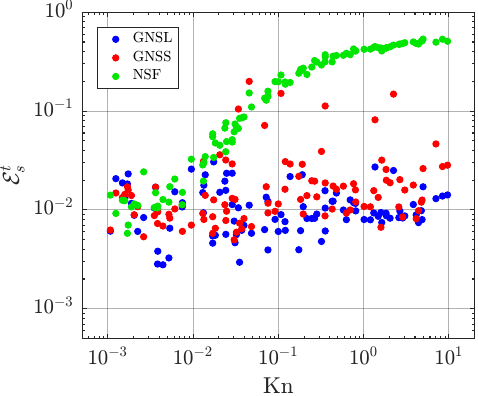}
        \caption{ $t = 0.15$}
    \end{subfigure}
    \begin{subfigure}[b]{0.24\textwidth}
        \centering
        \includegraphics[width=\linewidth]{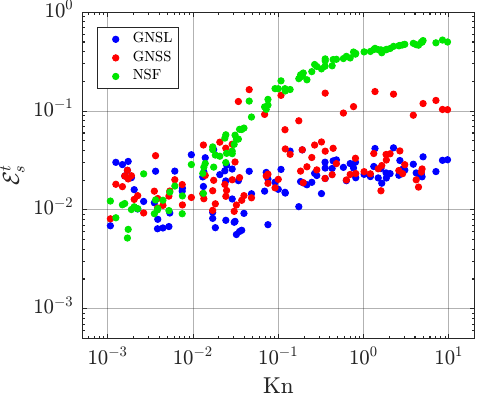}
        \caption{$t = 0.20$}
    \end{subfigure}    
    \caption{(1D3V Riemann problem in Sec.~\ref{sec:Riemann}) The distribution of the relative error $\mathcal{E}_{\rm s}^t$ (log-log form) with respect to the Knudsen number for NSF, GNSS, and GNSL at $t=0.05, 0.1, 0.15$ and $0.2$.}
    \label{fig:1D3VRiemann_error}
\end{figure}

The testing data are generated from $100$ samples with random initial conditions and Knudsen numbers following the same distribution as the training data. The final simulation time is extended to $t = 0.2$ to further test the generalization ability of GNS. The distribution of the relative error $\mE_{\rm s}^t$ for GNSL, GNSS, and NSF at $t = 0.05, 0.1, 0.15$, and $0.2$ is plotted in Fig.~\ref{fig:1D3VRiemann_error}. It shows that at all times, the relative error of NSF is much larger than those of GNSS and GNSL, especially for large Knudsen numbers, while the relative error of GNSL is a little smaller than that of GNSS. This indicates that GNSL has better behavior and better generalization ability than GNSS. This property is also validated quantitatively by Tab. ~\ref{tab:Riemann_error}, where the average relative errors for NSF, GNSS, and GNSL on the testing data are presented.




\begin{table}[htbp]
    \centering
    \def\arraystrech{1.5}
    {\footnotesize
    \begin{tabular}{lrrrr}
        \toprule
        \diagbox[width=6em]{model}{time}
        & $t = 0.05$ & $t = 0.10$ & $t = 0.15$ & $t = 0.20$ \\
        \midrule
        NSF  & $26.58$ & $23.78$ & $21.74$ & $20.12$ \\
        \midrule
        GNSS & $0.68$ & $1.27$ & $2.25$ & $3.86$ \\
        \midrule
        GNSL & $\bm{0.42}$ & $\bm{0.55}$ & $\bm{1.11}$ & $\bm{2.12}$ \\
        \bottomrule
    \end{tabular}}
    \caption{(1D3V Riemann problem in Sec.~\ref{sec:Riemann}) The average relative errors $\mE_{\rm s}^{t}$ (in percentages) of the testing data for different models at $t=0.05$, $0.1$, $0.15$, and $0.2$.}
    \label{tab:Riemann_error}
\end{table}

    
    
    

\begin{figure}[htbp]
    \centering
    \begin{subfigure}[b]{0.3\textwidth}
        \centering
        \includegraphics[width=\textwidth]{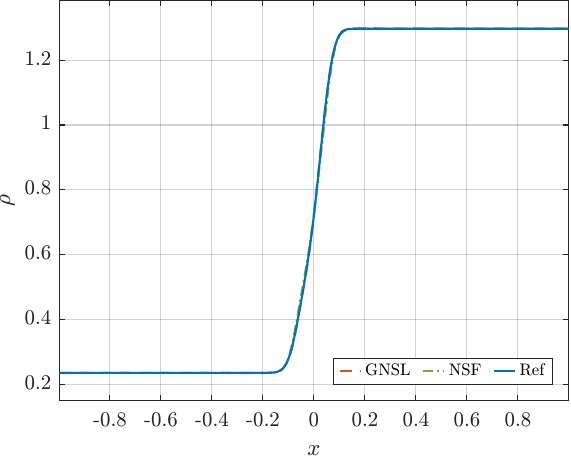}
        \caption{$\rho(\Kn = 0.0019, t = 0.05)$}
    \end{subfigure}
    \hfill
    \begin{subfigure}[b]{0.3\textwidth}
        \centering
        \includegraphics[width=\textwidth]{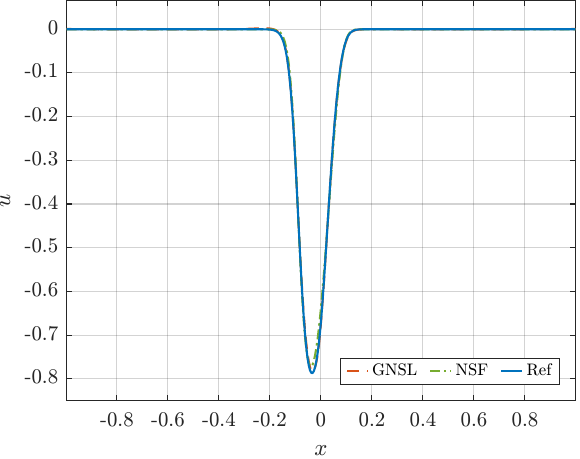}
        \caption{$u(\Kn = 0.0019, t = 0.05)$}
    \end{subfigure}
    \hfill
    \begin{subfigure}[b]{0.3\textwidth}
        \centering
        \includegraphics[width=\textwidth]{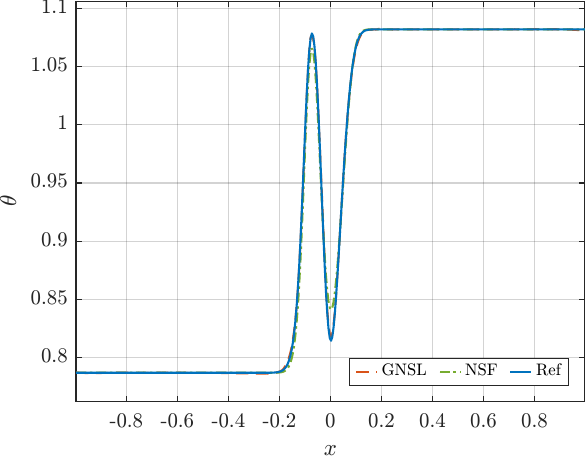}
        \caption{$\theta(\Kn = 0.0019, t = 0.05)$}
    \end{subfigure} \\
    \begin{subfigure}[b]{0.3\textwidth}
        \centering
        \includegraphics[width=\textwidth]{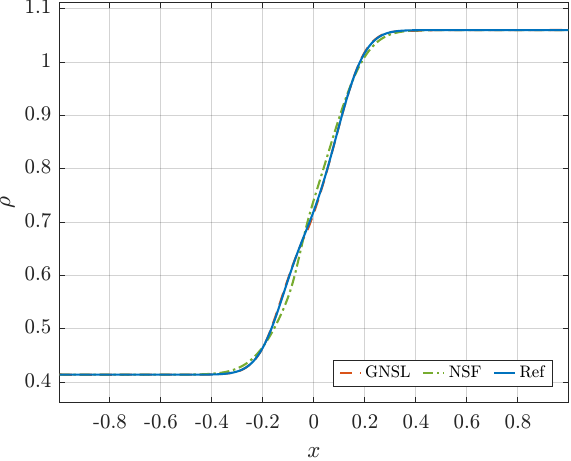}
        \caption{$\rho(\Kn = 0.0182, t = 0.1)$}
    \end{subfigure}
    \hfill
    \begin{subfigure}[b]{0.3\textwidth}
        \centering
        \includegraphics[width=\textwidth]{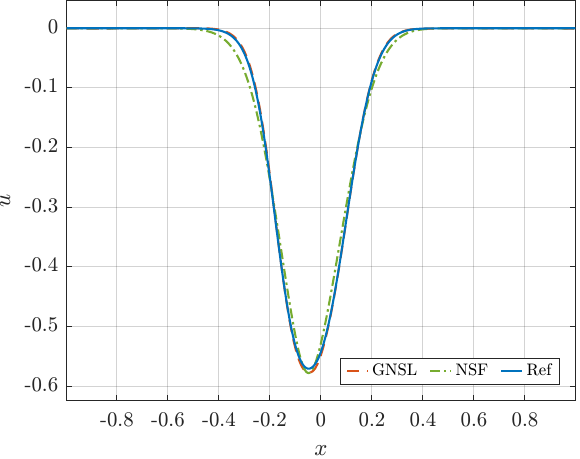}
        \caption{$u(\Kn = 0.0182, t = 0.1)$}
    \end{subfigure}
    \hfill
    \begin{subfigure}[b]{0.3\textwidth}
        \centering
        \includegraphics[width=\textwidth]{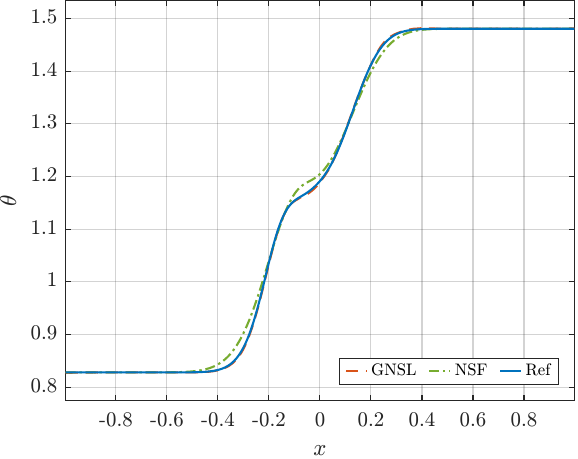}
        \caption{$\theta(\Kn = 0.0182, t = 0.1)$}
    \end{subfigure}
    \caption{(1D3V Riemann problem in Sec.~\ref{sec:Riemann}) The numerical solutions of the density $\rho$, macroscopic velocity $u$, and temperature $\theta$ obtained by GNSL, NSF, and the Boltzmann equation for $\Kn = 0.0019$ at $t = 0.05$, and $\Kn =  0.0182$ at $t = 0.1$.}
    \label{fig:Riemann_example_1}
\end{figure}

\begin{figure}[!htbp]
    \centering
    \begin{subfigure}[b]{0.3\textwidth}
        \centering
        \includegraphics[width=\textwidth]{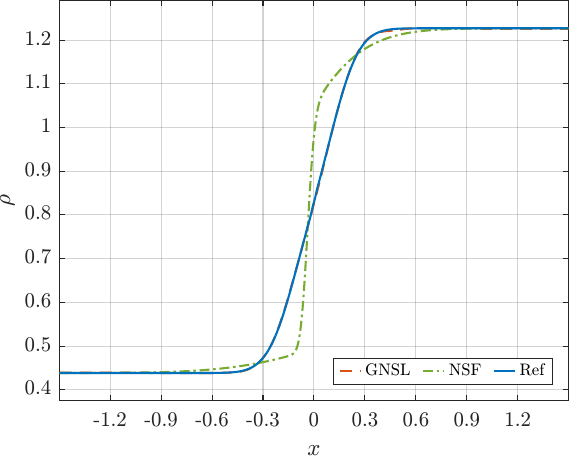}
        \caption{$\rho(\Kn = 0.1957,  t = 0.15)$}
    \end{subfigure}
    \hfill
    \begin{subfigure}[b]{0.3\textwidth}
        \centering
        \includegraphics[width=\textwidth]{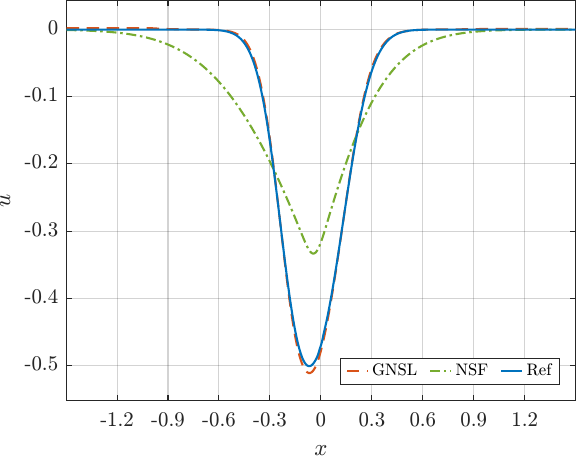}
        \caption{$u(\Kn = 0.1957, t = 0.15)$}
    \end{subfigure}
    \hfill
    \begin{subfigure}[b]{0.3\textwidth}
        \centering
        \includegraphics[width=\textwidth]{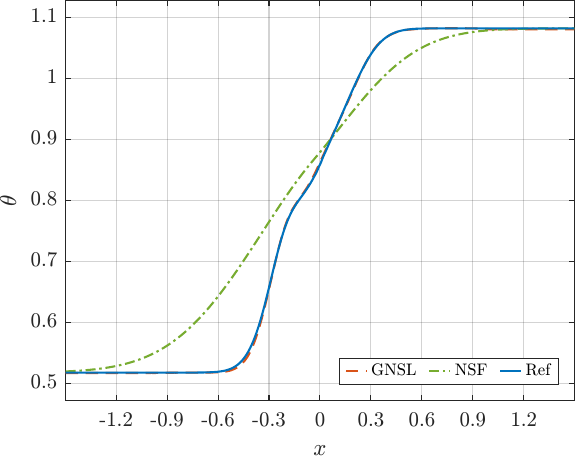}
        \caption{$\theta(\Kn = 0.1957, t = 0.15)$}
    \end{subfigure} \\
    \begin{subfigure}[b]{0.3\textwidth}
        \centering
        \includegraphics[width=\textwidth]{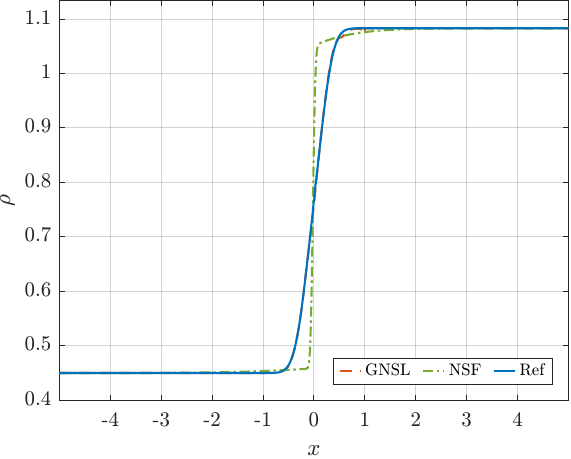}
        \caption{$\rho(\Kn = 1.4821, t = 0.2)$}
    \end{subfigure}
    \hfill
    \begin{subfigure}[b]{0.3\textwidth}
        \centering
        \includegraphics[width=\textwidth]{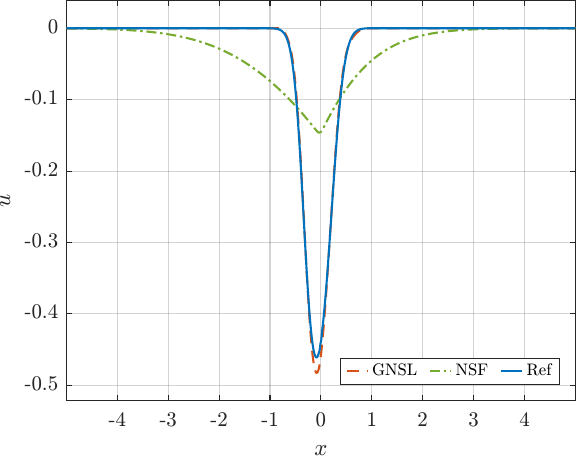}
        \caption{$u(\Kn = 1.4821, t = 0.2)$}
    \end{subfigure}
    \hfill
    \begin{subfigure}[b]{0.3\textwidth}
        \centering
        \includegraphics[width=\textwidth]{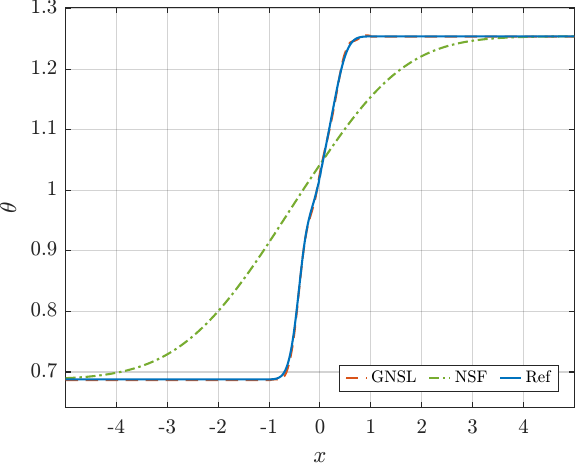}
        \caption{$\theta(\Kn = 1.4821, t = 0.2)$}
    \end{subfigure}
    \caption{(1D3V Riemann problem in Sec.~\ref{sec:Riemann}) The numerical solutions of the density $\rho$, macroscopic velocity $u$, and temperature $\theta$ obtained by GNSL, NSF, and the Boltzmann equation for $\Kn = 0.1957$ at $t = 0.15$, and $\Kn =  1.4821$ at $t = 0.2$.}    
    \label{fig:Riemann_example_2}
\end{figure}

Four specific samples with different Knudsen numbers and times are selected to show the performance of GNSL more clearly. The detailed information of these four samples is listed in Tab. \ref{tab:Riemann_error_1}, and their relative errors $\mE_{\rm GNSL}^t$ are all close to the average error levels at each time. The numerical solutions of density $\rho$, macroscopic velocity $u$, and temperature $\theta$ are plotted in Fig. \ref{fig:Riemann_example_1} and Fig. \ref{fig:Riemann_example_2}. The results show that when the Knudsen number is small, the numerical solutions of GNSL, NSF and Boltzmann match well with each other. However, with increasing $\Kn$, the NSF solution deviates from the Boltzmann solution, whereas the GNSL solution still agrees well with it, indicating the broad applicability of GNSL.

\begin{table}[htbp]
    \centering
    \def\arraystrech{1.5}
   {\footnotesize
    \begin{tabular}{crrrr}
        \toprule
         & $t=0.05$ & $t=0.10$ & $t=0.15$ & $t=0.20$ \\
        \midrule
        $\Kn$ & $0.0019$& $0.0182$ & $0.1957$ & $1.4821$ \\
        \midrule
       $\mE_{\rm GNSL}^t(\%)$ & $0.48$ & $0.58$ & $1.08$ & $2.06$ \\
        \bottomrule
    \end{tabular}}
     \caption{(1D3V Riemann problem in Sec. \ref{sec:Riemann}) The detailed information of the four selected samples from the testing data.}
    \label{tab:Riemann_error_1}
\end{table}

\subsection{2D3V isentropic vortex problem}
\label{sec:ist}
In this section, the GNS model is validated on the 2D3V isentropic vortex problem. The spatial domain is set to $[0,10]\times[0,10]$ with periodic boundary conditions, and the initial condition is given by a Maxwellian distribution $\mathcal{M}[\bm{U}]$, where the macroscopic variables $\bm{U}=(\rho,u_1,u_2,\theta)$ are defined as
\begin{equation}
\label{eq:isentropic_ini}
\begin{gathered}
\theta=\theta_{\infty}-\frac{(\gamma-1)b^2\sigma^2}{8\pi^2\gamma}\exp\left(1-\frac{r^2}{\sigma^2}\right), \qquad \rho=\rho_{\infty}\left(\frac{\theta}{\theta_{\infty}}\right)^{\frac{1}{\gamma-1}},\\
u_1=u_{1,\infty}-\frac{b(y-5)}{2\pi}\exp\left[\frac{1}{2}\left(1-\frac{r^2}{\sigma^2}\right)\right], \qquad
u_2=u_{2,\infty}+\frac{b(x-5)}{2\pi}\exp\left[\frac{1}{2}\left(1-\frac{r^2}{\sigma^2}\right)\right].
\end{gathered}
\end{equation}
Here, $\gamma=5/3$ and $r=\sqrt{(x-5)^2+(y-5)^2}$. The parameters $b\in[6.0,8.0]$, $\sigma^2\in[0.9,1.2]$, $\theta_{\infty}\in[1.3,1.7]$, and $\rho_{\infty}\in[0.8,1.2]$ are sampled independently and uniformly, while $u_{1,\infty}$ and $u_{2,\infty}$ are uniformly sampled from $[0.5,1.0]$.

The discretization of the spatial space is $N_x = N_y = 64$, and the microscopic velocity space is truncated to $[-7.5,7.5]^3$, with the mesh size $32\times 32\times 32$. The Knudsen number is uniformly sampled from the logarithmic range of $[0.001, 10]$. For the training data, the number of samples is $100$, and they are generated from the initial condition \eqref{eq:isentropic_ini} and the sampled Knudsen numbers. The final simulation time is $t = 0.2$ with the time step length $\Delta t = 0.0005$. The training data set consists of all macroscopic variables, including density $\rho$, macroscopic velocity $\bm{u}$, and temperature $\theta$ at each time step. 

\begin{figure}[htbp]
    \centering
    \begin{subfigure}[b]{0.24\textwidth}
        \centering
        \includegraphics[width=\linewidth]{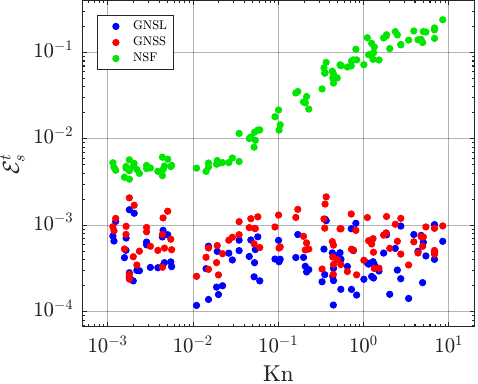}
        \caption{$t = 0.1$}
    \end{subfigure}
    \begin{subfigure}[b]{0.24\textwidth}
        \centering
        \includegraphics[width=\linewidth]{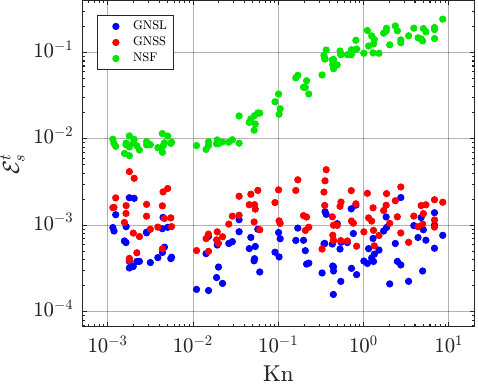}
        \caption{$t = 0.2$}
    \end{subfigure} 
    \begin{subfigure}[b]{0.24\textwidth}
        \centering
        \includegraphics[width=\linewidth]{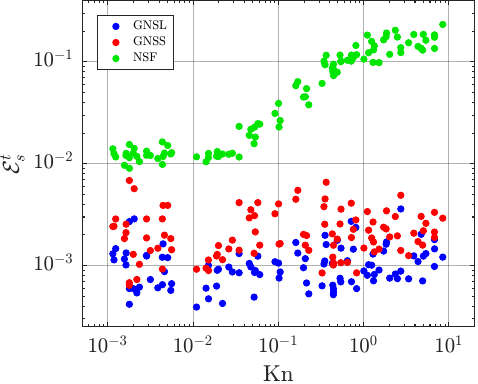}
        \caption{ $t = 0.3$}
    \end{subfigure}
    \begin{subfigure}[b]{0.24\textwidth}
        \centering
        \includegraphics[width=\linewidth]{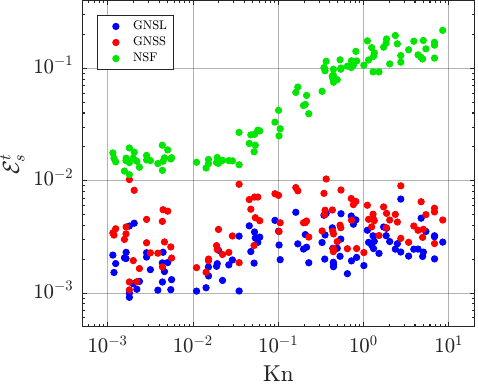}
        \caption{ $t = 0.4$}
    \end{subfigure}
    \caption{(2D3V isentropic vortex problem in Sec.~\ref{sec:ist}) The distribution of the relative error $\mathcal{E}_{\rm s}^t$ (log-log form) with respect to the Knudsen number for NSF, GNSS, and GNSL at $t=0.1, 0.2$, $0.3$ and $0.4$.}
    \label{fig:2D3Vist}
\end{figure}

For the testing data, the number of samples is also 100. These samples are generated with the random initial condition \eqref{eq:isentropic_ini} and Knudsen numbers drawn from the same distribution as that of the training data. The final simulation time is extended to $t = 0.4$ to further test the generalization ability of GNS. The average relative errors $\mE_{\rm s}^t$ for NSF, GNSS, and GNSL at $t = 0.1, 0.2, 0.3$, and $0.4$ are shown in Tab. \ref{tab:ist_error}. It shows that the relative error of NSF is much larger than those of GNSS and GNSL, and GNSL gives the smallest relative error among the three methods. The distribution of the relative error $\mE_{\rm s}^t$ with respect to the Knudsen number is shown in Fig. \ref{fig:2D3Vist}. We can clearly see that $\mE_{\rm NSF}^t$ is much larger than that of GNSS and GNSL for all $\Kn$. Moreover, $\mE_{\rm GNSL}^t$ is noticeably smaller than $\mE_{\rm GNSS}^t$, indicating the efficiency of GNSL. These results confirm that the neural network closure model effectively extends to two-dimensional scenarios. It maintains low relative errors across the tested Knudsen numbers, whereas the NSF equations exhibit much larger deviations. Additionally, the subsequent training with $L_{\text{long}}$ reduces the prediction error substantially compared to the model trained only with $L_{\text{short}}$.




\begin{table}[htbp]
    \centering
    \def\arraystrech{1.5}
    {\footnotesize
    \begin{tabular}{lrrrr}
        \toprule
        \diagbox[width=6em]{model}{time}
        & $t = 0.1$ & $t = 0.2$ & $t = 0.3$  & $t = 0.4$\\
        \midrule
        NSF  & $5.64$ & $6.88$ & $7.18$  & $7.08$\\
        \midrule
        GNSS & $0.08$ & $0.14$ & $0.23$  &  $0.43$\\
        \midrule
        GNSL & $\bm{0.05}$ & $\bm{0.07}$ & $\bm{0.11}$ & $\bm{0.26}$\\
        \bottomrule
    \end{tabular}}
    \caption{(2D3V isentropic vortex problem in Sec. \ref{sec:ist}) The average relative errors $\mE_{\rm s}^{t}$ (in percentages) of the testing data for different models at $t=0.1$, $0.2$, $0.3$ and $0.4$.}
    \label{tab:ist_error}
\end{table}

    
    
    

\begin{figure}[!htbp]
    \centering
    \begin{subfigure}[b]{0.24\textwidth}
        \centering
        \includegraphics[width=\textwidth]{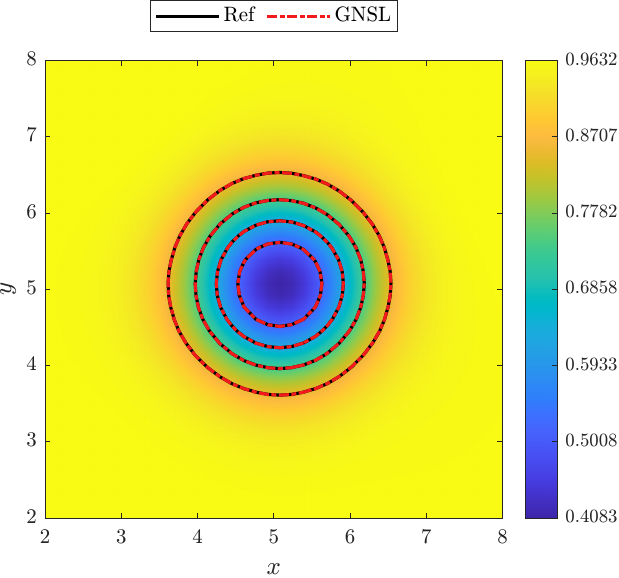}
        \caption{$\rho$}
    \end{subfigure}
    \hfill
    \begin{subfigure}[b]{0.24\textwidth}
        \centering
        \includegraphics[width=\textwidth]{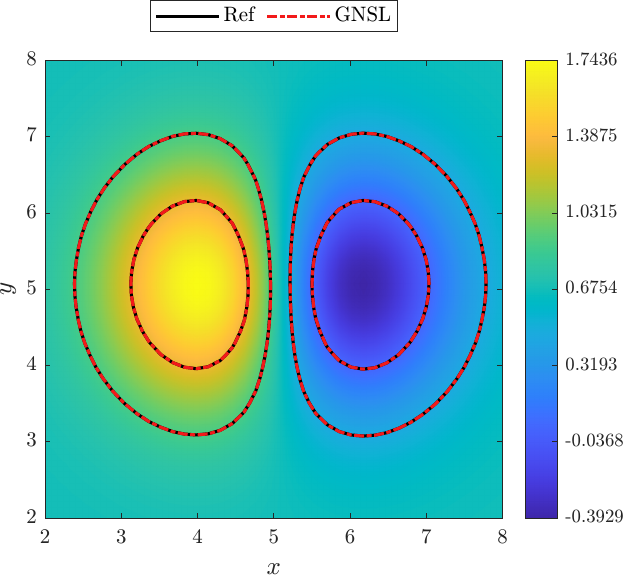}
        \caption{$u_1$}
    \end{subfigure}
    \hfill
    \begin{subfigure}[b]{0.24\textwidth}
        \centering
        \includegraphics[width=\textwidth]{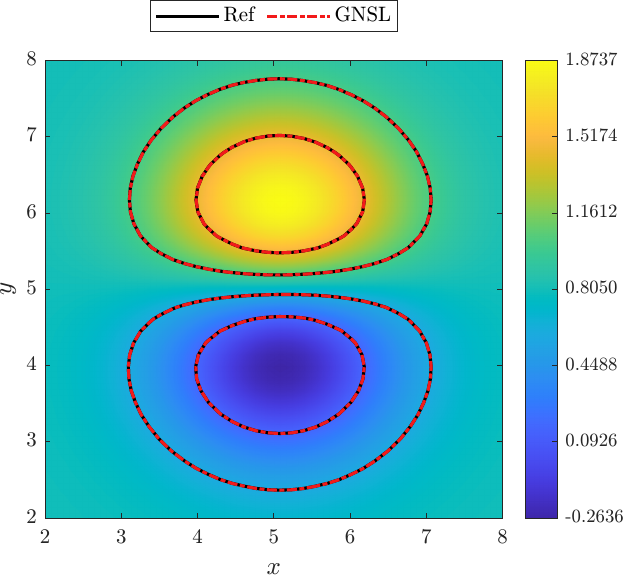}
        \caption{$u_2$}
    \end{subfigure} 
    \begin{subfigure}[b]{0.24\textwidth}
        \centering
        \includegraphics[width=\textwidth]{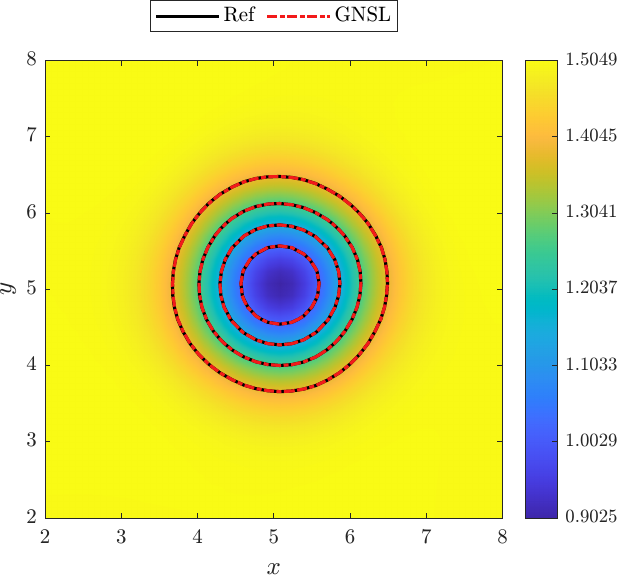}
        \caption{$\theta$}
    \end{subfigure}
    \hfill
    \caption{(2D3V isentropic vortex problem in Sec.~\ref{sec:ist})  The numerical solutions of the density $\rho$, macroscopic velocity $u_1, u_2$, and temperature $\theta$ obtained by GNSL and the Boltzmann equation for $\Kn = 0.5288$ at $t = 0.1$.}    
    \label{fig:vortex_example_1}
\end{figure}

\begin{figure}[!htbp]
    \centering
    \begin{subfigure}[b]{0.24\textwidth}
        \centering
        \includegraphics[width=\textwidth]{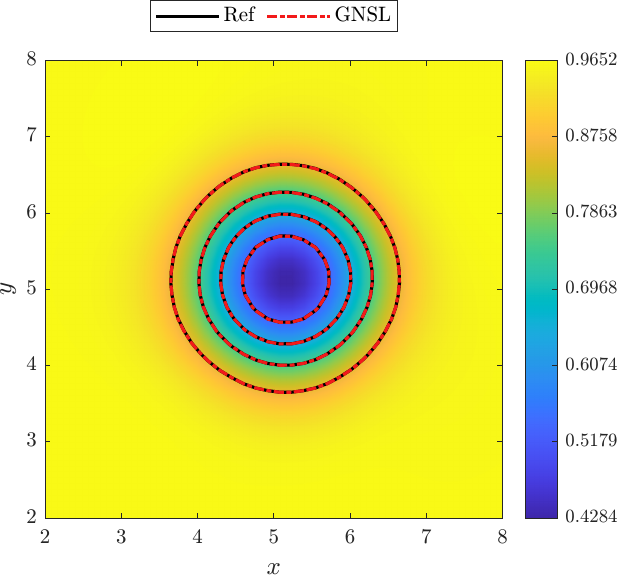}
        \caption{$\rho$}
    \end{subfigure}
    \hfill
    \begin{subfigure}[b]{0.24\textwidth}
        \centering
        \includegraphics[width=\textwidth]{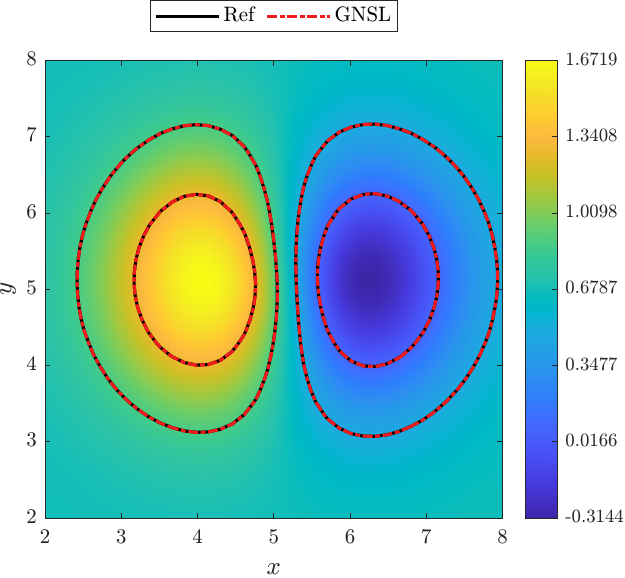}
        \caption{$u_1$}
    \end{subfigure}
    \hfill
    \begin{subfigure}[b]{0.24\textwidth}
        \centering
        \includegraphics[width=\textwidth]{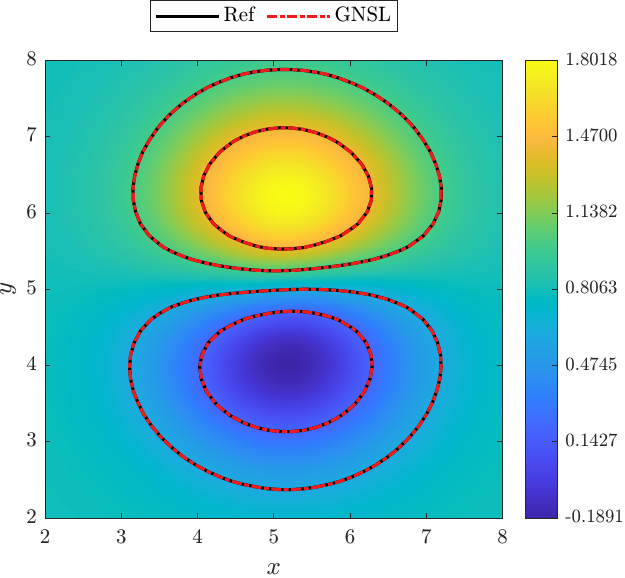}
        \caption{$u_2$}
    \end{subfigure} 
    \begin{subfigure}[b]{0.24\textwidth}
        \centering
        \includegraphics[width=\textwidth]{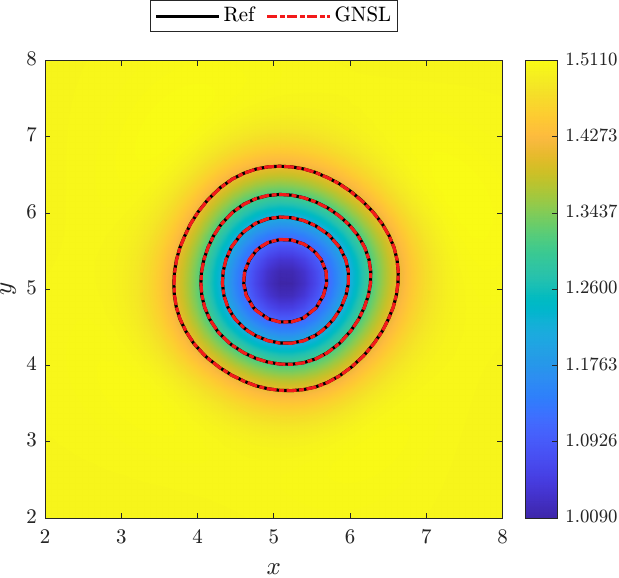}
        \caption{$\theta$}
    \end{subfigure}
    \hfill
    \caption{(2D3V isentropic vortex problem in Sec.~\ref{sec:ist})  The numerical solutions of the density $\rho$, macroscopic velocity $u_1,u_2$, and temperature $\theta$ obtained by GNSL and the Boltzmann equation for $\Kn = 5.3690$ at $t = 0.2$.}    
    \label{fig:vortex_example_2}
\end{figure}

\begin{figure}[!htbp]
    \centering
    \begin{subfigure}[b]{0.24\textwidth}
        \centering
        \includegraphics[width=\textwidth]{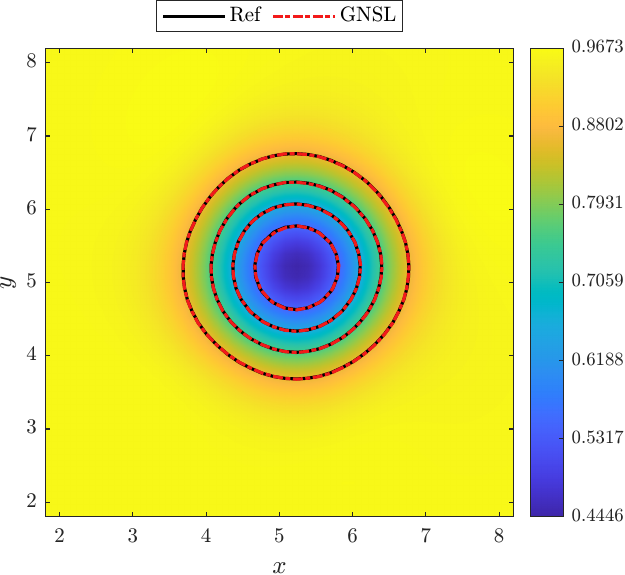}
        \caption{$\rho$}
    \end{subfigure}
    \hfill
    \begin{subfigure}[b]{0.24\textwidth}
        \centering
        \includegraphics[width=\textwidth]{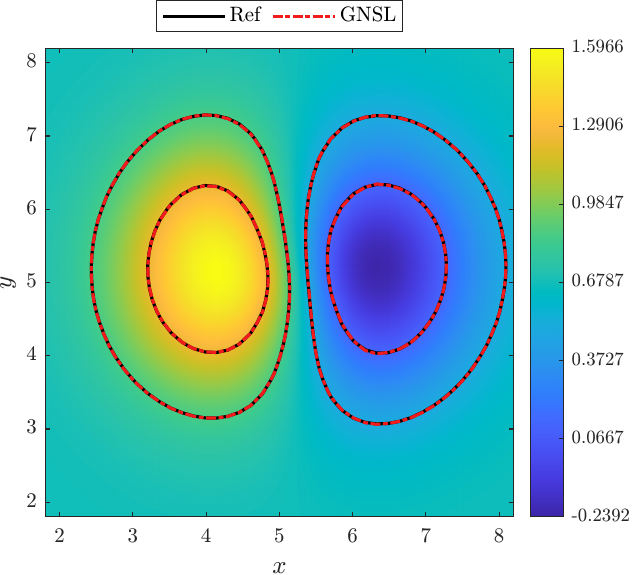}
        \caption{$u_1$}
    \end{subfigure}
    \hfill
    \begin{subfigure}[b]{0.24\textwidth}
        \centering
        \includegraphics[width=\textwidth]{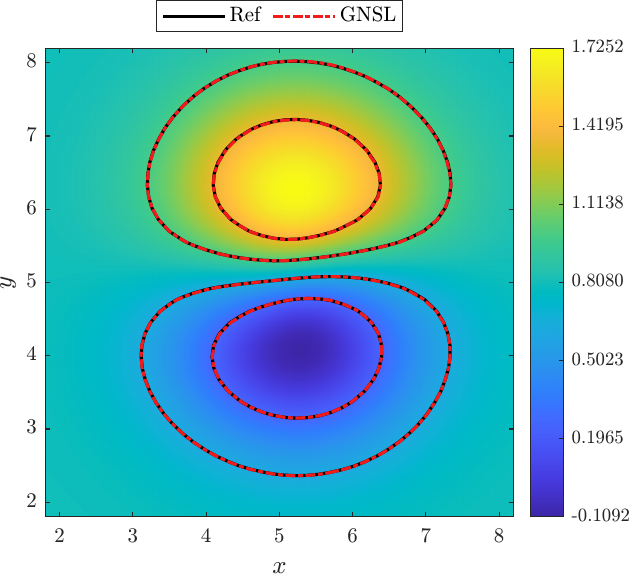}
        \caption{$u_2$}
    \end{subfigure} 
    \begin{subfigure}[b]{0.24\textwidth}
        \centering
        \includegraphics[width=\textwidth]{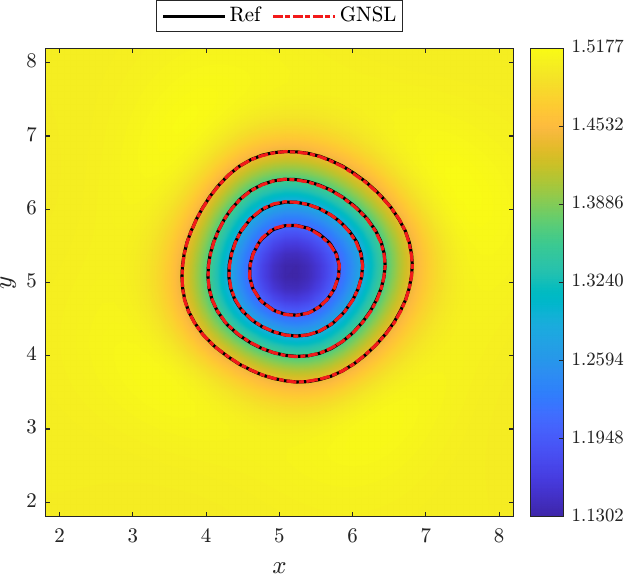}
        \caption{$\theta$}
    \end{subfigure}
    \hfill
    \caption{(2D3V isentropic vortex problem in Sec.~\ref{sec:ist})  The numerical solutions of the density $\rho$, macroscopic velocity $u_1,u_2$, and temperature $\theta$ obtained by GNSL and the Boltzmann equation for $\Kn = 0.0913$ at $t = 0.3$.}    
    \label{fig:vortex_example_3}
\end{figure}

\begin{figure}[!htbp]
    \centering
    \begin{subfigure}[b]{0.24\textwidth}
        \centering
        \includegraphics[width=\textwidth]{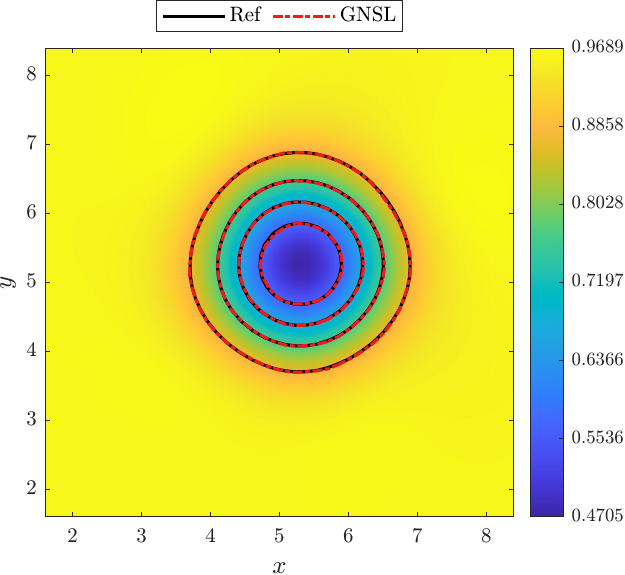}
        \caption{$\rho$}
    \end{subfigure}
    \hfill
    \begin{subfigure}[b]{0.24\textwidth}
        \centering
        \includegraphics[width=\textwidth]{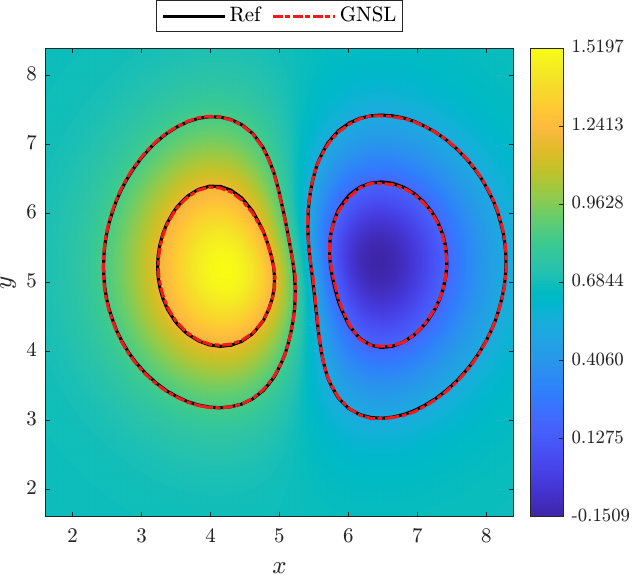}
        \caption{$u_1$}
    \end{subfigure}
    \hfill
    \begin{subfigure}[b]{0.24\textwidth}
        \centering
        \includegraphics[width=\textwidth]{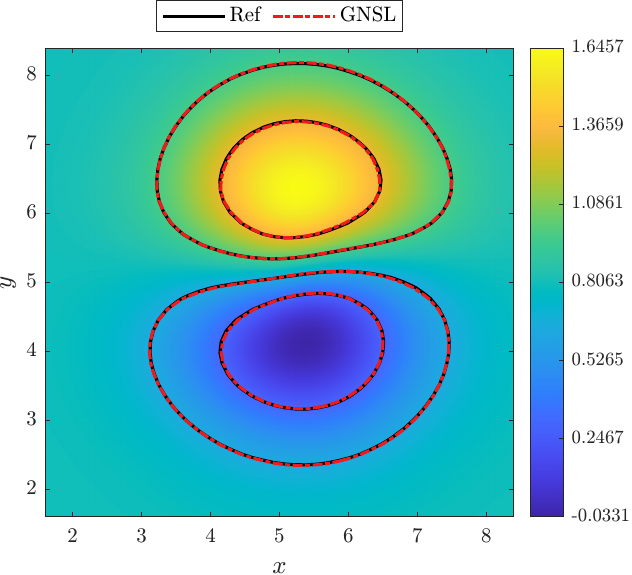}
        \caption{$u_2$}
    \end{subfigure} 
    \begin{subfigure}[b]{0.24\textwidth}
        \centering
        \includegraphics[width=\textwidth]{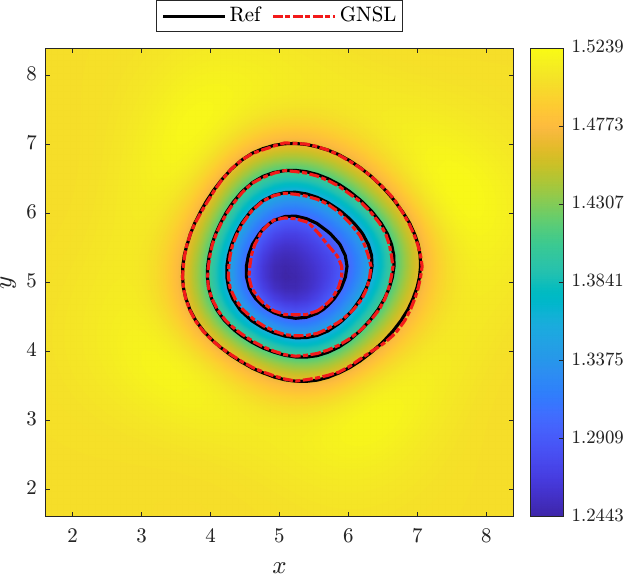}
        \caption{$\theta$}
    \end{subfigure}
    \hfill
    \caption{(2D3V isentropic vortex problem in Sec.~\ref{sec:ist})  The numerical solutions of the density $\rho$, macroscopic velocity $u_1,u_2$, and temperature $\theta$ obtained by GNSL and the Boltzmann equation for $\Kn = 0.2142$ at $t = 0.4$.}    
    \label{fig:vortex_example_4}
\end{figure}

To show the behavior of GNSL more precisely, four specific samples are selected, whose relative errors $\mE_{\rm GNSL}^t$ and corresponding Knudsen numbers are listed in Tab. \ref{tab:ist_error_1}. These samples are chosen so that they are all close to the average error levels at each time. The detailed numerical solutions for the density $\rho$, macroscopic velocity $u_1, u_2$, and the temperature $\theta$ are shown in Fig. \ref{fig:vortex_example_1} to \ref{fig:vortex_example_4}. It shows that for $t = 0.1, 0.2$ and $0.3$, the numerical solution and the reference solution match well with each other. For $t = 0.4$, there exist small differences between the numerical solution and the reference solution, which is consistent with the error listed in Tab. \ref{tab:ist_error_1}.

\begin{table}[htbp]
    \centering
    \def\arraystrech{1.5}
   {\footnotesize
    \begin{tabular}{crrrr}
        \toprule
         & $t=0.1$ & $t=0.2$ & $t=0.3$  & $t=0.4$ \\
        \midrule
        $\Kn$ & $0.5288$& $5.3690$ & $0.0913$ & $0.2142$\\
        \midrule
       $\mE_{\rm GNSL}^t(\%)$ & $0.05$ & $0.07$  & $0.11$ & $0.26$\\
        \bottomrule
    \end{tabular}}
     \caption{(2D3V isentropic vortex problem in Sec.~\ref{sec:ist}) The detailed information of the four selected samples from the testing data.}
    \label{tab:ist_error_1}
\end{table}

\subsection{2D3V Taylor-Green vortex problem}
\label{sec:TG}
In this section, the proposed neural network closure model is validated on the Taylor-Green vortex problem governed by the 2D3V Boltzmann equation.  The computational domain in the physical space is set to $[0,1]\times[0,1]$ with a periodic boundary condition. The initial condition is a Maxwellian function $\mathcal{M}[\bm{U}]$ with the macroscopic variables $\bm{U}=(\rho,u_1,u_2,\theta)$ given by
\begin{equation}
\label{eq:Taylor-Green}
\begin{gathered}
 \rho(x,y)=\rho_0, \qquad  \theta(x,y)=\theta_0-\frac{u_0^2}{4}(\cos(4\pi x)+\cos(4\pi y)), \\
 u_1(x,y)=-u_0\cos(2\pi x)\sin(2\pi y), \qquad 
 u_2(x,y)=u_0\sin(2\pi x)\cos(2\pi y), 
\end{gathered}
\end{equation}
where the parameters $\rho_0,\theta_0$ are uniformly sampled from $ [0.8,1.2]$ and $u_0$ is uniformly sampled from $[-0.5,0.5]$.

The same numerical method as in Sec.~\ref{sec:wave} is utilized to generate the training data. Precisely, the number of samples is $100$ with the spatial mesh size  $N_x = N_y = 48$ and the microscopic velocity mesh size $32\times 32\times 32$ in $[-7.5,7.5]^3$. The Knudsen number $\Kn$ is sampled from a log-uniform distribution over $[0.001, 10]$. The samples are generated from the initial condition \eqref{eq:Taylor-Green} and the sampled Knudsen numbers, and then evolved to the final simulation time $t = 0.1$ with the time step $\Delta t = 0.0005$ to obtain the training data. The training data set consists of the macroscopic variables at each time step. 

\begin{figure}[htbp]
    \centering
    \begin{subfigure}[b]{0.24\textwidth}
        \centering
        \includegraphics[width=\linewidth]{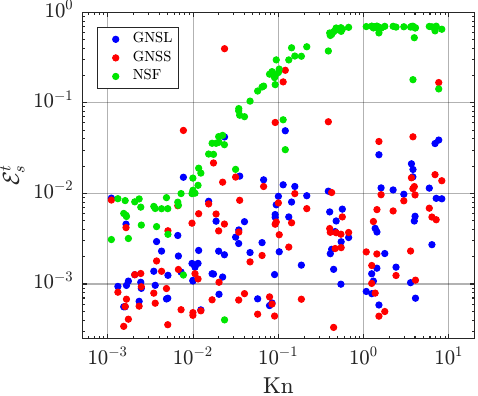}
        \caption{$t = 0.05$}
    \end{subfigure}
    \begin{subfigure}[b]{0.24\textwidth}
        \centering
        \includegraphics[width=\linewidth]{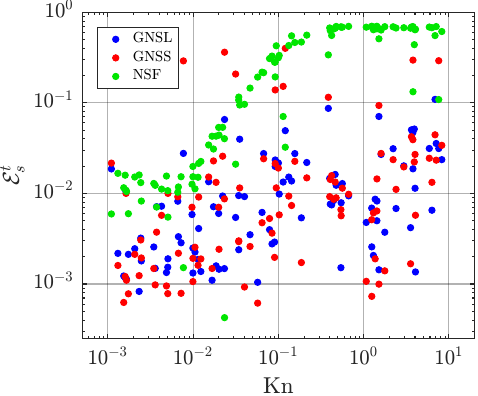}
        \caption{$t = 0.10$}
    \end{subfigure} 
    \begin{subfigure}[b]{0.24\textwidth}
        \centering
        \includegraphics[width=\linewidth]{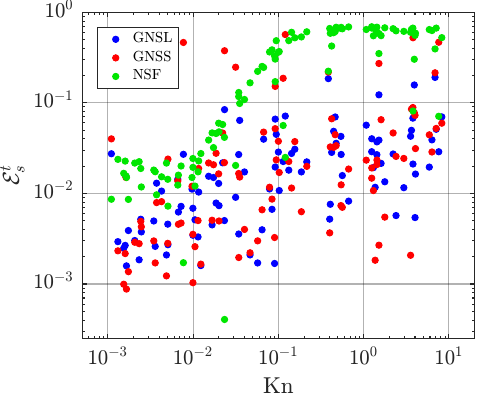}
        \caption{ $t = 0.15$}
    \end{subfigure}
    \begin{subfigure}[b]{0.24\textwidth}
        \centering
        \includegraphics[width=\linewidth]{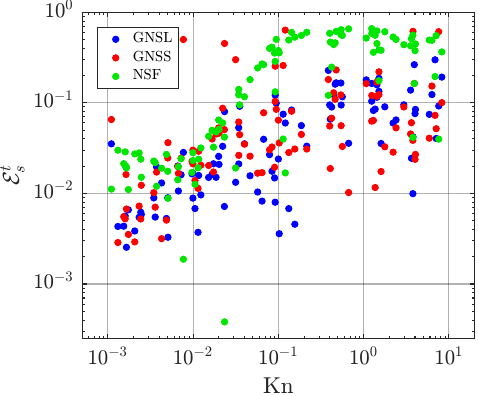}
        \caption{$t = 0.20$}
    \end{subfigure}    
    \caption{(2D3V Taylor-Green vortex problem in Sec.~\ref{sec:TG}) The distribution of the relative error $\mathcal{E}_{\rm s}^t$ with respect to the Knudsen number for NSF, GNSS, and GNSL at $t=0.05,0.10,0.15$ and $0.20$.}
    \label{fig:2D3VTG}
\end{figure}

For the testing data, the number of samples is also $100$, and they are generated following the same distribution as the training data. The final simulation time is extended to $t = 0.2$ to further test the generalization ability of this neural network closure model. The average relative errors $\mE_{\rm s}^t$ for NSF, GNSS and GNSL are listed in Tab.~\ref{tab:TG_error}, showing that for this 2D3V problem, the performance of GNSS and GNSL is also much better than that of NSF. The distribution of the relative error $\mE_{\rm s}^t$ with respect to the Knudsen number is shown in Fig.~\ref{fig:2D3VTG}. When the Knudsen number is small, the relative errors $\mE_{\rm s}^t$ for the three methods all keep a relatively low level. But that of NSF increases rapidly as $\Kn$ increases, while those of GNSS and GNSL remain small. These results indicate that the proposed neural network closure model successfully captures the complex dynamics of the decaying vortex. In contrast, the standard NSF equations produce substantial errors. Moreover, the subsequent training with $L_{\text{long}}$ provides a clear improvement in all test time steps.

\begin{table}[htbp]
    \centering
    \def\arraystrech{1.5}
    {\footnotesize
    \begin{tabular}{lrrrr}
        \toprule
        \diagbox[width=6em]{model}{time}
        & $t=0.05$ & $t=0.10$ & $t=0.15$ & $t=0.20$ \\
        \midrule
        NSF  & $28.80$ & $30.45$ & $29.21$ & $24.87$ \\
        \midrule
        GNSS & $1.60$ & $3.23$ & $5.48$ & $8.39$ \\
        \midrule
        GNSL & $\bm{0.60}$ & $\bm{1.40}$ & $\bm{2.64}$ & $\bm{5.83}$ \\
        \bottomrule
    \end{tabular}}
    \caption{(2D3V Taylor-Green vortex problem in Sec.~\ref{sec:TG}) The average relative errors $\mE_{\rm s}^{t}$ (in percentages) of the testing data for different methods at time $t=0.05$, $0.1$, $0.15$ and $0.2$.}
    \label{tab:TG_error}
\end{table}

\begin{figure}[!htbp]
    \centering
    \begin{subfigure}[b]{0.32\textwidth}
        \centering
        \includegraphics[width=\textwidth]
        {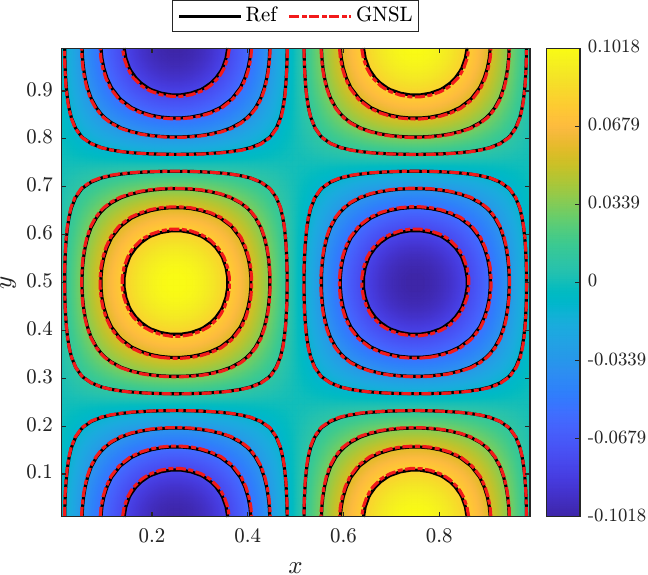}
        \caption{$u_1$}
    \end{subfigure}
    \hfill
    \begin{subfigure}[b]{0.32\textwidth}
        \centering
        \includegraphics[width=\textwidth]
        {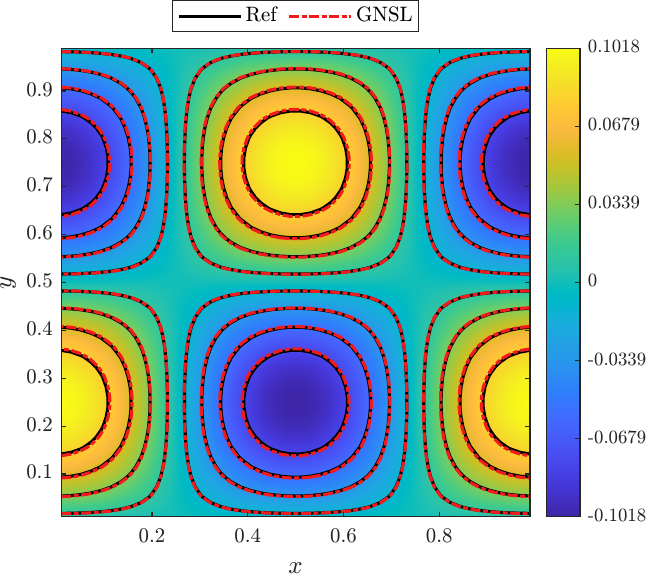}
        \caption{$u_2$}
    \end{subfigure}
    \hfill
    \begin{subfigure}[b]{0.32\textwidth}
        \centering
        \includegraphics[width=\textwidth]
        {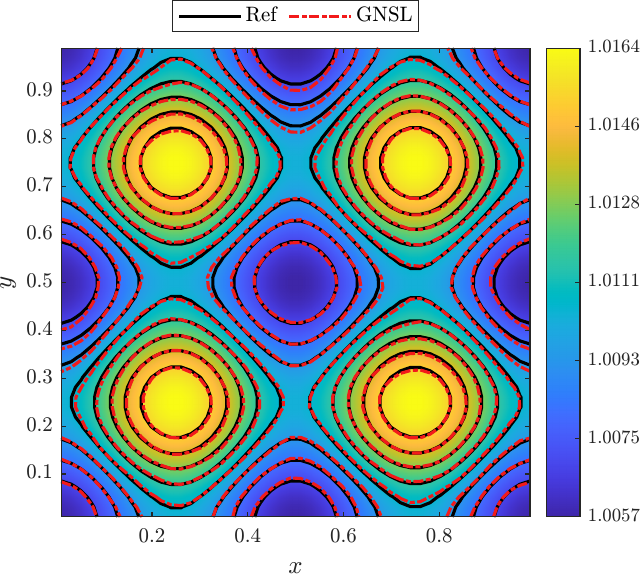}
        \caption{$\theta$}
    \end{subfigure}

    \caption{(2D3V Taylor-Green vortex problem in Sec.~\ref{sec:TG}) The numerical solutions of the macroscopic velocities $u_1$ and $u_2$ and temperature $\theta$ obtained by GNSL and the Boltzmann equation for $\Kn=0.0917$ at $t=0.05$.}
    \label{fig:TG_example_1}
\end{figure}

\begin{figure}[!htbp]
    \centering
    \begin{subfigure}[b]{0.32\textwidth}
        \centering
        \includegraphics[width=\textwidth]
        {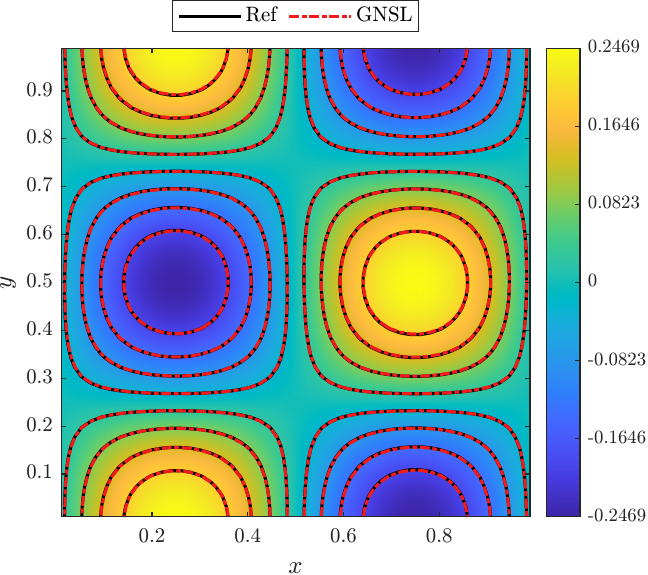}
        \caption{$u_1$}
    \end{subfigure}
    \hfill
    \begin{subfigure}[b]{0.32\textwidth}
        \centering
        \includegraphics[width=\textwidth]
        {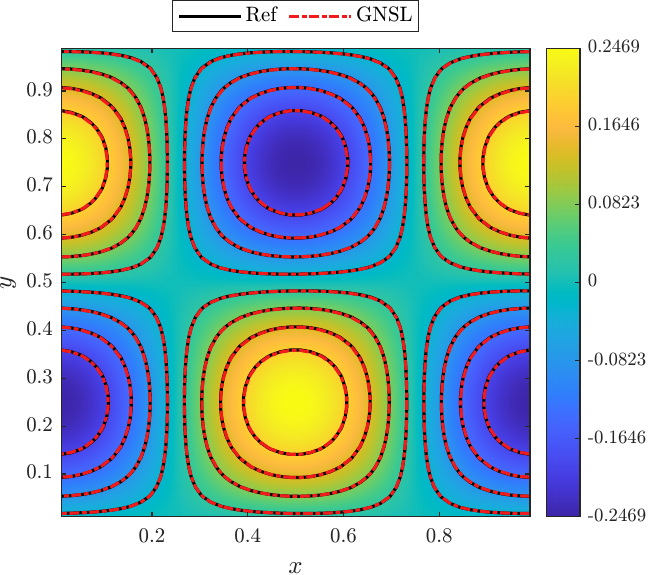}
        \caption{$u_2$}
    \end{subfigure}
    \hfill
    \begin{subfigure}[b]{0.32\textwidth}
        \centering
        \includegraphics[width=\textwidth]
        {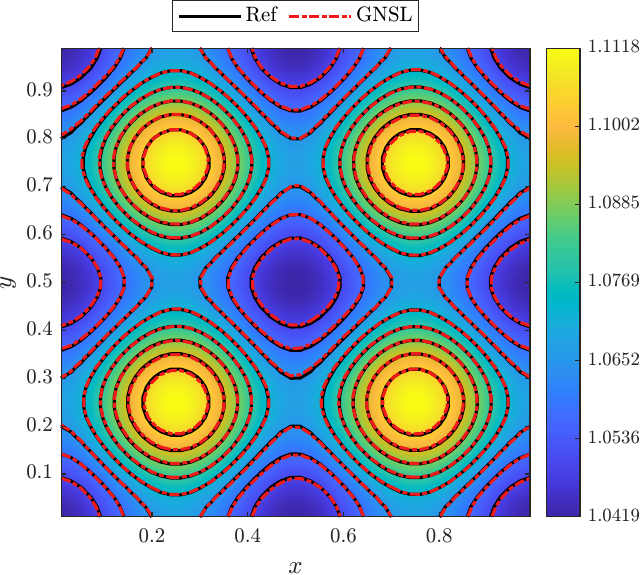}
        \caption{$\theta$}
    \end{subfigure}

    \caption{(2D3V Taylor-Green vortex problem in Sec.~\ref{sec:TG}) The numerical solutions of the macroscopic velocities $u_1$ and $u_2$ and temperature $\theta$ obtained by GNSL and the Boltzmann equation for $\Kn=0.1429$ at $t=0.10$.}
    \label{fig:TG_example_2}
\end{figure}

\begin{figure}[!htbp]
    \centering
    \begin{subfigure}[b]{0.32\textwidth}
        \centering
        \includegraphics[width=\textwidth]
        {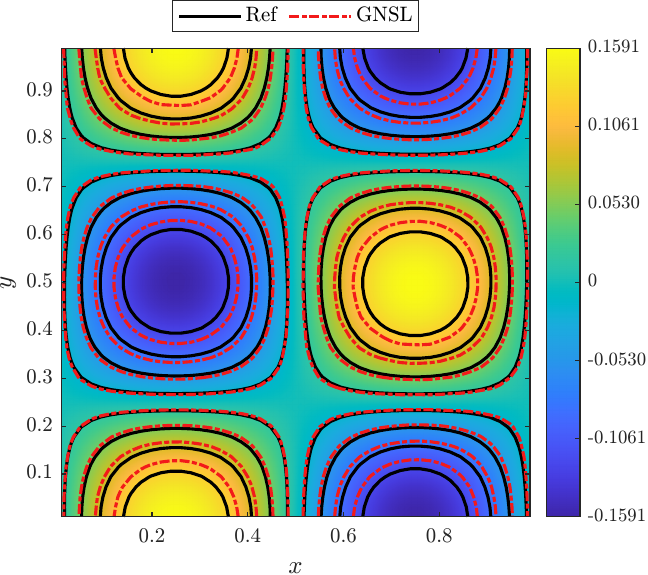}
        \caption{$u_1$}
    \end{subfigure}
    \hfill
    \begin{subfigure}[b]{0.32\textwidth}
        \centering
        \includegraphics[width=\textwidth]
        {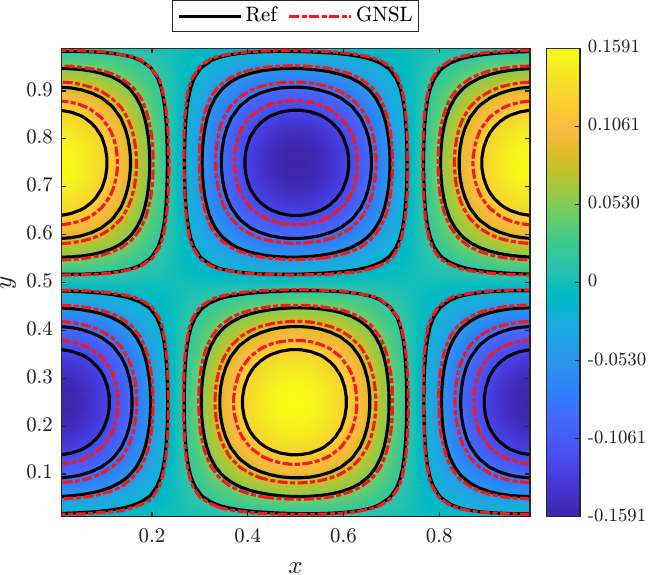}
        \caption{$u_2$}
    \end{subfigure}
    \hfill
    \begin{subfigure}[b]{0.32\textwidth}
        \centering
        \includegraphics[width=\textwidth]
        {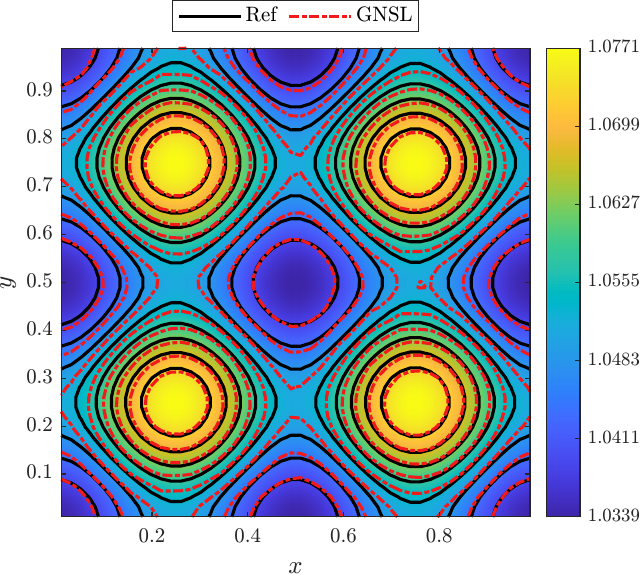}
        \caption{$\theta$}
    \end{subfigure}

    \caption{(2D3V Taylor-Green vortex problem in Sec.~\ref{sec:TG}) The numerical solutions of the macroscopic velocities $u_1$ and $u_2$ and temperature $\theta$ obtained by GNSL and the Boltzmann equation for $\Kn=1.4294$ at $t=0.15$.}
    \label{fig:TG_example_3}
\end{figure}

\begin{figure}[!htbp]
    \centering
    \begin{subfigure}[b]{0.32\textwidth}
        \centering
        \includegraphics[width=\textwidth]
        {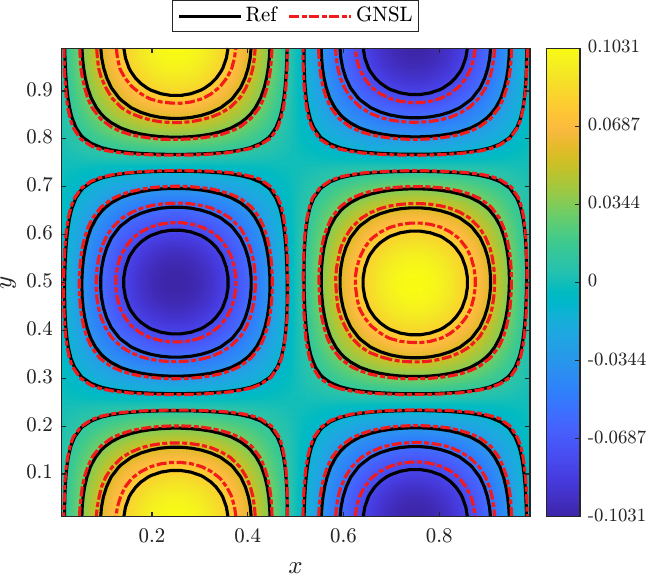}
        \caption{$u_1$}
    \end{subfigure}
    \hfill
    \begin{subfigure}[b]{0.32\textwidth}
        \centering
        \includegraphics[width=\textwidth]
        {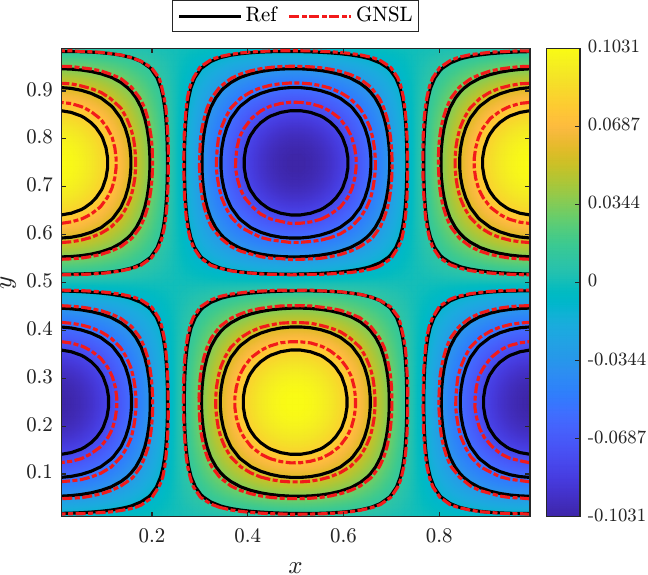}
        \caption{$u_2$}
    \end{subfigure}
    \hfill
    \begin{subfigure}[b]{0.32\textwidth}
        \centering
        \includegraphics[width=\textwidth]
        {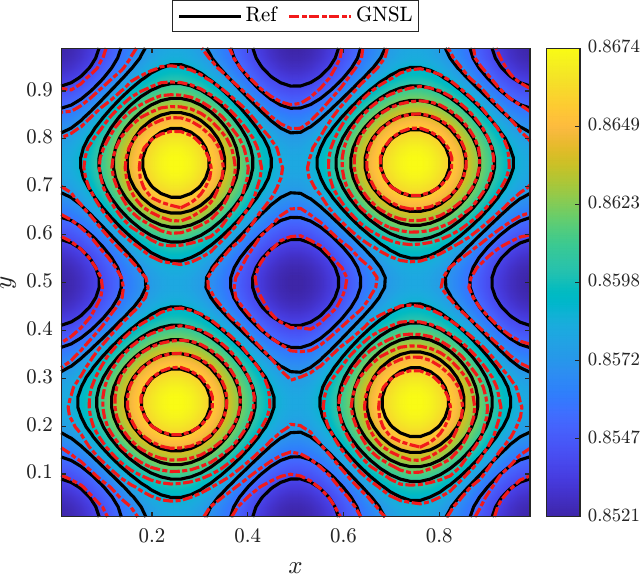}
        \caption{$\theta$}
    \end{subfigure}
    \caption{(2D3V Taylor-Green vortex problem in Sec.~\ref{sec:TG}) The numerical solutions of the macroscopic velocities $u_1$ and $u_2$ and temperature $\theta$ obtained by GNSL and the Boltzmann equation for $\Kn=0.1861$ at $t=0.20$.}
    \label{fig:TG_example_4}
\end{figure}

To show the behavior of GNSL more precisely, four specific samples are selected, whose relative errors $\mE_{\rm GNSL}^t$ and corresponding Knudsen numbers are listed in Tab. \ref{tab:TG_error_1}. The numerical solutions of the macroscopic velocity $u_1, u_2$, and the temperature $\theta$ are shown in Fig. \ref{fig:TG_example_1}, \ref{fig:TG_example_2}, \ref{fig:TG_example_3} and \ref{fig:TG_example_4}. They show that when $t = 0.05$ and $0.1$, the numerical solution of GNSL and the reference solution match well with each other, while when $t = 0.15$ and $0.2$, there are clear differences, especially when $t = 0.2$. This is also consistent with the numerical error listed in Tab. \ref{tab:TG_error_1}, indicating that the generalization ability of GNSL needs to be improved, especially when capturing the structure of the complex vortex.

\begin{table}[htbp]
    \centering
    \def\arraystrech{1.5}
   {\footnotesize
    \begin{tabular}{crrrr}
        \toprule
         & $t=0.05$ & $t=0.10$ & $t=0.15$  & $t=0.20$ \\
        \midrule
        $\Kn$ & $0.0917$& $0.1429$ & $1.4294$ & $0.1861$\\
        \midrule
       $\mE_{\rm GNSL}^t(\%)$ & $0.58$ & $1.37$  & $2.68$ & $5.61$\\
        \bottomrule
    \end{tabular}}
     \caption{(2D3V Taylor-Green vortex problem in Sec.~\ref{sec:TG}) The detailed information of the four selected samples from the testing data.}
    \label{tab:TG_error_1}
\end{table}

\paragraph{Studying the effect of the trajectory length}
For the last part of the numerical test, the effect of the short-term trajectory length $K_s$ is studied with the example Taylor-Green vortex problem, while $K_l=2$ is fixed for the subsequent $L_{\text{long}}$ training. The setting for the training and testing data is the same as in Sec. \ref{sec:TG}. The average relative errors for GNSS and GNSL with $K_s = 4$ and $K_s = 8$ are listed in Tab. \ref{tab:TG_error_2}. The comparison reveals that the relative error $\mE_{\rm GNSL}^t$ decreases as the trajectory length $K_s$ increases during the initial $L_{\text{short}}$ training. However, this accuracy gain is relatively minor compared to the improvement achieved by the subsequent $L_{\text{long}}$ training. Notably, combining both strategies by conducting the $L_{\text{long}}$ training with $K_s=8$ yields the highest accuracy among all cases. This confirms that optimizing over longer trajectories improves the long-term accuracy of the model. However, increasing $K_s$ will make the training process more expensive, and the balance between the trajectory length $K_s$ and the relative error will be left to future work. 

\begin{table}[htbp]
    \centering   
    \begin{tabular}{lcccc}
        \toprule
          & $t=0.05$ & $t=0.10$ & $t=0.15$ & $t=0.20$ \\
        \midrule
        GNSS, $K_s=4$ & $1.60$ & $3.23$ & $5.48$ & $8.39$ \\
        \midrule
        GNSS, $K_s=8$ & $0.98$ & $1.90$ & $3.31$ & $7.02$ \\
        \midrule
        GNSL, $K_s=4$ & $0.60$ & $1.40$ & $2.64$ & $5.83$ \\
        \midrule
        GNSL, $K_s=8$ & $\bm{0.40}$ & $\bm{0.96}$ & $\bm{2.17}$ & $\bm{5.27}$   \\
        \bottomrule
    \end{tabular}
     \caption{(2D3V Taylor-Green vortex problem in Sec.~\ref{sec:TG}) The average relative errors (in percentages) of the testing data for different methods at time $t=0.05, 0.1, 0.15$ and $0.2$.}
    \label{tab:TG_error_2}
\end{table}

\subsection{Study of the computational efficiency}
\label{sec:computational_efficiency}
The computational efficiency for the GNS model is studied in this section, where the 1D3V wave problem in Sec. \ref{sec:wave} and the 2D3V Taylor-Green vortex problem in Sec. \ref{sec:TG} are studied. Here, the time to numerically solve the GNS is recorded, and the time to train the neural network is excluded. For the Boltzmann equation, it is simulated using a CPU with $36$ cores and a parallel implementation. The GNS model is simulated on a single GPU and NSF on a single CPU. For all the numerical results, the relative error and the simulation time are averaged over $10$ independent runs.
\begin{table}[htbp]
    \centering
    \def\arraystrech{1.35}
    {\footnotesize
    \resizebox{\linewidth}{!}{%
    \begin{tabular}{rccrrrrrr}
    \toprule
    & & \multicolumn{1}{c}{Boltzmann} & \multicolumn{2}{c}{NSF} & \multicolumn{2}{c}{GNSS} & \multicolumn{2}{c}{GNSL} \\
    \cmidrule(lr){3-3}\cmidrule(lr){4-5}\cmidrule(lr){6-7}\cmidrule(lr){8-9}$\Kn$ & $t$ & time (s) & error (\%)  & time (s) & error (\%)  & time (s) & error (\%) & time (s) \\
    \midrule
    \multirow{2}{*}{$0.001$} & $0.10$ & $364.84$ & $3.12$ & $0.12$ & $4.74$ & $1.97$ & $3.59$ & $1.97$ \\ 
     & $0.20$ & $729.43$ & $2.60$ & $0.25$ & $6.29$ & $3.96$ & $4.45$ & $3.97$ \\
    \midrule
    \multirow{2}{*}{$0.01$} & $0.10$ & $365.26$ & $8.08$ & $0.13$ & $2.27$ & $1.97$ & $1.88$ & $1.98$ \\ 
     & $0.20$ & $730.47$ & $10.13$ & $0.25$ & $2.51$ & $3.96$ & $1.82$ & $3.97$ \\
    \midrule
    \multirow{2}{*}{$0.1$} & $0.10$ & $364.71$ & $30.97$ & $0.13$ & $3.53$ & $1.97$ & $2.63$ & $1.97$ \\ 
     & $0.20$ & $730.62$ & $27.90$ & $0.25$ & $6.73$ & $3.97$ & $2.84$ & $3.97$ \\
    \midrule
    \multirow{2}{*}{$1.0$} & $0.10$ & $365.35$ & $39.47$ & $0.12$ & $1.75$ & $1.98$ & $2.02$ & $1.97$ \\
     & $0.20$ & $730.84$ & $22.59$ & $0.25$ & $9.85$ & $3.97$ & $3.39$ & $3.96$ \\
    \midrule
    \multirow{2}{*}{$10.0$} & $0.10$ & $365.63$ & $\mathbf{39.74}$ & $0.12$ & $1.46$ & $1.97$ & $\mathbf{1.03}$ & $1.98$ \\
     & $0.20$ & $\mathbf{731.36}$ & $31.39$ & $0.25$ & $10.77$ & $3.97$ & $6.21$ & $\mathbf{3.97}$ \\
    \bottomrule
    \end{tabular}}}
    \caption{(Study of the computational efficiency in Sec. \ref{sec:computational_efficiency}) The efficiency study of GNS for the 1D3V wave problem.}
    \label{tab:wave_efficiency}
\end{table}

To compare the computational efficiency in different kinetic regimes, the Knudsen numbers $\Kn = 0.001, 0.01, 0.1, 1.0$ and $10$ are tested. The relative errors obtained by \eqref{eq:relative_error} and computational costs at $t = 0.1$ and $0.2$ are shown in Tab. \ref{tab:wave_efficiency}. In Sec. \ref{sec:wave}, the final time for the training data is $t = 0.1$. Thus, the numerical results at $t = 0.1$ show the modeling ability of GNS, and those at $t = 0.2$ indicate the generalization ability of GNS. Tab. \ref{tab:wave_efficiency} shows that the computational time for the Boltzmann equation is much larger for all $\Kn$ compared to other models. When the Knudsen number is small, as $\Kn = 0.001$, the relative error for NSF is a little smaller than that of GNS, with the computational time a little less than that of GNS. However, with increasing $\Kn$, the relative error for NSF grows fast, while that of GNS model still remains small. Especially, when $\Kn = 10$, the relative error of NSF at $t = 0.1$ is almost $20$ times larger than that of GNS, and that at $t = 0.2$ is $2$ times larger than that of GNSS, and $4$ times larger than that of GNSL.

The numerical results for the 2D3V Taylor-Green vortex problem are listed in Tab. \ref{tab:TG_efficiency}. The behavior of GNS is very similar to that of the 1D3V problems. When $\Kn$ is large, although the computational time for GNS is a little larger than that of NSF, the relative error is much smaller, indicating the validity of GNS in the transitional regime. Moreover, the relative error of GNSL is smaller than that of GNSS for all Knudsen numbers, showing the stronger robustness of GNSL compared to GNSS. Besides, we can also expect that GNS model will be faster than that of NSF for the 3D3V problems.




\begin{table}[htbp]
    \centering
    \def\arraystrech{1.35}
    {\footnotesize
    \resizebox{\linewidth}{!}{%
    \begin{tabular}{rcrrrrrrr}
        \toprule
        & & \multicolumn{1}{c}{Boltzmann} & \multicolumn{2}{c}{NSF} & \multicolumn{2}{c}{GNSS} & \multicolumn{2}{c}{GNSL} \\
        \cmidrule(lr){3-3}\cmidrule(lr){4-5}\cmidrule(lr){6-7}\cmidrule(lr){8-9} $\Kn$ & $t$ & time (s) & error (\%) & time (s) & error (\%) & time (s) & error (\%) & time (s) \\
        \midrule
       \multirow{2}{*}{$0.001$}  & $0.10$ & $776.84$ & $0.46$ & $1.09$ & $1.36$ & $4.57$ & $0.44$ & $4.58$ \\ 
         & $0.20$ & $1557.93$ & $0.83$ & $2.18$ & $2.83$ & $9.15$ & $0.90$ & $9.17$ \\ 
        \midrule
       \multirow{2}{*}{$0.01$} & $0.10$ & $781.25$ & $1.43$ & $1.09$ & $0.26$ & $4.60$ & $0.19$ & $4.57$ \\ 
         & $0.20$ & $1562.34$ & $2.11$ & $2.19$ & $2.57$ & $9.20$ & $0.54$ & $9.16$ \\ 
        \midrule
       \multirow{2}{*}{$0.1$} & $0.10$ & $778.63$ & $29.19$ & $1.10$ & $2.16$ & $4.58$ & $1.32$ & $4.60$ \\ 
       & $0.20$ & $1559.46$ & $30.06$ & $2.18$ & $2.30$ & $9.16$ & $0.61$ & $9.20$ \\ 
        \midrule
       \multirow{2}{*}{$1.0$} & $0.10$ & $783.17$ & $62.73$ & $1.09$ & $1.89$ & $4.60$ & $1.72$ & $4.59$ \\ 
     & $0.20$ & $1564.52$ & $33.55$ & $2.19$ & $10.78$ & $9.20$ & $7.71$ & $9.18$ \\
        \midrule
        \multirow{2}{*}{$10.0$} & $0.10$ & $780.42$ & $\mathbf{46.34}$ & $1.10$ & $2.04$ & $4.57$ & $\mathbf{1.07}$ & $4.60$ \\ 
         & $0.20$ & $\mathbf{1561.87}$ & $20.42$ & $2.19$ & $10.22$ & $9.15$ & $4.15$ & $\mathbf{9.20}$ \\
        \bottomrule
    \end{tabular}}}
    \caption{(Study of the computational efficiency in Sec. \ref{sec:computational_efficiency}) The efficiency study of GNS for the 2D3V Taylor-Green vortex problem.}
    \label{tab:TG_efficiency}
\end{table}
\section{Conclusion}
\label{sec:con}


In this work, a generalized Navier-Stokes (GNS) model is constructed to enlarge the application scope of the classical NS equations, where a data-driven method is proposed to learn a new relationship between the non-equilibrium variables, such as the stress tensor and heat flux in the NS equations, and the equilibrium macroscopic variables such as the density, velocity, and temperature. The neural network is specially designed to preserve Galilean invariance, while a U-Net architecture is employed to capture nonlocal spatial features. In addition, a long-term end-to-end learning strategy is developed to reduce accumulated errors and improve long-term predictive accuracy.

Numerical experiments on both 1D3V and 2D3V problems are performed to validate the efficiency of this new GNS model. The timing comparisons further show that the GNS model substantially reduces computational cost compared with solving the high-dimensional Boltzmann equation. Using alternative sources of training data, including experimental data and numerical data generated by other kinetic models, and the applications to higher-dimensional problems will be left for future work.

\section*{Acknowledgments}
This work is partially supported by the Science Challenge program under the Grant No. TZ2026004 (Y. W.), the National Natural Science Foundation of China under the Grant No. 12325112 (S. S., Z. X.), and the High-performance Computing Platform of Peking University.

\FloatBarrier 
\bibliographystyle{plain}
\bibliography{references.bib}
\end{document}